\documentclass[final,3p,fleqn,11pt]{elsarticle}

\usepackage{amssymb}
\usepackage{multicol}
\usepackage[svgnames]{xcolor}

\usepackage{booktabs}
\usepackage{array}
\newcommand{\PreserveBackslash}[1]{\let\temp=\\#1\let\\=\temp}
\newcolumntype{C}[1]{>{\PreserveBackslash\centering}p{#1}}
\newcolumntype{R}[1]{>{\PreserveBackslash\raggedleft}p{#1}}
\newcolumntype{L}[1]{>{\PreserveBackslash\raggedright}p{#1}}

\journal{Computer Methods in Applied Mechanics and Engineering (CMAME)}

\usepackage{amsmath}
\usepackage{amsfonts}
\usepackage{amsthm}
\usepackage{mathtools}
\usepackage{subfig}
\usepackage{lineno}
\usepackage[hidelinks]{hyperref}
\usepackage{textcomp}
\usepackage{algorithm,algpseudocode,algorithmicx}

\usepackage{tabularx}

\usepackage{physics}
\usepackage{MnSymbol}

\usepackage{float}
\usepackage{array,multirow}
\usepackage{color}
\usepackage{tikz}
\usetikzlibrary{patterns}
\usetikzlibrary{spy}
\usetikzlibrary{spy,decorations.fractals,shapes.misc}

\usepackage{pgfplots}
\pgfplotsset{compat=1.16}
\newsavebox\Axis

\definecolor{lightgray}{gray}{0.80}

\usetikzlibrary{decorations.pathreplacing}
\usepackage[listings,theorems,skins,breakable]{tcolorbox}
\tcbset{skin=enhanced}
\newtcolorbox{lbracebox}[1][Word]{%
   frame hidden,enlarge left by=2cm,width=\linewidth-2cm,%
  overlay unbroken = {\draw [decorate,decoration={brace,amplitude=10pt},]%
                     (frame.south west)-- (frame.north west)
                    node [black,midway,left,xshift=-.6cm] {#1};},%
}

\newcommand{\Bezier}{B\'ezier~}

\definecolor{grey1}{rgb}{0.5, 0.5, 0.5}
\definecolor{green1}{rgb}{0.4660, 0.6740, 0.1880} 
\definecolor{blue1}{rgb}{0, 0.4470, 0.7410} 
\definecolor{red1}{rgb}{0.8500, 0.3250, 0.0980}
\definecolor{yellow1}{rgb}{0.9290, 0.6940, 0.1250}
\definecolor{purple1}{rgb}{0.4940, 0.1840, 0.5560}
\definecolor{lightblue1}{rgb}{0.3010, 0.7450, 0.9330}
\definecolor{bordeaux1}{rgb}{0.6350, 0.0780, 0.1840}

\definecolor{burntorange}{rgb}{0.74902,0.341176,0}

\newcommand{\changed}[1]{\textcolor{%
black}{#1}}
\newcommand{\changedV}[1]{\textcolor{
black}{#1}}

\makeatletter
\renewcommand*\env@matrix[1][*\c@MaxMatrixCols c]{%
  \hskip -\arraycolsep
  \let\@ifnextchar\new@ifnextchar
  \array{#1}}
\makeatother

\theoremstyle{plain}
\newtheorem{theorem}{Theorem}[section]

\newtheorem{remark}[theorem]{Remark}

\theoremstyle{definition}

\newcommand{\vect}[1]{\boldsymbol{#1}} 									
\newcommand{\mat}[1]{\mathbf{#1}} 											
\newcommand{\eigenvec}{U}             
\newcommand{\laplace}{\Delta}         
\newcommand{\npa}{n_{\text{pa}}}      
\newcommand{\nele}{n_{\text{ele}}}    
\newcommand{\noMode}{n}               
\newcommand{\normalvect}{\boldsymbol{\nu}}    
\newcommand{\npai}{\mathcal{E}}       

\DeclareMathOperator*{\argmax}{arg\,max}

\begin{document}

\begin{frontmatter}

\title{A variational approach based on perturbed eigenvalue analysis for improving spectral properties of isogeometric multipatch discretizations}

\corref{cor1}
\author[address1]{Thi-Hoa Nguyen}
\ead{nguyen@mechanik.tu-darmstadt.de}

\author[address1]{Ren\'e R. Hiemstra}
\ead{hiemstra@mechanik.tu-darmstadt.de}

\author[address2]{Stein K. F. Stoter}
\ead{k.f.s.stoter@tue.nl}

\author[address1]{Dominik Schillinger}
\ead{schillinger@mechanik.tu-darmstadt.de}

\cortext[cor1]{Corresponding author}

\address[address1]{Institute of Mechanics, Computational Mechanics Group, Technical University of Darmstadt, Germany}

\address[address2]{Department of Mechanical Engineering, Eindhoven University of Technology, The Netherlands}

\begin{abstract}
A key advantage of isogeometric discretizations is their accurate and well-behaved eigenfrequencies and eigenmodes. For degree two and higher, however, optical branches of spurious outlier frequencies and modes may appear due to boundaries or reduced continuity at patch interfaces.
In this paper, we introduce a variational approach based on perturbed eigenvalue analysis that eliminates outlier frequencies without negatively affecting the accuracy in the remainder of the spectrum and modes.
We then propose a pragmatic iterative procedure that estimates the  
perturbation parameters in such a way that the outlier frequencies are effectively reduced.
We demonstrate that our approach allows for a much larger critical time-step size in explicit dynamics calculations. In addition, we show that the critical time-step size obtained with the proposed approach does not depend on the polynomial degree of spline basis functions. 
\end{abstract}

\begin{highlights}
	\item We improve spectral properties of isogeometric multipatch discretizations via a variational approach based on perturbed eigenvalue analysis.
	\item We present an iterative procedure to estimate optimal scaling parameters such that the outlier frequencies are effectively reduced. 
	\item Our approach allows for a much larger critical time-step size in explicit dynamics, independently of the polynomial degree.
\end{highlights}


\end{frontmatter}


\section{Introduction}\label{sec:intro}

\begin{figure}[h!]
    \centering
    \subfloat[Normalized frequencies]{{\includegraphics[width=0.5\textwidth]{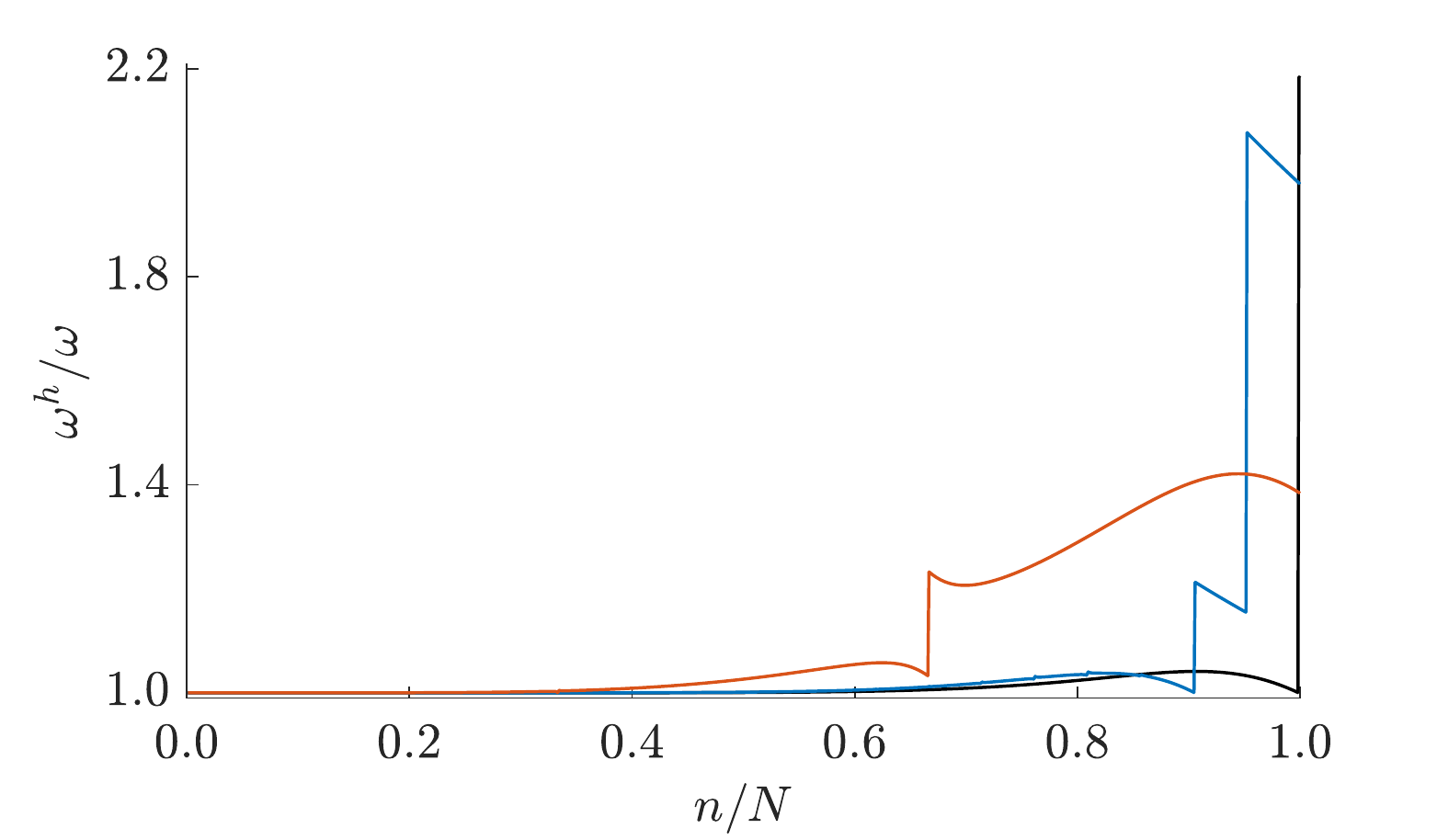} }}
    \subfloat[$L^2$ errors in the mode shapes]{{\includegraphics[width=0.5\textwidth]{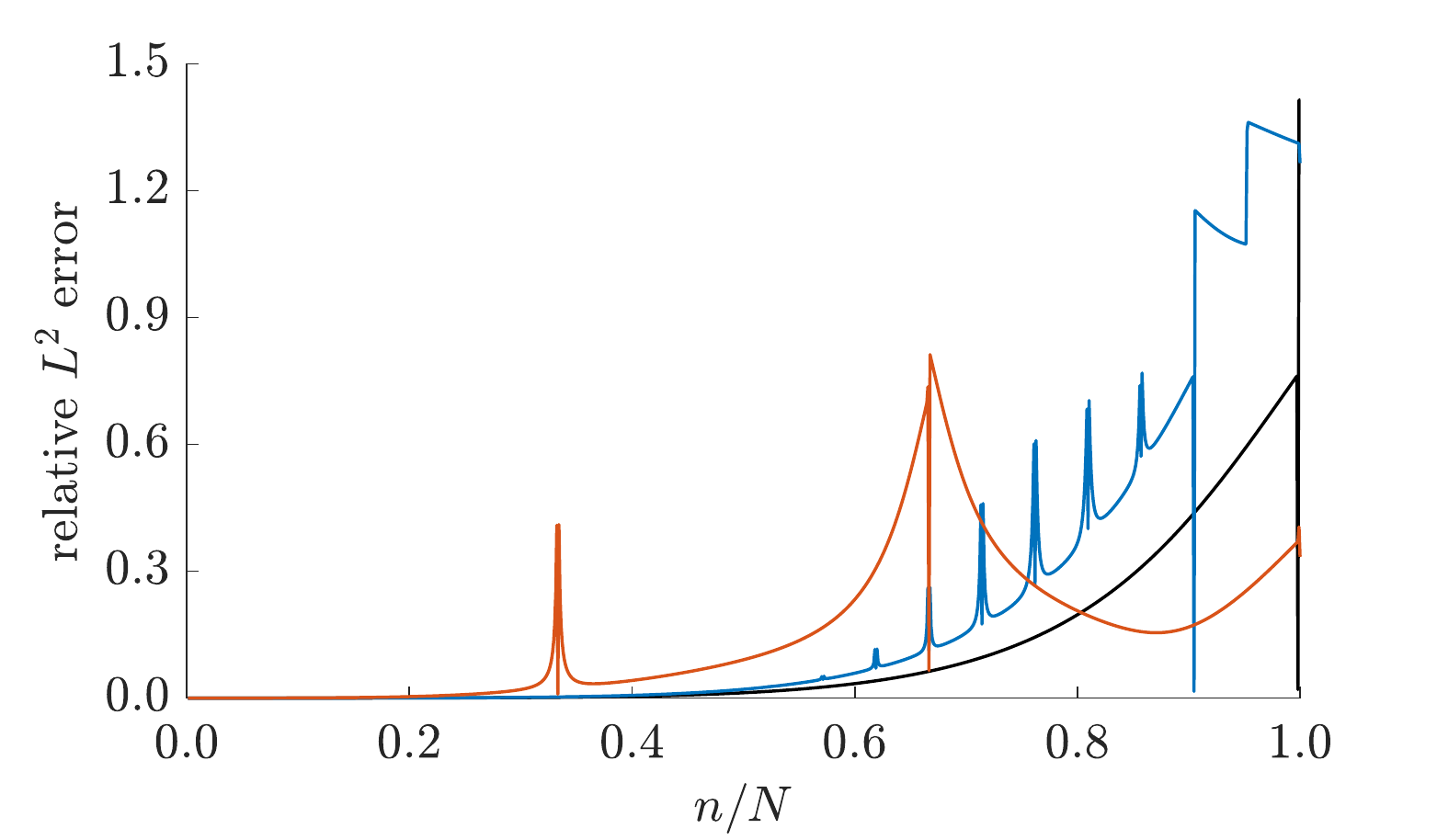} }}
    \vspace{0.3cm}
    \begin{tikzpicture}
        \filldraw[black,line width=1pt, solid] (0.0,0) -- (0.5,0);
		\filldraw[black,line width=1pt] (0.5,0) node[right]{\footnotesize $1$ patch};
		\filldraw[blue1,line width=1pt, solid] (3.0,0) -- (3.5,0);
		\filldraw[blue1,line width=1pt] (3.5,0) node[right]{\footnotesize $50$ patches};
		\filldraw[red1,line width=1pt, solid] (6.0,0) -- (6.5,0);
		\filldraw[red1,line width=1pt] (6.5,0) node[right]{\footnotesize $350 \;(\nele)$ patches};				
	\end{tikzpicture}
    \caption{Normalized frequencies and $L^2$ errors in the mode shapes of a \changed{free vibrating} \textbf{bar with free boundary conditions}, unit length and unit material parameters, computed with \textbf{cubic $C^2$ B-splines ($p=3$), $N=1050$ modes} and \textbf{increasing number of patches with $C^0$ continuity at patch interfaces}. The rigid body mode is excluded.}
    \label{fig:spectra_1Dbar_increasing_patches}
\end{figure}
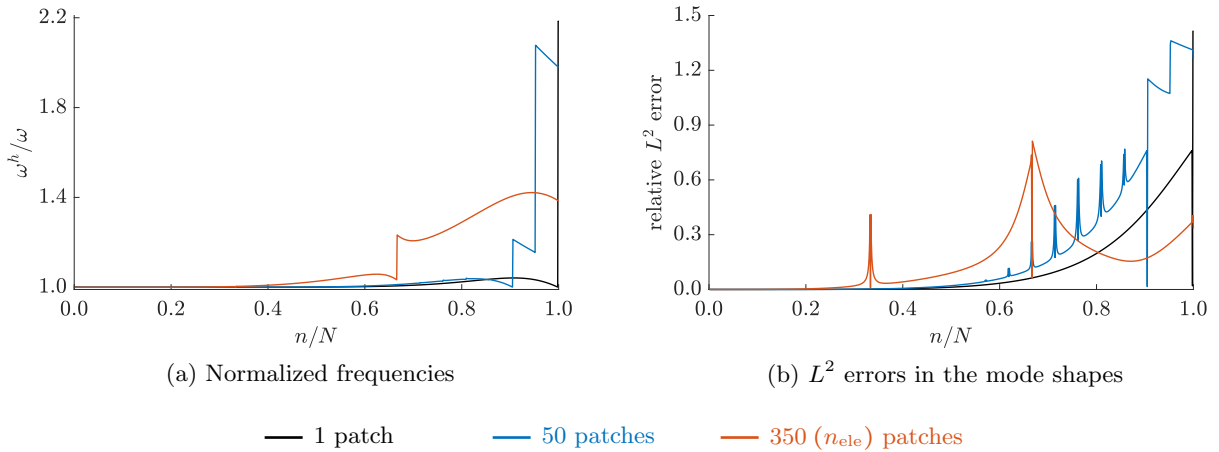

\begin{figure}[h!]
    \centering
    \subfloat[mode $N-1$]{{\includegraphics[width=0.43\textwidth]{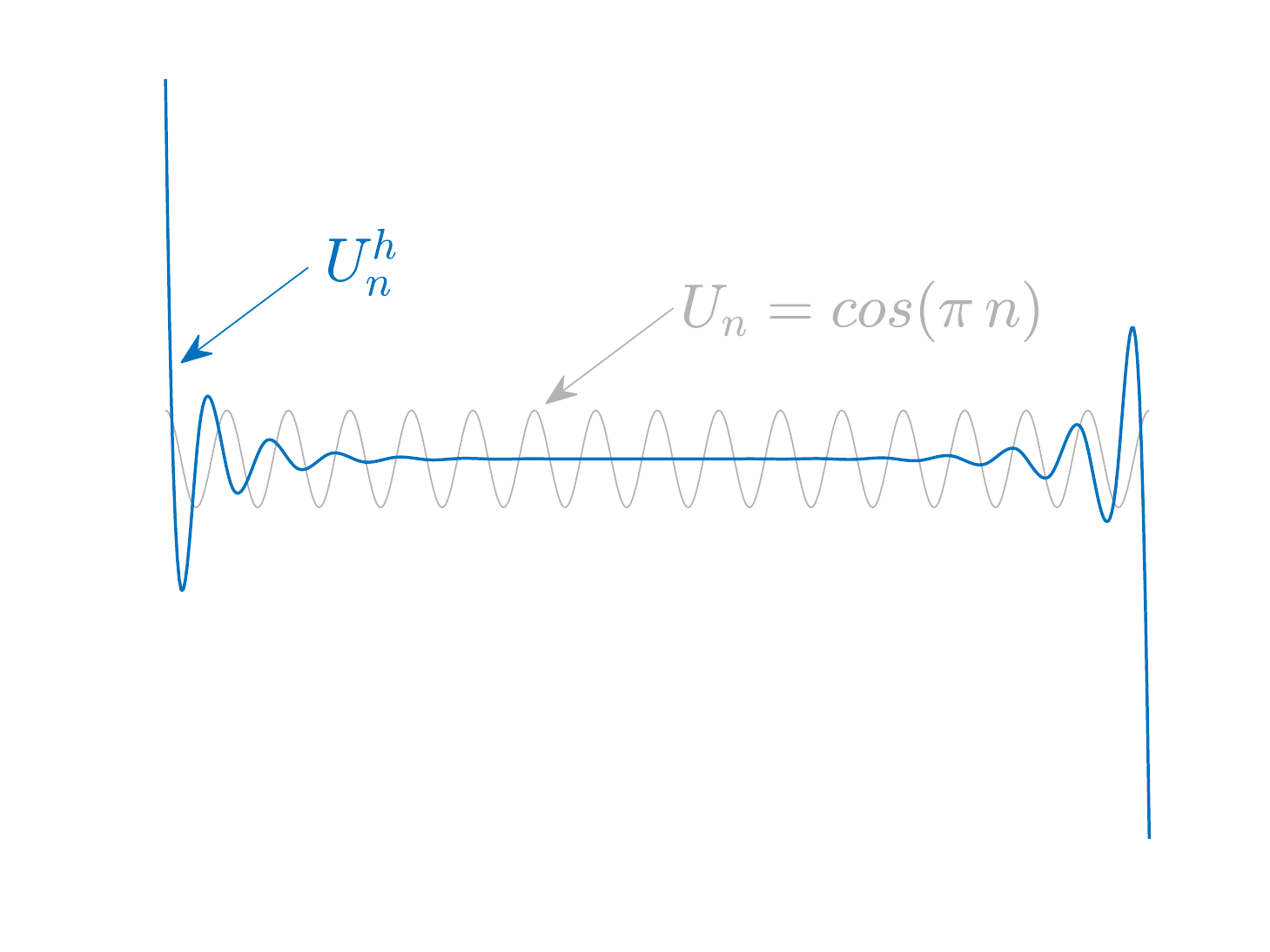} }}
    \subfloat[mode $N$]{{\includegraphics[width=0.43\textwidth]{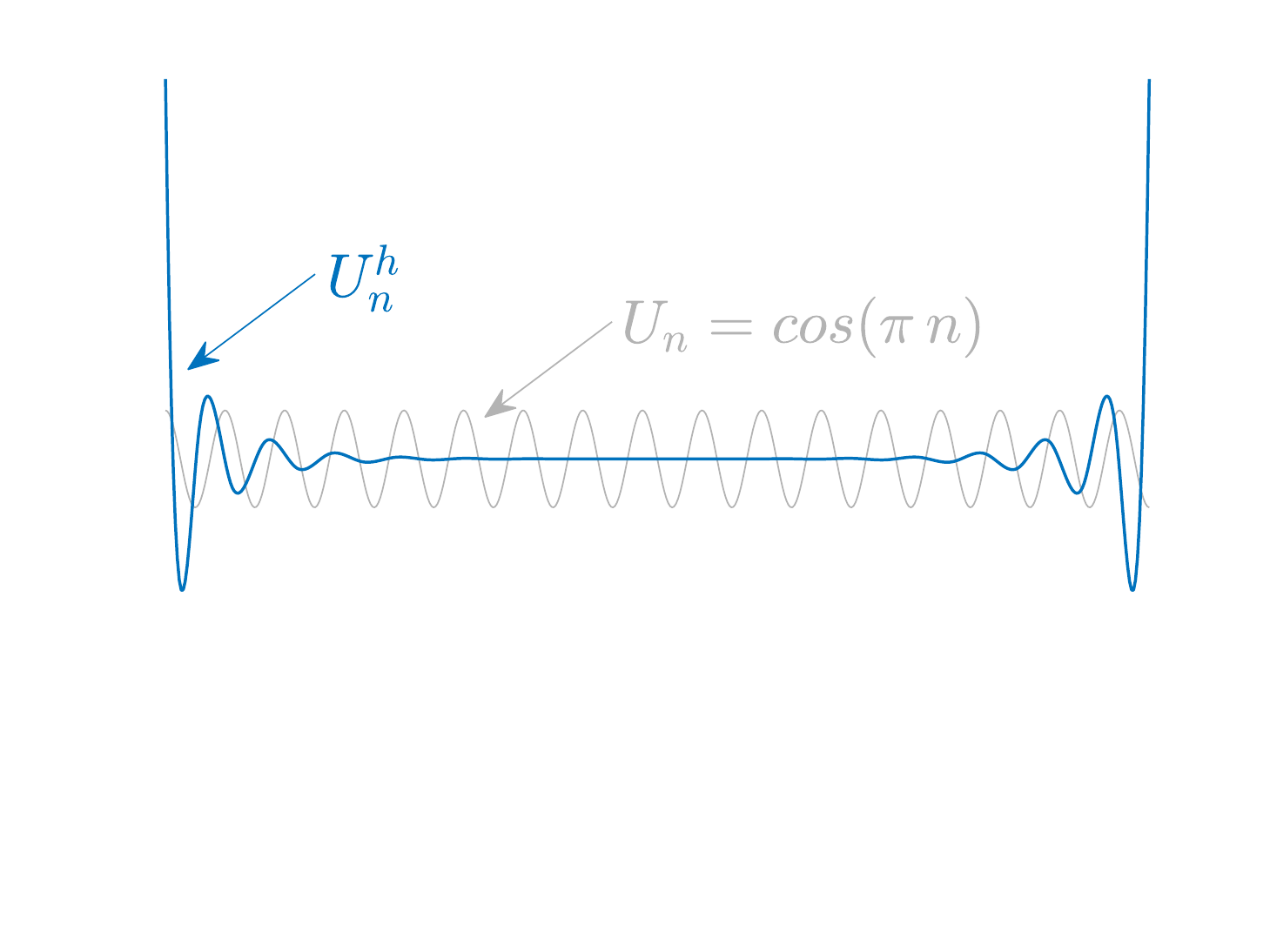} }}
    \caption{Discrete outlier modes $\eigenvec_n^h$ (\textbf{blue}) corresponding to the example of \textbf{the free bar studied in Figure \ref{fig:spectra_1Dbar_increasing_patches}}, computed with \textbf{one patch of cubic $C^2$ B-splines ($p=3$)} and discretized with $30$ elements. These modes behave significantly differently \changed{compared to} the exact modes (\textbf{gray}) $\eigenvec_n = \cos (\pi \, \noMode), \noMode = 1, \ldots, N$, where $N=33$ is the number of modes.}
    \label{fig:outlier_mode_shapes}
\end{figure}

\begin{figure}[h!]
    \centering
    \subfloat[mode $N-3$]{{\includegraphics[width=0.43\textwidth]{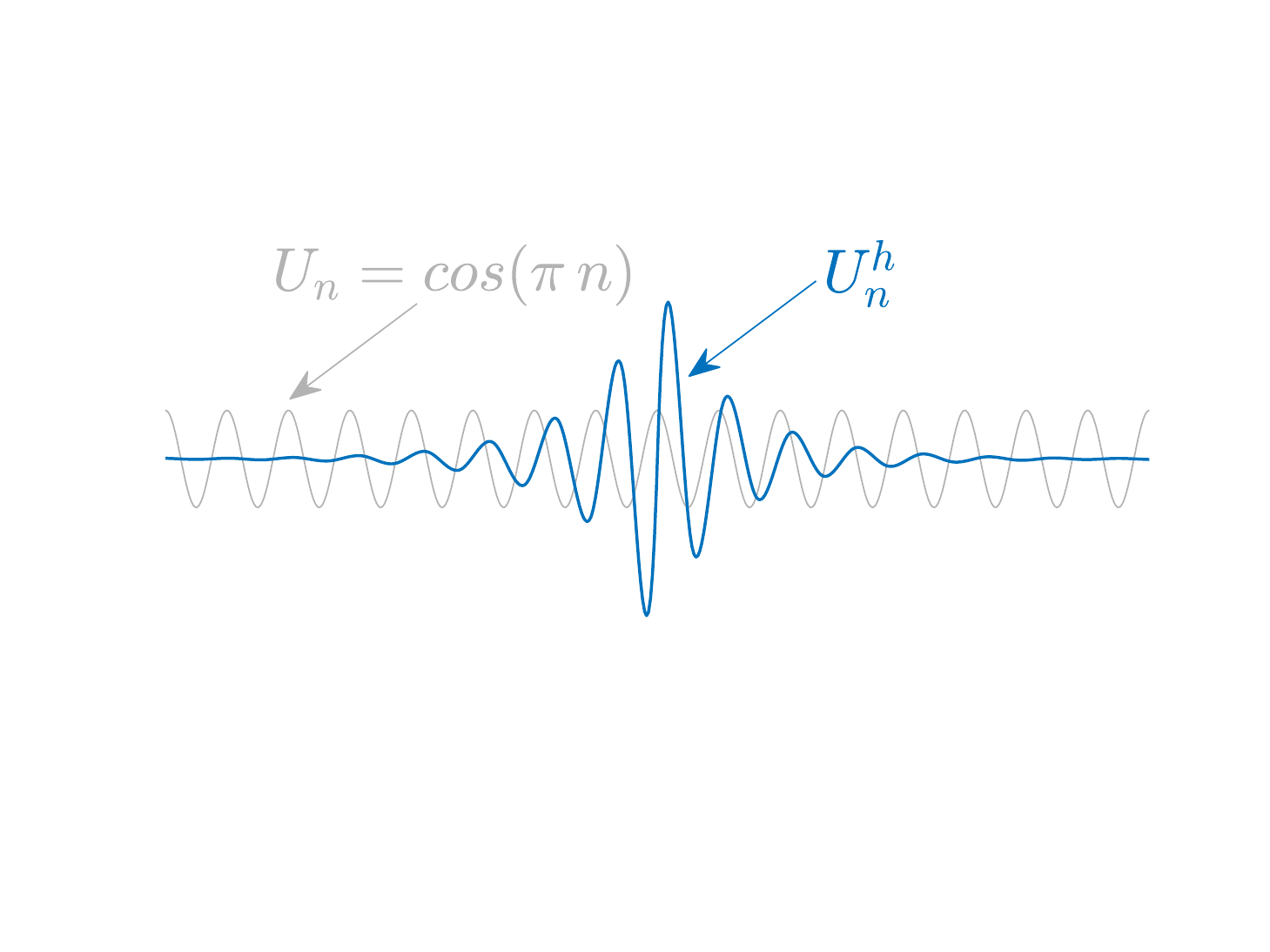} }}
    \subfloat[mode $N-2$]{{\includegraphics[width=0.43\textwidth]{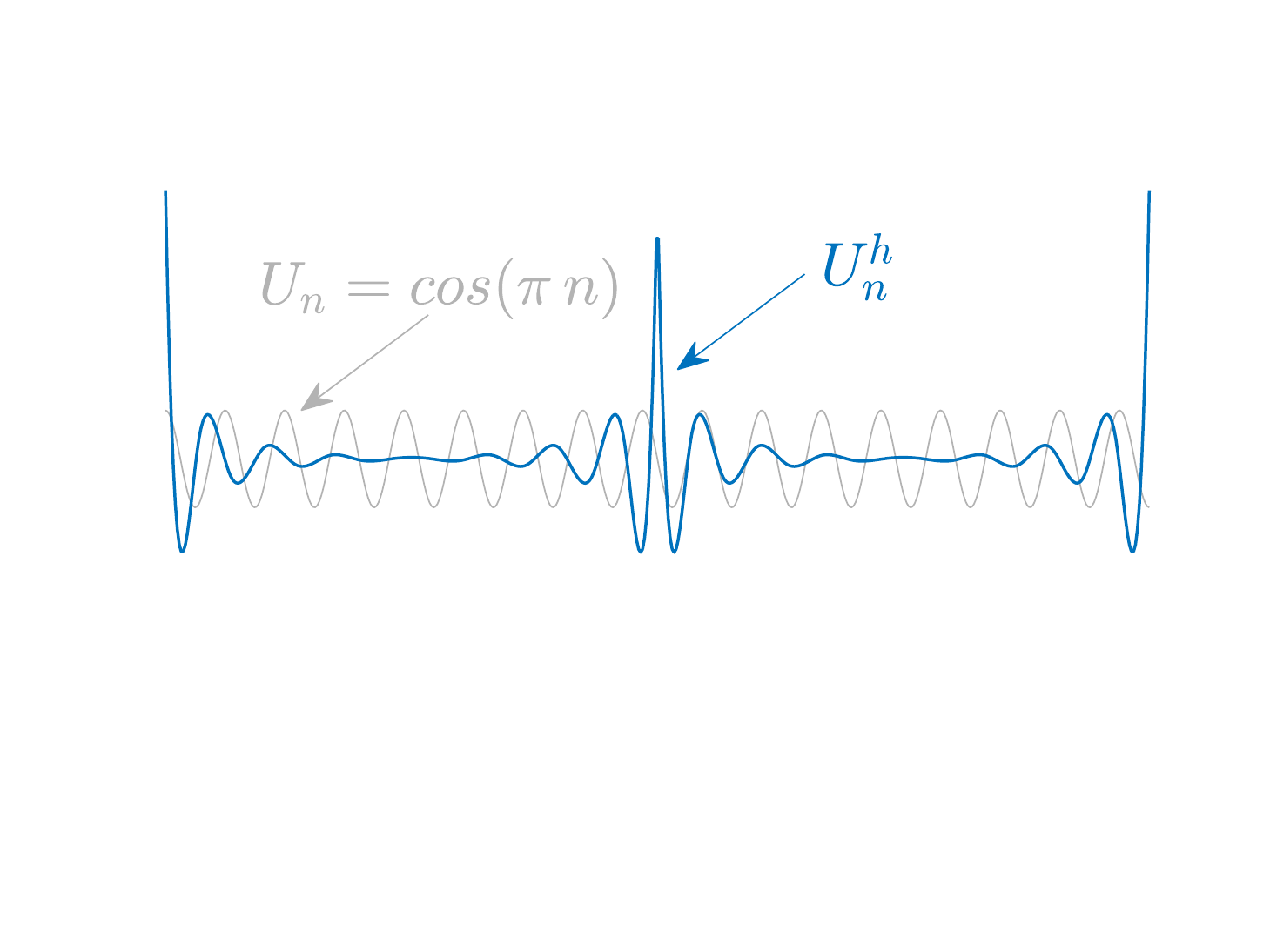} }}

    \subfloat[mode $N-1$]{{\includegraphics[width=0.43\textwidth]{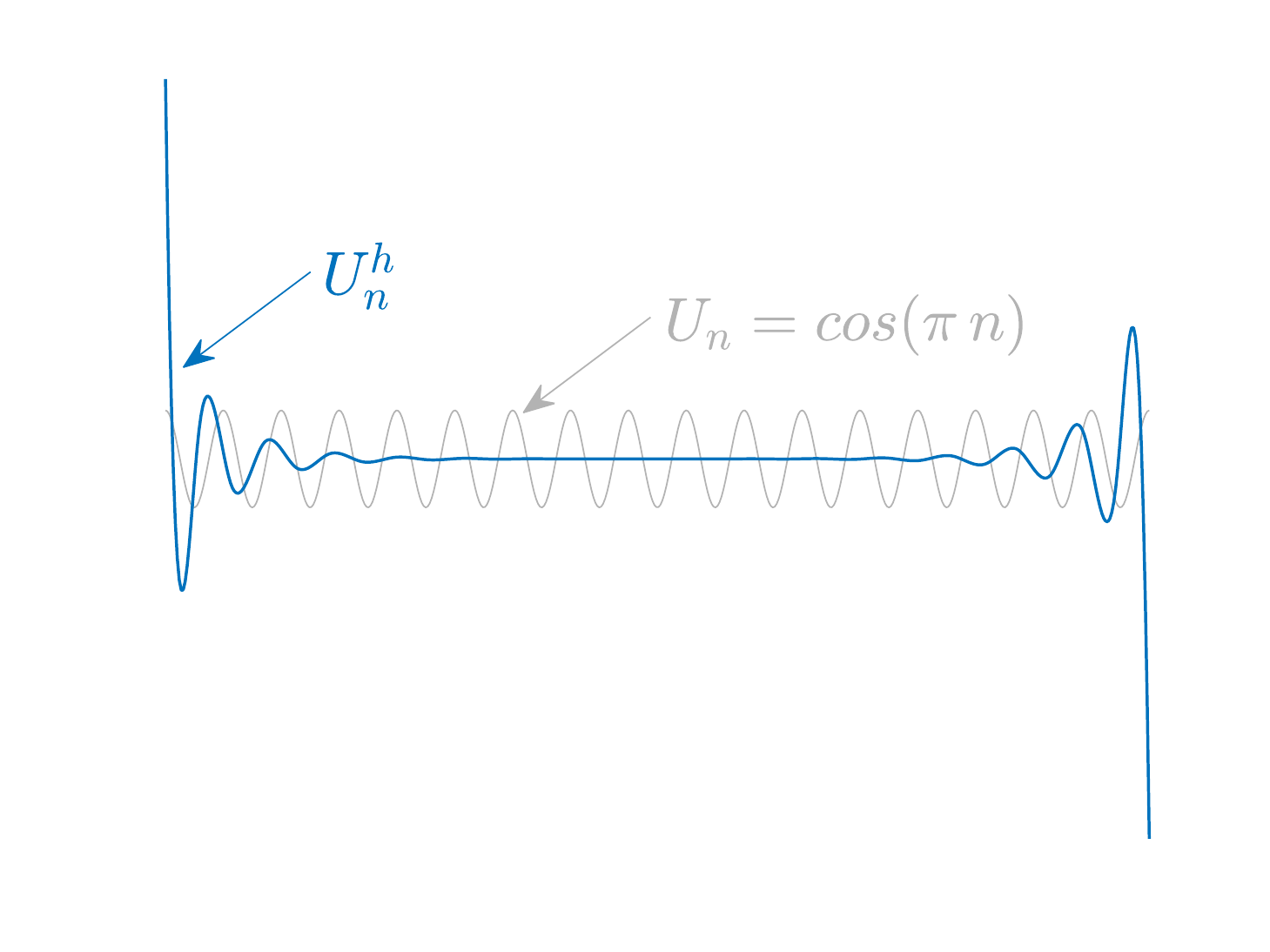} }}
    \subfloat[mode $N$]{{\includegraphics[width=0.43\textwidth]{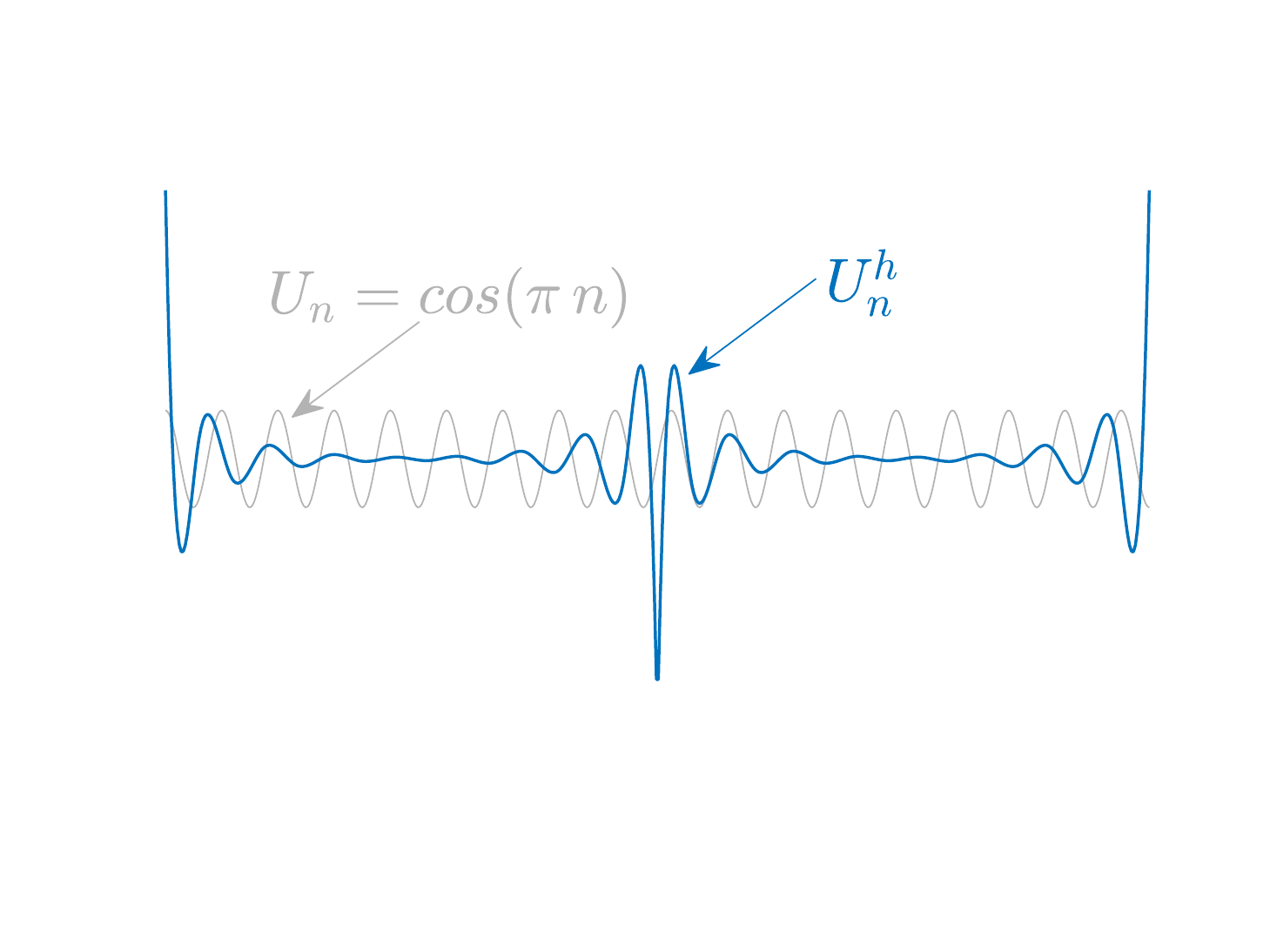} }}

    \caption{Discrete outlier modes $\eigenvec_n^h$ (\textbf{blue}) corresponding to the example of \textbf{the free bar studied in Figure \ref{fig:spectra_1Dbar_increasing_patches}}, computed with \textbf{two patches of cubic $C^2$ B-splines ($p=3$), and $C^0$ patch continuity}. Each patch is discretized with $15$ elements. These modes behave significantly different \changed{compared to} the exact modes (\textbf{gray}) $\eigenvec_n = \cos (\pi \, \noMode), \noMode = 1, \ldots, N$, where $N=35$ is the number of modes.}
    \label{fig:outlier_mode_shapes2}
\end{figure}

\changed{Isogeometric analysis (IGA) was first introduced in 2005 as a computational framework to improve the integration of computer-aided design (CAD) and finite element analysis (FEA) \cite{hughes_isogeometric_2005}. Compared with classical $C^0$ FEA, isogeometric discretizations exhibit better spectral properties} \cite{cottrell_studies_2007, hughes_finite_2014, hughes_duality_2008, puzyrev_spectral_2018}. While the upper part of the spectrum in classical FEA is inaccurate \cite{strang_analysis_2008, hughes_finite_2003} and \changed{the errors diverge with increasing interpolation degree $p$}, almost the entire spectrum converges with increasing $p$ in the case of smooth isogeometric discretizations \cite{cottrell_studies_2007, hughes_finite_2014, hughes_duality_2008, puzyrev_spectral_2018}. A small number of modes, however, form the so-called \textit{optical branch} at the end of the spectra and are denoted as ``outliers'' \cite{cottrell_isogeometric_2006, cottrell_isogeometric_2009}. 
\changedV{Isogeometric discretizations using multiple patches with lower smoothness at patch interfaces increase the number of spurious outlier frequencies \cite{puzyrev_spectral_2018}, as illustrated in Figure \ref{fig:spectra_1Dbar_increasing_patches}.} In the case where the number of patches equals the number of elements (dark red curves), the outliers form the well-known three branches of cubic $C^0$ finite elements in classical FEA, see e.g.\ \cite{hughes_finite_2014}.
We can also see that the highest frequencies are significantly overestimated, which may negatively affect
the stable critical \changed{time-step size} in explicit dynamics calculations.
The outlier modes, \changed{illustrated in Figure \ref{fig:outlier_mode_shapes} for a single-patch discretization and Figure \ref{fig:outlier_mode_shapes2} for a two-patch discretization}, behave in a spurious manner and may negatively affect the \changedV{solution accuracy} and robustness, particularly in hyperbolic problems \cite{hughes_finite_2014}.
We observe two types of outliers: one with all \changed{strain energy} located near the boundaries (\changedV{Figures} \ref{fig:outlier_mode_shapes}, \ref{fig:outlier_mode_shapes2}c) and one near the patch interface (\changedV{Figure} \ref{fig:outlier_mode_shapes2}a), to which we refer as boundary and interior outliers, respectively. For an in-depth discussion of these outlier types, we refer to \cite{puzyrev_spectral_2018} (interior outliers) and \cite{cottrell_isogeometric_2006,hiemstra_outlier_2021} (boundary outliers).  
Some of the outlier modes result from the combination of boundary and interior outliers, as illustrated in Figures \ref{fig:outlier_mode_shapes2}b and d.
\changed{We list the number of interior outliers in Tables \ref{tab:number_of_outliers_1Dbar} and \ref{tab:number_of_outliers_1Dbeam} for univariate multipatch discretizations of second-order and fourth-order problems, respectively, after removing the boundary outliers using the reduced spline basis introduced in \cite{hiemstra_outlier_2021}.} 
We count the number of interior outliers as a function of the degree $p$ and the number of patches $\npa$.
We also note that the total number of outliers is obtained by adding the number of interior outliers when the boundary ones are removed, and the number of boundary outliers when the interior outliers are not present.
For counting boundary outliers dependent on the boundary types and the degree $p$, we refer to \cite{hiemstra_outlier_2021}. 
For multivariate tensor-product discretizations,
the authors of \cite{hiemstra_outlier_2021} also provide formulas for the number of boundary outliers for multivariate tensor-product discretizations, which can be straightforwardly extended to discretizations with interior outliers.

\begin{table}[h!]
	\centering
	{\small
	\begin{tabularx}{0.9\linewidth}{ >{\centering\arraybackslash\hsize=.22\hsize}X | >{\centering\arraybackslash\hsize=.12\hsize}X >{\centering\arraybackslash\hsize=.12\hsize}X >{\centering\arraybackslash\hsize=.12\hsize}X >{\centering\arraybackslash\hsize=.12\hsize}X >{\centering\arraybackslash\hsize=.2\hsize}X }
	\toprule
	\multirow{2}{*}{Number of patches} & \multicolumn{5}{c}{Polynomial degree} \\
	  	  		& $2$ 		 & $3$ 			 & $4$  		 & $5$  		 & $p$ \\
	\midrule
	 $2$  		& $1$ 		 & $2$ 			 & $3$  		 & $4$  		 & $p-1$    \\
	 $3$  		& $2$ 		 & $4$ 			 & $6$  		 & $8$  		 & $2(p-1)$   \\
	 $\npa$   & $\npa-1$ & $2(\npa-1)$ & $3(\npa-1)$ & $4(\npa-1)$ & $(\npa-1)(p-1)$  \\
	\bottomrule
\end{tabularx}}
\caption{Number of interior outlier modes in one-dimensional \textbf{multipatch discretizations with $C^0$ patch continuity} applied to a second-order problem, e.g.\ the axial vibration of \textbf{a bar}. 
}
\label{tab:number_of_outliers_1Dbar}
\end{table}

\begin{table}[h!]
	\centering
	{\small
	\begin{tabularx}{0.9\linewidth}{ >{\centering\arraybackslash\hsize=.22\hsize}X | >{\centering\arraybackslash\hsize=.12\hsize}X >{\centering\arraybackslash\hsize=.12\hsize}X >{\centering\arraybackslash\hsize=.12\hsize}X >{\centering\arraybackslash\hsize=.12\hsize}X >{\centering\arraybackslash\hsize=.2\hsize}X }
	\toprule
	\multirow{2}{*}{Number of patches} & \multicolumn{5}{c}{Polynomial degree} \\
					& $3$ 			 & $4$  		 & $5$  		& $6$	 		 & $p$ \\
	\midrule
	 $2$  		 	& $1$ 			 & $2$  		 & $3$  		& $4$	 		 & $p-2$    \\
	 $3$  		 	& $2$ 			 & $4$  		 & $6$  		& $8$	 		 & $2(p-2)$   \\
	 $\npa$    		& $\npa-1$    	 & $2(\npa-1)$   & $3(\npa-1)$ 	& $4(\npa-1)$    & $(\npa-1)(p-2)$  \\
	\bottomrule
\end{tabularx}}
\caption{Number of interior outlier modes in one-dimensional \textbf{multipatch discretizations with $C^1$ patch continuity} applied to a fourth-order problem, e.g.\ the transverse vibration of \textbf{a beam}.}
\label{tab:number_of_outliers_1Dbeam}
\end{table}

The first idea on how to remove boundary outliers in an IGA context was based on the nonlinear parameterization of the domain via a uniform distribution of the control points \cite{cottrell_isogeometric_2009}. 
On the one hand, this approach removes the outliers for arbitrary degree $p$. On the other hand, a nonlinear parameterization changes the original geometry representation from CAD, which contradicts the isogeometric paradigm of using the same geometry in design and analysis. Furthermore, in \cite{hiemstra_outlier_2021}, the authors verified that this approach leads to a loss of spatial accuracy of the low modes and frequencies.
\changedV{In more recent contributions \cite{hiemstra_outlier_2021,deng_outlier_2021,Manni2021}, the authors imposed additional boundary constraints arising from higher-order eigenvalue problems, either by building constraints into the basis \cite{hiemstra_outlier_2021,Manni2021}} or via penalization \cite{deng_outlier_2021}, to improve the spectral properties of isogeometric discretizations.  
The strong approach entirely removes the outlier frequencies and modes for arbitrary degree $p$ in one- and multidimensional settings. \changedV{The penalty approach reduces the outlier frequencies, but does not remove the corresponding spurious outlier
modes.} An alternative treatment is the penalty approach introduced in \cite{horger_penalty_2019} that imposes additional higher-order continuity constraints at interfaces of multipatch discretizations, as well as the first-order derivative at the Neumann boundary, i.e.\ a penalty treatment of both boundary and interior outliers.
The authors reported an improvement of the eigenvalue spectra after removing ``unphysical'' modes using a cut-off normalized eigenvalue.
The penalty approach of \cite{horger_penalty_2019}, using large values of penalty parameters, serves \changed{purely} as an indicator for the outlier eigenvalues / frequencies such that these are identified easily, as verified in subsequent sections in this paper. 
Another approach is to weakly enforce the continuity constraints at patch interfaces in the framework of mortar methods \cite{Dittmann2019,Schuss2019}, or to apply  
the \textit{optimally-blended} quadratures \cite{ainsworth_blended_quadrature_2010} that can suppress the boundary outlier \cite{puzyrev_quadrature_2017, calo_quadrature_2019,deng_outlier2_2021} as well as the interior outlier frequencies \cite{puzyrev_spectral_2018}.

Alternatively, the highest frequencies may be reduced via mass scaling. Its idea is to add artificial terms to the mass matrix in such a way that high frequencies are affected and any negative impact on the lower frequencies and modes is kept at a minimum. A widely used variant scales the density in combination with mass lumping, see e.g. \cite{hartmann_mass_scaling_2015}. 
One existing mass-scaling technique is to add a weighting of some form of stiffness matrix as a mass scaling \cite{Macek1995,olovsson_mass_scaling_2005,olovsson_mass_scaling_2006,Askes2011}. The approach is then called \textit{selective mass scaling} when it targets specific frequencies and modes \cite{olovsson_mass_scaling_2005}. In another technique, the added mass follows from a penalized Hamilton's principle \cite{tkachuk_mass_scaling_2014, schaeuble_mass_scaling_2017}, which is a variationally consistent approach. Further artificial mass terms are also developed to optimize accuracy and efficiency in  e.g. \cite{cocchetti_mass_scaling_2013,gonzalez_mass_tailoring_2020}.

In this paper, we introduce a novel variational approach based on perturbed eigenvalue analysis that improves the spectral properties of isogeometric multipatch discretizations.
We combine the ideas of penalizing both the stiffness and \changed{the mass matrix} \cite{deng_outlier_2021} and of adding \changed{higher-order continuity constraints at patch interfaces \cite{horger_penalty_2019}} to arrive at an improved suppression of the interior outlier frequencies.
In particular, we add scaled perturbation terms that weakly enforce the patch continuity constraints of \cite{horger_penalty_2019} to both the stiffness and \changed{the mass matrix}. 
We note that this approach results in modified left- and right-hand sides of the standard formulation in explicit dynamics, where the term involving the stiffness matrix affects the right-hand side residual. This differs from a mass scaling approach which modifies only the mass matrix. 
To remove boundary outliers, we combine the proposed variational approach with the methodology introduced in \cite{hiemstra_outlier_2021}. 
The proposed approach is consistent given that the analytical solution satisfies the patch continuity constraints of \cite{horger_penalty_2019}, i.e.\ \changed{the solution is sufficiently smooth}, as well as the additional boundary constraints of \cite{hiemstra_outlier_2021}.
Moreover, we introduce an approach for estimating optimal scaling parameters of the perturbation term, in the sense that the outlier frequencies are effectively reduced and \changedV{the accuracy in the remainder of the spectrum and modes} is not negatively affected. We also show how this approach can be cast into a pragmatic iterative procedure that can be readily implemented in any IGA framework.

We discuss different perturbation variants, such as perturbation of the stiffness matrix only (equivalent with the approach of \cite{horger_penalty_2019}), perturbation of the mass matrix only (equivalent with a selective mass scaling approach), and perturbation of both stiffness and mass matrices.
The proposed iterative procedure can also be applied to approximate the optimal scaling parameters for the approach of \cite{deng_outlier_2021} such that the boundary outliers are optimally suppressed. 
We verify numerically via spectral analysis of second- and fourth-order problems that the proposed approach improves spectral properties of isogeometric multipatch discretizations in one- and multidimensional settings. For the examples of membrane and plate structures in an explicit dynamics setting, we confirm that our approach maintains spatial accuracy and enables a larger critical time-step size. We also demonstrate that it is effective irrespective of 
the polynomial order $p$.

The structure of the paper is as follows: 
In Section \ref{sec:var-form}, we derive the variational formulation based on perturbed eigenvalue analysis. 
In Section \ref{sec:1d-study}, we motivate the iterative scheme for parameter estimation, focusing on a one-dimensional problem. 
In Section \ref{sec:2d-study}, we generalize our approach 
to multidimensional discretizations, including its practical implementation, and demonstrate its effectiveness for discretizations of  
second- and fourth-order problems. 
In Section \ref{sec:results-dyn}, we discuss its advantages in explicit dynamics. 
In Section \ref{sec:conclusion}, we summarize our results and draw conclusions.

\section{Variational formulation}\label{sec:var-form}

We start this section with a brief review of the equations of motion that govern free vibrations of an undamped linear structural system, and derive the corresponding generalized eigenvalue problem in the continuous and discrete settings.
We then  
introduce a variational formulation based on a perturbed eigenvalue problem that weakly enforces additional continuity constraints at patch interfaces. These additional constraints suppress only the interior outliers and do not negatively affect the important low frequency part.

\subsection{Natural frequencies and modes}

\changed{The equation of motion that governs the free vibration of an undamped linear structural system is:}
\begin{align}
    \mathcal{K} \, u(\vect{x}, t) + \mathcal{M} \, \frac{d^2}{d \, t^2} \, u(\vect{x}, t) \; = \;  0 \, , \quad \vect{x} \in \Omega\, , \quad t > 0 \, . \label{eom}
\end{align}
Here, $\mathcal{M}$ and $\mathcal{K}$ are the mass and stiffness operators, respectively, and $u(\vect{x}, t)$ is the displacement \changed{of a point $\vect{x}$ in} the domain $\Omega$.
Using separation of variables, the displacement can be expanded in terms of the eigenmodes $\eigenvec_\noMode(\vect{x})$ and the time-dependent coefficients $T_\noMode(t)$, that is $u(\vect{x}, t) = \sum_\noMode \eigenvec_\noMode(\vect{x}) \; T_\noMode(t)$. 
\changed{Substitution in \eqref{eom} leads to two results. \changedV{Firstly,}  
$T_\noMode(t) = C_+ \, e^{i\,\omega_\noMode\,t} + C_- \, e^{-i\,\omega_\noMode\,t}$, which satisfies the equation:} 
\begin{align}
    \frac{d^2}{d\,t^2} T_\noMode(t) + \omega_\noMode^2 \, T_\noMode(t) \; = \;  0 \, ,
\end{align}
and describes an oscillation at a frequency $\omega_\noMode$. \changed{Here, $C_+$ and $C_-$ are constants determined from initial conditions}.
Secondly, it results in 
the strong form of the generalized eigenvalue problem in the continuous setting, that is: 
find eigenmodes $\eigenvec_\noMode(\vect{x})$ and eigenfrequencies $\omega_\noMode$, 
$(U_\noMode, \omega_\noMode) \in \mathcal{V} \times \mathbb{R}^+$ such that:
\begin{equation}
    \left(\mathcal{K} - \omega_\noMode^2 \mathcal{M}\right) \, \eigenvec_\noMode(\vect{x}) \; = \;  0 \, . 
\label{gep} 
\end{equation}
Here, $\mathcal{V}$ is the space of functions with sufficient regularity that allows the differential operators in $\mathcal{M}$ and $\mathcal{K}$ to be applied.

Applying the standard Galerkin method and subsequently discretizing with $N$ finite element basis functions $B_i (\vect{x})$ results in the following semi-discrete system of equations:
\begin{align}
    \mat{K} \, \mat{u}^h(t) + \mat{M} \, \frac{d^2}{d\,t^2} \, \mat{u}^h(t) \; = \;  0 \, , \label{deom}
\end{align}
where $\mat{K}$ and $\mat{M}$ denote the stiffness and consistent mass matrix, respectively, and $\mat{u}^h(t)$ is the unknown time-dependent displacement vector, \changed{such that:} 
\begin{align}
    u^h(\vect{x},t) = \begin{bmatrix}
          B_1 (\vect{x}) \; \ldots \; B_N (\vect{x}) 
      \end{bmatrix} \, \mat{u}^h(t) \, , \qquad u^h(\vect{x},t) \in \mathcal{V}^h \subset \mathcal{V} \, . \nonumber
\end{align}
Here, $u^h(\vect{x},t)$ is the discrete displacement and \changedV{$\mathcal{V}^h$ is a discrete space spanned by sufficiently smooth basis functions $B_i(\vect{x}), \; i=1, \ldots, N$.}
The corresponding discrete eigenvalue problem can be expressed in the following matrix equation:
\begin{equation} \label{dgep}
    \mat{K} \, \mat{\eigenvec}_\noMode^h \; = \; \left(\omega_\noMode^h \right)^2 \, \mat{M} \, \mat{\eigenvec}_\noMode^h \, ,
\end{equation}
where $\mat{\eigenvec}_\noMode^h$ denotes the vector of unknown coefficients corresponding to the $\noMode^{\text{th}}$ discrete eigenmode $\eigenvec_\noMode^h$, and 
$\omega_\noMode^h$ is the $\noMode^{\text{th}} $ discrete eigenfrequency.

\subsection{Perturbed second- and fourth-order eigenvalue problems}

In this paper, we focus on second-order eigenvalue problems involving rods and membranes, and fourth-order problems involving  
square plate structures. 
These are of unit size and unit material parameters and include either homogeneous Dirichlet or homogeneous Neumann boundary conditions. 
Since we consider spaces with sufficient regularity, as discussed in the previous subsection, the minimum patch continuity of multipatch discretizations is $C^0$ and $C^1$ for second- and fourth-order problems, respectively. 
We employ spaces of $C^{p-1}$ B-splines of degree $p$ that are free of boundary outliers, i.e.\ spaces with outlier removal constraints that are strongly enforced at the boundary \cite{hiemstra_outlier_2021}.

\begin{remark}
    The choice of unit material parameters, i.e. unit mass and stiffness, is not a realistic scenario and does not represent the true conditioning of the problem. In this paper, however, this choice is trivial since we mostly look at the normalized eigenfrequency.     
    We plan to further study realistic scenarios and the problem conditioning in future research.
\end{remark}

The discrete eigenvalue problem \eqref{dgep} can be expressed in the following variational form:
find $(U^h_\noMode, \omega^h_\noMode) \in \mathcal{V}^h \times \mathbb{R}^+$, for $n = 1,2,\ldots,N$, such that:
\begin{align}\label{vdgep}
    a(\eigenvec_\noMode^h, v^h) = \left(\omega^h_\noMode \right)^2 \, b(\eigenvec_\noMode^h, v^h) \qquad \forall \, v^h \in \mathcal{V}^h \subset \mathcal{V} \, .
\end{align}
Here, the bilinear form $b(\cdot, \cdot)$ is:
\begin{align}\label{eq:mass_bilinear_form}
    b(u^h, w^h) = \int_{\Omega} \, u^h \, w^h \, \mathrm{d} x \, .
\end{align}
The bilinear form $a(\cdot, \cdot)$ of a discrete second-order eigenvalue problem is: 
\begin{align}
    a(u^h, w^h) = \int_{\Omega} \, \nabla \, u^h \, \cdot \, \nabla \, w^h \, \mathrm{d} x \, , \label{eq:stiffness_bilinear_form_2nd_order}
\end{align}
which corresponds to the second-order stiffness operator of the strong form \eqref{gep}:
\begin{align}
    \mathcal{K} := - \laplace = - \sum_{k=1}^d \frac{\partial^2}{\partial x_k^2} \qquad d = 1,2,3 \,, \label{eq:2nd_order_operator}
\end{align}
where $d$ is the dimension of the problem.
 
For fourth-order problems, the bilinear form $a(\cdot, \cdot)$ is:
\begin{align}
    a(u^h, w^h) = \int_{\Omega} \, \Delta \, u^h \; \Delta \, w^h \, \mathrm{d} x \, , \label{eq:stiffness_bilinear_form_4th_order}
\end{align}
which corresponds to the bi-harmonic operator:
\begin{align}
    \mathcal{K} = \Delta^2 = \Delta \, \Delta \, . \label{eq:4th_order_operator}
\end{align}
In two-dimensional settings, i.e. the case of a vibrating plate, we consider the simply supported boundary conditions in this paper, since only the analytical solution of this case is known \cite{szilard_theories_2004}.

For improving the spectral properties of isogeometric multipatch discretizations, focussing on the interior outliers, 
we weakly enforce the \changed{following $C^{p-1}$ continuity} constraints at patch interfaces:
\begin{align}\label{eq:continuity_constraint}
    \left \lsem \partial^l_{\normalvect} \eigenvec_\noMode(\vect{x}) \right \rsem = 0 \qquad \text{on} \;\; \Gamma^e \; , \; l = 1,\ldots p-1 \;\; ,
\end{align}
where $\Gamma^e$ denotes the \changed{$e^{\text{th}}$ patch interface, $e = 1,\ldots,\npai$ and $\npai$ is the number of patch interfaces}. $\normalvect$ denotes the outward unit normal to the patch interface, and $\left \lsem \cdot \right \rsem$ denotes the jump across the interface, $\left \lsem w \right \rsem = w_+ - w_-$. 
The $l^\text{th}$ constraint of \eqref{eq:continuity_constraint}, i.e.\ the $C^l$ continuity constraint, at the $e^{\text{th}}$ patch interface 
corresponds to the bilinear form: 
\begin{align}
    c^l_e(u^h, w^h) = \int_{\Gamma^e} \, \left \lsem \partial^l_{\normalvect} u^h \right \rsem \, \left \lsem \partial^l_{\normalvect} w^h \right \rsem \, \mathrm{d} x \, . \label{eq:perturbation_bilinear_form}
\end{align}

We propose the following variational formulation for perturbing the eigenvalue problem \eqref{vdgep} that weakly enforces \eqref{eq:continuity_constraint}:
find $(\tilde{\eigenvec}^h_\noMode, \tilde{\omega}^h_\noMode) \in \mathcal{V}^h \times \mathbb{R}^+$, for $n = 1,2,\ldots,N$, such that $\forall \, v^h \in \mathcal{V}^h \subset \mathcal{V}$:
\begin{align}
    a(\tilde{\eigenvec}_\noMode^h, v^h) + \sum_{l=1}^{p-1} \, \sum_{e=1}^{\npai} \, \alpha^l_e \, c^l_e(\tilde{\eigenvec}_\noMode^h, v^h) = \left(\tilde{\omega}_\noMode^h\right)^2 \, \left[b(\tilde{\eigenvec}_\noMode^h, v^h) + \sum_{l=1}^{p-1} \, \sum_{e=1}^{\npai} \, \beta^l_e \, c^l_e(\tilde{\eigenvec}_\noMode^h, v^h) \right] \, , \label{eq:vdgep_perturbed0}
\end{align}
where $\alpha^l_e$ and $\beta^l_e$ are scaling factors of the perturbation $c^l_e(\cdot,\cdot)$.  
The tilde in the superscript of $\tilde{\eigenvec}^h_\noMode$ and $\tilde{\omega}^h_\noMode$ distinguishes the eigenmode and frequency corresponding to the perturbed eigenvalue problem \eqref{eq:vdgep_perturbed} from those corresponding to the standard non-perturbed problem \eqref{vdgep}.

Based on empirical observation, we find that, \changedV{for uniform discretizations, $\alpha^l_1 = \ldots = \alpha^l_{\npai} = \alpha^l$ and $\beta^l_1 = \ldots = \beta^l_{\npai} = \beta^l$.}
The variational formulation \eqref{eq:vdgep_perturbed0} then becomes: 
find $(\tilde{\eigenvec}^h_\noMode, \tilde{\omega}^h_\noMode) \in \mathcal{V}^h \times \mathbb{R}^+$, for $n = 1,2,\ldots,N$, such that $\forall \, v^h \in \mathcal{V}^h \subset \mathcal{V}$:
\begin{align}
    a(\tilde{\eigenvec}_\noMode^h, v^h) + \sum_{l=1}^{p-1} \, \alpha^l \, c^l(\tilde{\eigenvec}_\noMode^h, v^h) = \left(\tilde{\omega}_\noMode^h\right)^2 \, \left[ b(\tilde{\eigenvec}_\noMode^h, v^h) + \sum_{l=1}^{p-1} \, \beta^l \, c^l(\tilde{\eigenvec}_\noMode^h, v^h) \right] \, , \label{eq:vdgep_perturbed}
\end{align}
where
\begin{align}
    c^l(u^h, w^h) = \sum_{e=1}^{\npai} \, c^l_e(u^h, w^h) \, . \label{eq:perturbation_bilinear_form2}
\end{align}
The matrix equation of the perturbed eigenvalue problem \eqref{eq:vdgep_perturbed} is:
\begin{align}
    \left( \, \mat{K} + \sum_{l=1}^{p-1} \, \alpha^l \, \mat{K}_{\Gamma}^l  \right) \, \mat{\tilde{\eigenvec}}_\noMode^h \; = \; \left(\tilde{\omega}_\noMode^h \right)^2 \, \left( \, \mat{M} + \sum_{l=1}^{p-1} \, \beta^l \, \mat{K}_{\Gamma}^l \right) \, \mat{\tilde{\eigenvec}}_\noMode^h \, , \label{eq:dgep_perturbed}
\end{align}
where the stiffness matrix $\mat{K}$ and the consistent mass matrix $\mat{M}$ correspond to the bilinear forms \changed{$a(\cdot,\cdot)$ from \eqref{eq:stiffness_bilinear_form_2nd_order} or \eqref{eq:stiffness_bilinear_form_4th_order} and $b(\cdot,\cdot)$ from \eqref{eq:mass_bilinear_form}}, respectively, \changed{which} are symmetric positive definite matrices;
and the perturbation matrix $\mat{K}_{\Gamma}^l$ corresponds to the bilinear form $c^l(\cdot,\cdot)$ \changed{from \eqref{eq:perturbation_bilinear_form2}}, \changedV{which is a symmetric positive semi-definite matrix}. 

We note that the proposed approach \eqref{eq:vdgep_perturbed0} is applicable \changed{to} non-uniform discretizations while the approach \eqref{eq:vdgep_perturbed} is developed for uniform mesh discretizations. In this paper, we focus on uniform discretizations.

\subsection{A note on consistency vs.\ variational consistency}

Consistency and variational consistency are important properties of a finite element formulation. They play key roles in error analysis and \changed{are necessary} for ensuring optimal orders of convergence of the method \cite{strang_analysis_2008,hughes_finite_2003,ern_fem_2004}. Additionally, they guarantee that the method yields the true solution if that solution lies in the trial function space. 
The variational approach proposed in the previous subsection is consistent when the true solution is sufficiently smooth, but it is not variationally consistent. 

\textit{Variational consistency} relates to the equivalence of the strong and weak forms in the limit $h\rightarrow 0$, with $h$ being the characteristic mesh size. If the formulation is also stable, and when the data (body force and boundary conditions) are sufficiently regular, then variational consistency ensures that the finite element approximation converges to the strong solution with mesh refinement \cite{strang_analysis_2008,hughes_finite_2003,ern_fem_2004}. The variational formulation \eqref{eq:vdgep_perturbed0} is not equivalent to the strong formulation of the eigenvalue problem \eqref{gep} 
as $h\rightarrow 0$. Instead, we obtain the following alternative strong formulation after performing integration by parts on a patch level:
\begin{subequations}
    \begin{align}
        & \left(\mathcal{K} - \left(\omega_\noMode^h\right)^2 \mathcal{M}\right) \, \eigenvec_\noMode(\vect{x}) \; = \;  0 & \qquad \text{in} \;\; \Omega \, , \\
        & \left \lsem \partial^l_{\normalvect} \, \eigenvec_\noMode(\vect{x}) \right \rsem \; = \; 0 & \qquad \text{on} \;\; \Gamma \; , \; l = 1,\ldots p-1 \; ,
    \end{align}
\end{subequations}
where $\Gamma$ is a collection of all patch interfaces. 
This corresponds to the following strong form of the governing equation of free vibrations: 
\begin{subequations} \label{varcons}
    \begin{align} 
        & \mathcal{K} \, u(\vect{x}, t) + \mathcal{M} \, \frac{d^2}{d \, t^2} \, u(\vect{x}, t) \; = \;  0 & \qquad \text{in} \;\; \Omega \, , \\
        & \left \lsem \partial^l_{\normalvect} \, u(\vect{x},t) \right \rsem \; = \; 0 & \qquad \text{on} \;\; \Gamma \; , \; l = 1,\ldots p-1 \; ,
    \end{align}
\end{subequations}
which thus involves more constraints on jumps of higher order derivatives across the patch interface compared to the strong form \eqref{eom}. 

In a broader sense, \textit{consistency} refers to some specific solution. A formulation is consistent with respect to the solution $u_{true}$ if the variational statement is satisfied when $u_{true}$ is substituted in the trial function slot. This property immediately implies a form of Galerkin orthogonality with respect to $u_{true}$, i.e.:
\begin{align}
a(u_{true}-u^h, v^h) = 0 \quad \forall\, v^h \in \mathcal{V}^h \subset \mathcal{V}
\end{align}
C\'ea's lemma then only requires that the formulation satisfies a stability criterion (i.e., coercivity) to yield the optimal order of convergence to the solution $u_{true}$ in the natural norm \cite{strang_analysis_2008,ern_fem_2004}:
\begin{align}
|| u_{true} - u^h ||_{H^1} \leq \frac{C_b}{C_c}\min\limits_{v^h\in \mathcal{V}^h} || u_{true} - v^h ||_{H^1}
\end{align}
where $C_b$ is the boundedness coefficient and $C_c$ is the coercivity coefficient.

Our formulation \eqref{eq:vdgep_perturbed0} is consistent with respect to solutions that satisfy the conditions of \eqref{varcons}. This means that we can expect to converge optimally to solutions that satisfy these conditions. Put simply, we can expect to converge optimally for solutions in \changed{$H^p(\Omega)$}. Such solutions often correspond to system response governed by free vibrations. For an isogeometric discretization, these are generally the solutions we focus on, since the optimal order of convergence \changed{to a function in $H^m(\Omega)$} in the $H^q(\Omega)$ norm is $\min(p+1,m)-q$, with $0 \leq q = \min(p+1,m,s)$, where $p$ is the polynomial degree of the B-spline basis function, $m$ is the order of smoothness of the solution, and $s$ is the minimum global regularity of the basis functions \cite{tagliabue_error_2014,bazilevs_error_2006}. 
We thus see that our formulation is consistent with respect to solutions for which the B-spline discretization yields optimal convergence and not for those solutions for which the B-spline discretization convergences suboptimally anyway.

\section{Parameter estimation for a one-dimensional case study}\label{sec:1d-study}

In the next step, we will address the systematic choice of the open parameters $\alpha$ and $\beta$ in \eqref{eq:dgep_perturbed}. To approach this aspect, we consider the perturbed eigenvalue problem \eqref{eq:dgep_perturbed} of a fixed bar (unit length, unit axial stiffness, unit mass), discretized by a univariate multipatch discretization with quadratic $C^1$ B-splines ($p=2$) and sufficient regularity, i.e.\ $C^0$ patch continuity. 
The perturbed eigenvalue problem \eqref{eq:dgep_perturbed} then simplifies to:
\begin{align}
    \left[ \mat{K} + \alpha \mat{K}_\Gamma \right] \, \mat{\tilde{\eigenvec}}_{n}^h = \left(\tilde{\omega}^h_{n}\right)^2 \, \left[ \mat{M} + \beta \mat{K}_\Gamma \right] \, \mat{\tilde{\eigenvec}}_{n}^h \, , \qquad n = 1,\ldots,N \; . \label{eq:dgep_perturbed_1dp2}
\end{align}
We consider a discretization of two patches ($\npa=2$). 
The resulting spectrum consists of one interior outlier frequency (see Table \ref{tab:number_of_outliers_1Dbar}) which is the maximum frequency.

\subsection{First-order approximation of the perturbation}

The objective of the perturbation
is to ensure that $\tilde{\omega}^h_\noMode$ is a good approximation to the true eigenfrequency over the complete spectrum. Specifically, the parameters $\alpha$ and $\beta$ need to be chosen to reduce the severely over-estimated maximum $\omega^h_\noMode$ without compromising the accuracy of the lower $\omega^h_\noMode$'s. Focusing on the highest frequency mode, relation \eqref{eq:dgep_perturbed_1dp2} becomes:
\begin{align}
    \mat{K}\mat{\tilde{\eigenvec}}_{max}^h + ( \alpha - \beta\left(\tilde{\omega}^h_{max}\right)^2 ) \mat{K}_\Gamma\mat{\tilde{\eigenvec}}_{max}^h = \left(\tilde{\omega}^h_{max}\right)^2 \mat{M}\mat{\tilde{\eigenvec}}_{max}^h \, , \label{step1}
\end{align}
and after premultiplying by $\mat{\tilde{\eigenvec}}_{max}^{h\, T}$:
\begin{align}
    \mat{\tilde{\eigenvec}}_{max}^{h\, T} \mat{K}\mat{\tilde{\eigenvec}}_{max}^h + ( \alpha - \beta\left(\tilde{\omega}^h_{max}\right)^2 ) \mat{\tilde{\eigenvec}}_{max}^{h\, T} \mat{K}_\Gamma \, \mat{\tilde{\eigenvec}}_{max}^h = \left(\tilde{\omega}^h_{max}\right)^2 \, , \label{step2}
\end{align}
where we assume that the eigenmodes are normalized with respect to the unperturbed mass matrix (i.e., $\mat{\tilde{\eigenvec}}_{n}^{h\, T} \mat{M} \mat{\tilde{\eigenvec}}_m^h = \delta_{nm}$). 

A first order approximation of a perturbed eigenvalue problem reveals that $\left(\omega^h_\noMode\right)^2$ changes with order $\mathcal{O}(\norm{\eigenvec^h_{n}}^2)$
and eigenmode $\eigenvec^h_{n}$ changes with order $\mathcal{O}(\norm{\eigenvec^{h}_{n}}^2/\left(\omega^{h}_{max}\right)^2)$ \cite[Section~15.4]{li_perturbation_2014}.
Since $\omega^{h}_{max}$ is large, the relative change in eigenmodes can be expected to be much smaller than the relative change in frequencies. This implies that, for small $\alpha$ and $\beta$, we may approximate $\mat{\tilde{\eigenvec}}_{max}^h \approx \mat{\eigenvec}^{h}_{max}$ and thus also:
\begin{align}
    \mat{\tilde{\eigenvec}}_{max}^{h\, T} \mat{K}\mat{\tilde{\eigenvec}}_{max}^h \approx \mat{\eigenvec}_{max}^{h\, T} \mat{K}\mat{\eigenvec}^{h}_{max} = \left(\omega^{h}_{max}\right)^2 \, \mat{\eigenvec}_{max}^{h\, T}\mat{M} \mat{\eigenvec}_{max}^h = \left(\omega^{h}_{max}\right)^2 \, . \label{eq:1st_order_approx}
\end{align}

\subsection{Identifying (un)suitable parameter windows}

To help identify suitable ranges of parameter values, we choose to write $\beta$ in terms of $\alpha$ and a scaling factor $f$ as:
\begin{align}
    \beta = f \, \frac{1}{\left(\tilde{\omega}^{h}_{max}\right)^2} \, \alpha \, . \label{eq:beta}
\end{align}
Substitution of \changed{\eqref{eq:1st_order_approx} and \eqref{eq:beta} into \eqref{step2}} gives:
\begin{align}
    \alpha\, ( 1  -  f )  \approx \frac{\left(\tilde{\omega}^{h}_{max}\right)^2 -\left(\omega^{h}_{max}\right)^2 }{ \mat{\eigenvec}^{h\, T}_{max} \mat{K}_\Gamma\mat{\eigenvec}^{h}_{max}}  \, . \label{eq:perturbed_alpha}
\end{align}
Depending on our choice of the scaling factor $f$, we can identify different parameter windows that we will briefly discuss in the following.\\

\textbf{Case $\boldsymbol f = \boldsymbol 0$} (no perturbation of the mass matrix): The left-hand side is positive for positive choice of $\alpha$. This necessarily means that $\tilde{\omega}_{max}^h > \omega^{h}_{max}$, which is precisely not our goal: the frequency $\omega^{h}_{max}$ needs to be reduced. Adding perturbations only to the stiffness matrix is thus unsuitable for improving the spectrum. 

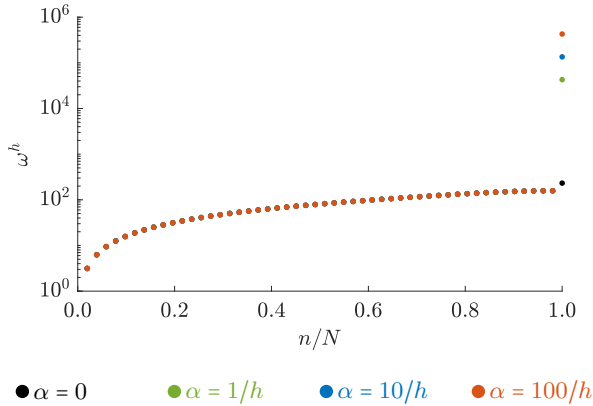
\begin{figure}[h!]
    \centering 
    \includegraphics[width=0.5\textwidth]{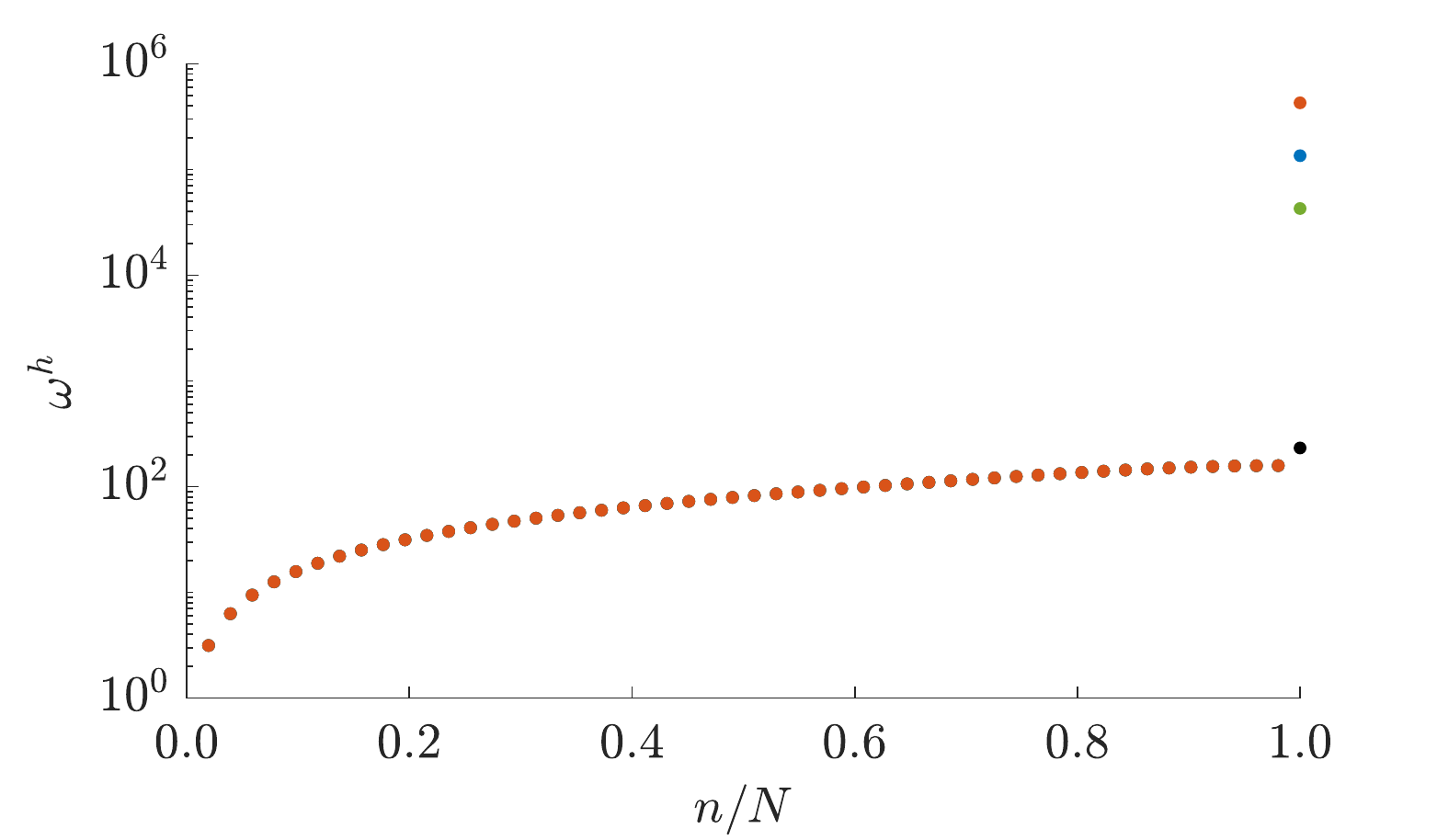} \\
    \vspace{0.2cm}
    \begin{tikzpicture}
		\filldraw[black,line width=1pt] (2,0) circle (2pt);
		\filldraw[black,line width=1pt] (2,0) node[right]{\footnotesize $\alpha=0$};
        \filldraw[green1,line width=1pt] (4,0) circle (2pt);
		\filldraw[green1,line width=1pt] (4,0) node[right]{\footnotesize $\alpha=1/h$};	
        \filldraw[blue1,line width=1pt] (6,0) circle (2pt);
		\filldraw[blue1,line width=1pt] (6,0) node[right]{\footnotesize $\alpha=10/h$};	
        \filldraw[red1,line width=1pt] (8,0) circle (2pt);
		\filldraw[red1,line width=1pt] (8,0) node[right]{\footnotesize $\alpha=100/h$};
	\end{tikzpicture}
    \caption{Discrete frequencies of a freely vibrating fixed bar, computed with $\boldsymbol f = \boldsymbol 0$ and different values of $\alpha$ chosen according to \cite{horger_penalty_2019}. We apply two patches of quadratic $C^1$ B-splines and discretize each patch with $25$ elements ($h=0.01$).}
    \label{fig:study_f0}
\end{figure}

We note that this case is equivalent to the formulation introduced in \cite{horger_penalty_2019}, excluding the penalty terms on the Neumann boundary and at cross-points. 
In Figure \ref{fig:study_f0}, we illustrate the discrete frequencies
of a fixed bar with values for $\alpha$ that are scaled with $1/h$ based on \cite{horger_penalty_2019}, without removing any outlier frequency. 
We note that 
$\alpha$ scaled with $1/h$ is not a consistent scaling factor of the matrices $\mat{K}$ and $\mat{K}_\Gamma$ (to the contrary, a consistent scaling factor would be $\alpha=h$).
As expected, we observe in Figure \ref{fig:study_f0} that the outlier frequency $\tilde{\omega}_{N}^h$ increases with increasing $\alpha$. \\

\begin{figure}[h!]
    \centering
    \includegraphics[width=0.5\textwidth]{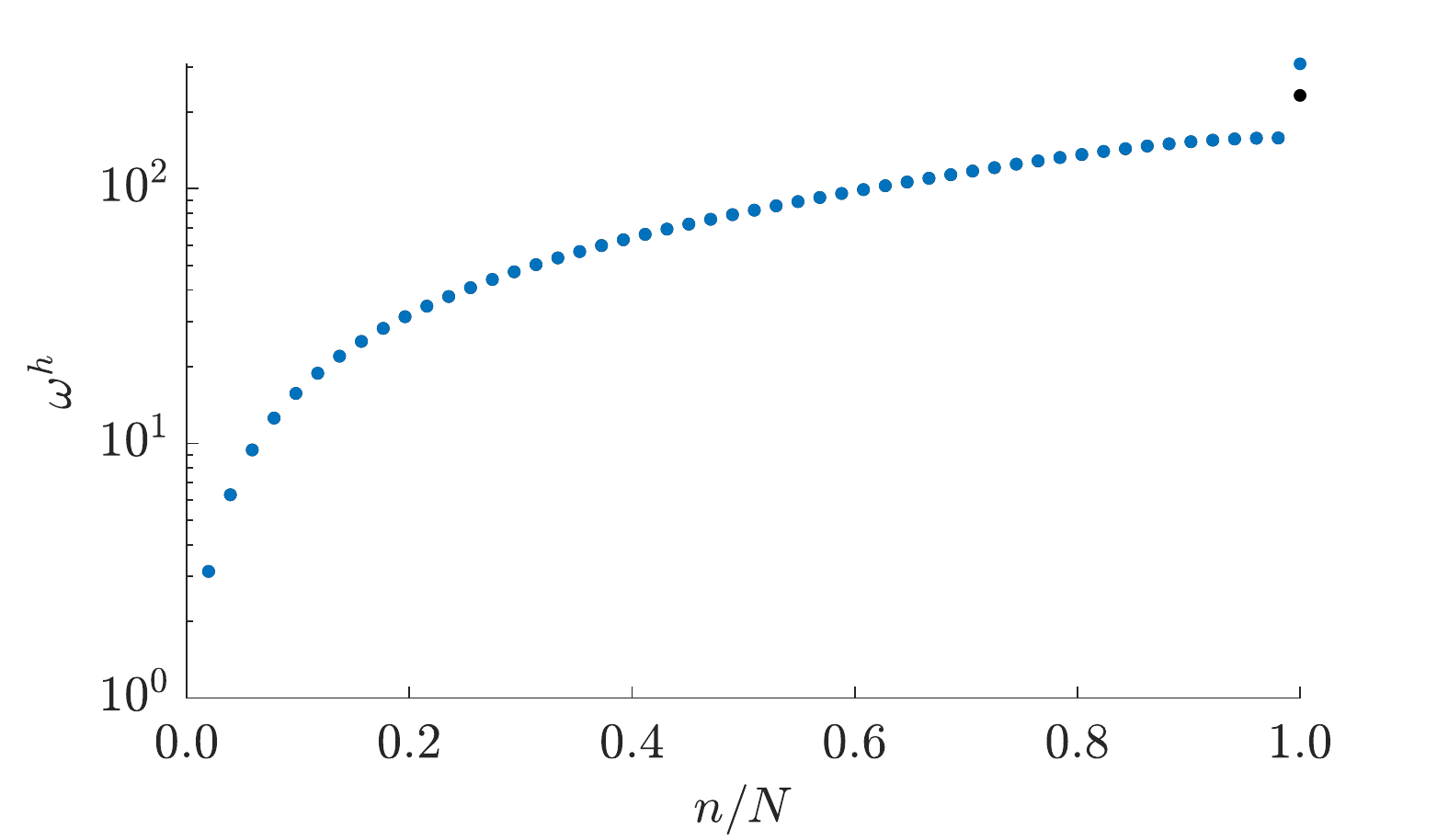} \\
    \vspace{0.2cm}
    \begin{tikzpicture}
		\filldraw[black,line width=1pt] (0,0) circle (2pt);
		\filldraw[black,line width=1pt] (0,0) node[right]{\footnotesize standard spectrum ($\alpha = \beta = 0$)};
        \filldraw[blue1,line width=1pt] (6,0) circle (2pt);
		\filldraw[blue1,line width=1pt] (6,0) node[right]{\footnotesize $f=0.5,\alpha=h$, $\beta = \alpha \, f / (\omega_N^h)^2$};	
	\end{tikzpicture}
    \caption{Discrete frequencies of a freely vibrating fixed bar, computed with $\boldsymbol f = \boldsymbol{0.5}$. We apply two patches of quadratic $C^1$ B-splines and discretize each patch with $25$ elements.}
    \label{fig:study_f_cases1}
\end{figure}
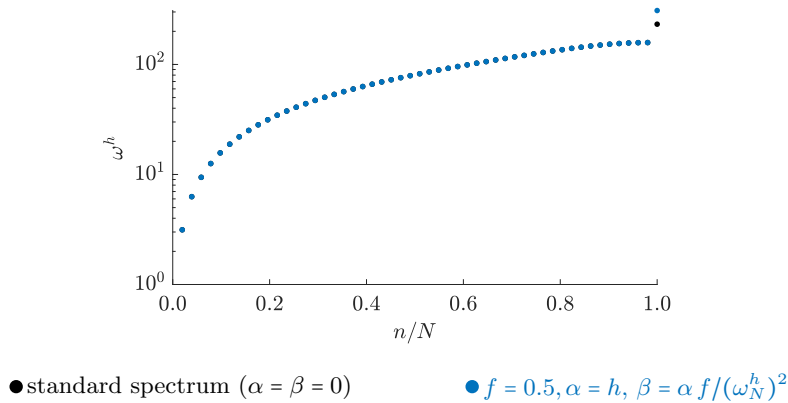

\begin{figure}[h!]
    \centering
    \includegraphics[width=0.5\textwidth]{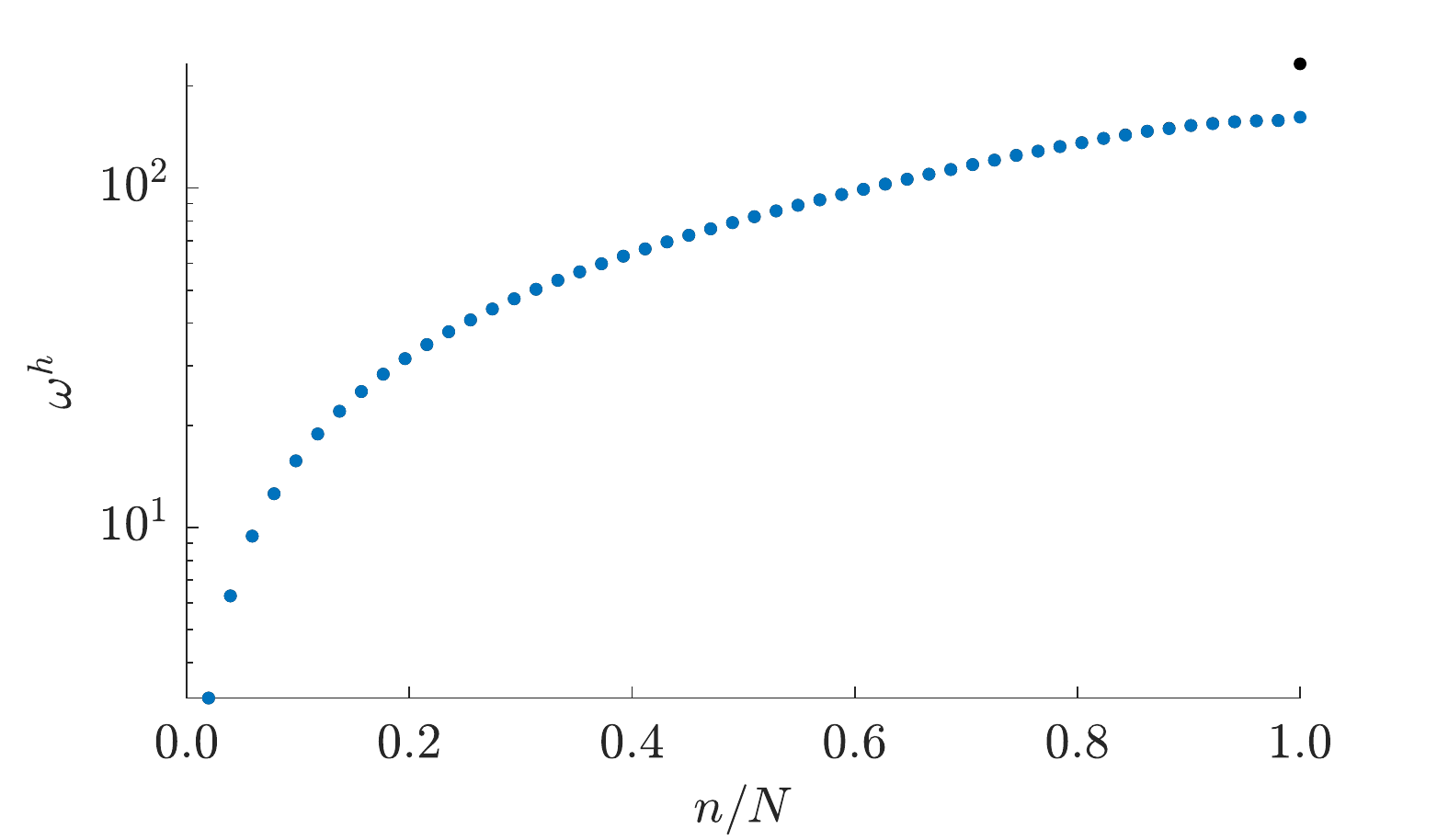} \\
    \vspace{0.2cm}
    \begin{tikzpicture}
		\filldraw[black,line width=1pt] (0,0) circle (2pt);
		\filldraw[black,line width=1pt] (0,0) node[right]{\footnotesize standard spectrum ($\alpha = \beta = 0$)};
        \filldraw[blue1,line width=1pt] (6,0) circle (2pt);
		\filldraw[blue1,line width=1pt] (6,0) node[right]{\footnotesize $f=2$};	
	\end{tikzpicture}
    \caption{Discrete frequencies of a freely vibrating fixed bar, computed with $\boldsymbol f \boldsymbol = \boldsymbol 2$. We apply two patches of quadratic $C^1$ B-splines and discretize each patch with $25$ elements.}
    \label{fig:f_study}
\end{figure}
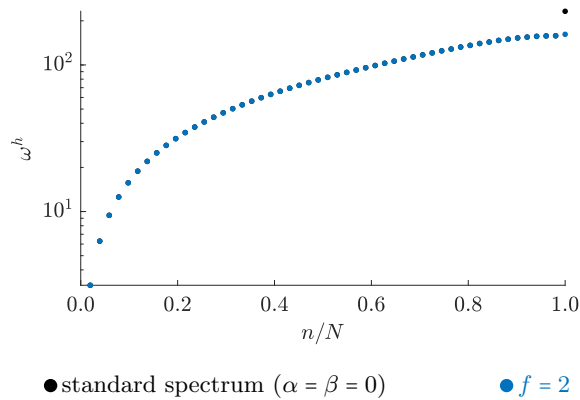

\begin{figure}[h!]
    \centering
    \includegraphics[width=0.5\textwidth]{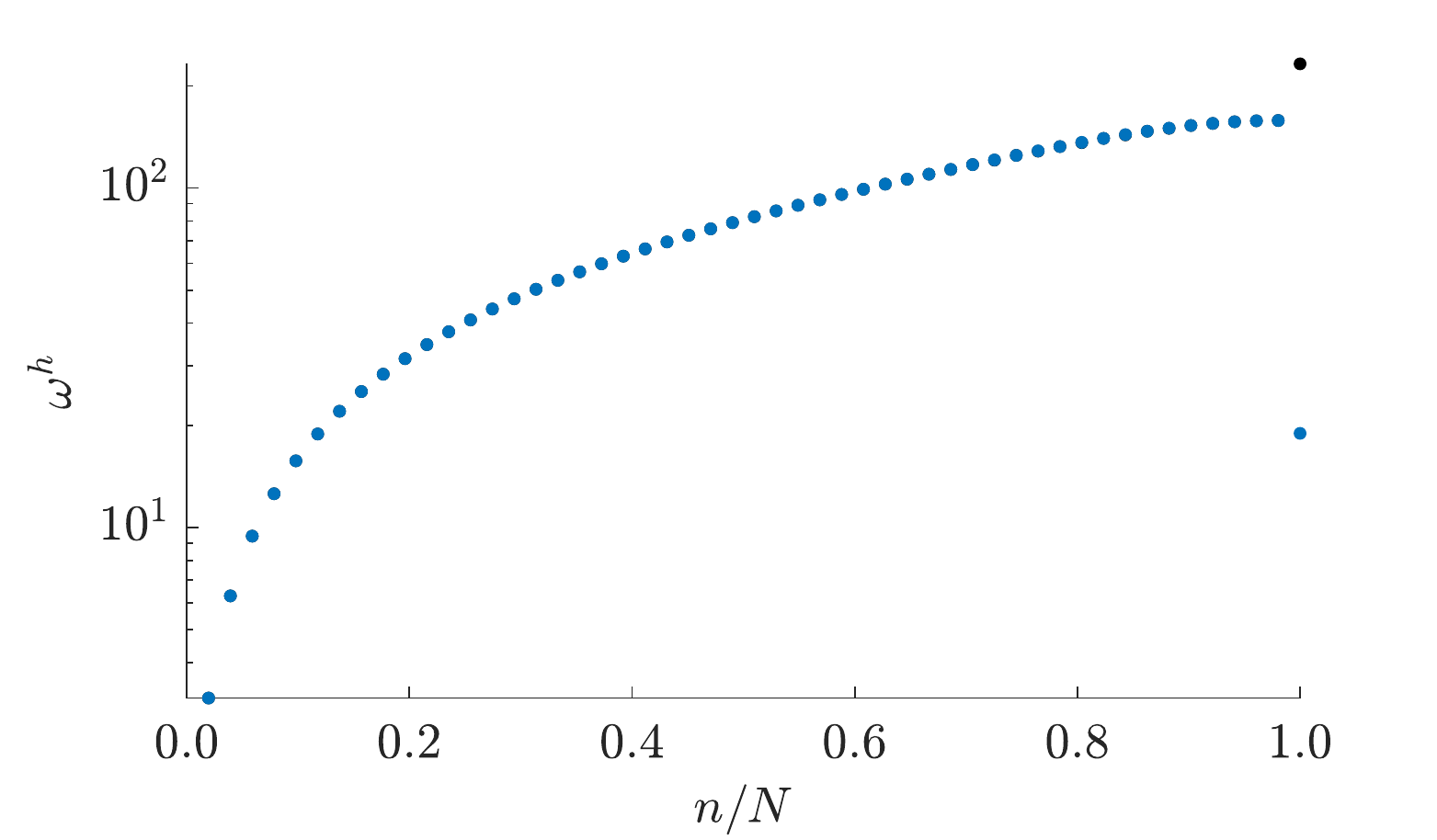} \\ 
    \vspace{0.2cm}
    \begin{tikzpicture}
		\filldraw[black,line width=1pt] (0,0) circle (2pt);
		\filldraw[black,line width=1pt] (0,0) node[right]{\footnotesize standard spectrum ($\alpha = \beta = 0$)};
        \filldraw[blue1,line width=1pt] (6,0) circle (2pt);
		\filldraw[blue1,line width=1pt] (6,0) node[right]{\footnotesize $\alpha = 0, \beta = h^3$};	
	\end{tikzpicture}
    \caption{Normalized frequencies of a freely vibrating fixed bar, computed with $\boldsymbol \alpha \boldsymbol = \boldsymbol 0, \boldsymbol \beta = h^3$. We apply two patches of quadratic $C^1$ B-splines and discretize each patch with $25$ elements.}
    \label{fig:study_f_cases2}
\end{figure}
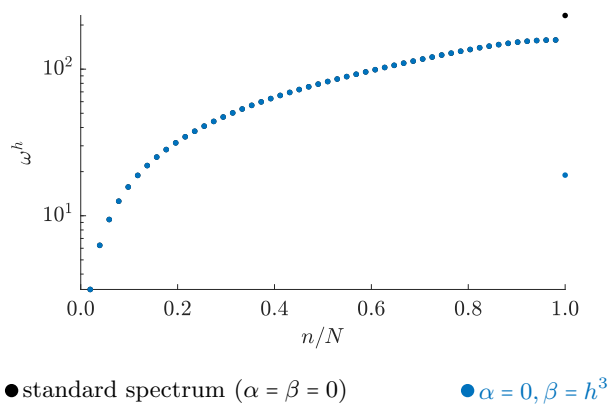

\begin{remark} \label{rm_ordering}
    \changedV{Instead of ordering the discrete frequencies in ascending order, as is typically done, the frequencies in Figures \ref{fig:study_f0} to \ref{fig:study_f_cases2} are ordered in such a way that the corresponding analytical frequencies are ascending. The corresponding pairs of discrete and analytical frequencies are identified by inspecting the corresponding mode shapes. In particular, we find the discrete mode that best fits, in the $L^2$ sense, a certain analytical mode.}
\end{remark}

\textbf{Case $\boldsymbol 0 \boldsymbol < \boldsymbol f \boldsymbol < \boldsymbol 1$}: The left-hand side of \eqref{eq:perturbed_alpha} is also positive when $\alpha$ is chosen larger than zero, i.e. $\tilde{\omega}_{max}^h > \omega^{h}_{max}$. Thus, a choice of $f$ in $(0,1)$ does not improve the spectrum. In Figure \ref{fig:study_f_cases1}, we illustrate the discrete frequencies of a fixed bar using $f=0.5$. We choose $\alpha=h$ and $\beta = \alpha \, f / (\omega_N^h)^2$, where $\omega_N^h$ is the $N^{\text{th}}$ frequency of the non-perturbed eigenvalue problem, i.e. the maximum outlier frequency. As expected, $\tilde{\omega}_{N}^h$ increases. \\

\textbf{Case $\boldsymbol f \boldsymbol > \boldsymbol 1$}: For a choice of $\alpha$ larger than zero, the left-hand side of \eqref{eq:perturbed_alpha} is negative, such that $\tilde{\omega}_{max}^h < \omega^{h}_{max}$. We thus observe that it is the mass matrix to which we should add the perturbation to improve the spectrum. 
In Figure \ref{fig:f_study}, we illustrate the discrete frequencies of a fixed bar computed \changed{with a value $f = 2$}. We compute $\alpha$ and $\beta$ \changed{using \eqref{eq:beta}-\eqref{eq:perturbed_alpha}}. 
\changed{As expected, the resulting frequency $\tilde{\omega}^h_N$ is reduced.} \\

\textbf{Case $\boldsymbol \alpha \boldsymbol = \boldsymbol 0$ and $\boldsymbol \beta \boldsymbol > \boldsymbol 0$}: Based on the previous observation, we may consider only adding the perturbation to the mass matrix, 
which can be interpreted as an approach of selective mass scaling \cite{olovsson_mass_scaling_2005,tkachuk_mass_scaling_2014}. 
In this case, values of $\beta$ that are scaled with $h^3$ are 
consistent scaling factors of the matrices $\mat{M}$ and $\mat{K}_\Gamma$. 
In Figure \ref{fig:study_f_cases2}, we illustrate the discrete frequencies of a fixed bar using a factor $\beta = h^3$. We observe that the outlier frequency \changedV{is also effectively reduced.}\\

\subsection{An iterative scheme based on first-order perturbed eigenvalue analysis}

\changedV{In this work, we focus on the parameter estimation in the case of perturbing both the stiffness and mass matrix that reduces the outlier frequency, i.e. the case of $\boldsymbol f \boldsymbol > \boldsymbol 1$.}

\subsubsection{\changed{Suppressing a single outlier frequency}}

Based on the first-order perturbation discussed above, we propose to use the following iterative procedure for approximating $\alpha$ and $\beta$ in the case $\boldsymbol f \boldsymbol > \boldsymbol 1$:
\begin{align}
    &\alpha^{(i)}  = \frac{\left(\omega^{h\,*}_{max}\right)^2 -\mat{\tilde{\eigenvec}}^{h\,(i-1)\,T}_{max} \mat{K} \mat{\tilde{\eigenvec}}^{h\,(i-1)}_{max} }{ \mat{\tilde{\eigenvec}}^{h\,(i-1)\,T}_{max} \, \mat{K}_\Gamma \, \mat{\tilde{\eigenvec}}^{h\,(i-1)}_{max} } \, \frac{1}{1 - f} \, , \label{eq:iterative_alpha} \\[0.2cm]  
    &\beta^{(i)} = f\, \frac{1}{\left(\omega^{h\,*}_{max}\right)^2} \alpha^{(i)} \, , \label{eq:iterative_beta}
\end{align}
where 
$\mat{\tilde{\eigenvec}}^{h\,(i-1)}_{max}$ correspond to the eigenvalue problem \eqref{eq:dgep_perturbed_1dp2} that is:
\begin{align}
    \left( \, \mat{K} + \alpha^{(i-1)} \, \mat{K}_{\Gamma}  \right) \, \mat{\tilde{\eigenvec}}_\noMode^{h\,(i-1)} \; = \; \left(\tilde{\omega}_\noMode^{h\,(i-1)} \right)^2 \, \left( \, \mat{M} + \beta^{(i-1)} \, \mat{K}_{\Gamma} \right) \, \mat{\tilde{\eigenvec}}_\noMode^{h\,(i-1)} \, ,
\end{align}
and $\omega^{h\,*}_{max}$ is the target maximum frequency.  
A practical choice of the target maximum frequency is a fraction of the unperturbed outlier frequency $\omega^{h}_{max}$. 
In this section, to illustrate the effectiveness of the proposed iterative scheme, we choose as target value $\omega^{h\,*}_{max}$ the maximum analytical frequency. For the studied benchmarks of bars and beams, the analytical solution is well-known. We find that only three or four iterations are required to obtain sufficiently converged values for $\alpha$ and $\beta$, i.e. these parameters converge within a small number of iterations. In Section \ref{sec:2d-study}, we propose an alternate strategy that avoids the need for a known analytical solution.

\subsubsection{Suppressing multiple outlier frequencies}

We \changed{now} discuss the parameter estimation of an example using cubic $C^2$ B-splines ($p=3$) \changedV{with $C^0$ patch continuity.}
The resulting spectrum of a two-patches-discretization ($\npa=2$) consists of two interior outliers (see Table \ref{tab:number_of_outliers_1Dbar}). 
The perturbed eigenvalue problem \eqref{eq:dgep_perturbed} then simplifies to:
\begin{align}
    \left[ \mat{K} + \alpha^1 \mat{K}_\Gamma^1 + \alpha^2 \mat{K}_\Gamma^2 \right] \, \mat{\tilde{\eigenvec}}_{n}^h = \left(\tilde{\omega}^h_{n}\right)^2 \, \left[ \mat{M} + \beta^1 \mat{K}_\Gamma^1 + \beta^2 \mat{K}_\Gamma^2 \right] \, \mat{\tilde{\eigenvec}}_{n}^h \, . \label{eq:dgep_perturbed_1dp3}
\end{align}

To iteratively estimate the parameters $\alpha^l$ and $\beta^l$, $l=1,2$, in each iteration, we first identify the two outlier modes corresponding to each continuity constraint of \eqref{eq:continuity_constraint} as follows:
\begin{align}
    n^l = \argmax_{N^l} \left(\mat{\tilde{\eigenvec}}_{n}^{h\,T} \, \mat{K}_\Gamma^l \, \mat{\tilde{\eigenvec}}_{n}^{h} \right) \; , \quad N^l = \left\{ 1, \ldots, N \right\} \backslash \left\{n^1, \ldots, n^{l-1}\right\} \, .
	\label{eq:maximum}
\end{align}
Inserting these outlier modes in \eqref{eq:dgep_perturbed_1dp3}, aiming at the corresponding target frequencies, and performing the steps \eqref{step1}-\eqref{step2}, as well as expressing $\beta^l$ in terms of $\alpha^l$ and the target frequencies, we obtain the following system of equations:
\begin{subequations}\label{eq:parameter_1dp3}
    \begin{align}
        & \alpha^1 \, \left(1 - f^1 \right) \,\mat{\tilde{\eigenvec}}_{n^1}^{h\,T} \, \mat{K}_\Gamma^1 \, \mat{\tilde{\eigenvec}}_{n^1}^{h} + \alpha^2 \, \left(1 - f^2 \frac{\omega^{h\,*}_{n^1}}{\omega^{h\,*}_{n^2}}\right) \, \mat{\tilde{\eigenvec}}_{n^1}^{h\,T} \, \mat{K}_\Gamma^2 \, \mat{\tilde{\eigenvec}}_{n^1}^{h} = \left(\tilde{\omega}^{h\,*}_{n^1}\right)^2 - \mat{\tilde{\eigenvec}}_{n^1}^{h\,T} \, \mat{K} \, \mat{\tilde{\eigenvec}}_{n^1}^{h} \, , \\
        & \alpha^1 \, \left(1 - f^1 \frac{\omega^{h\,*}_{n^2}}{\omega^{h\,*}_{n^1}} \right) \, \mat{\tilde{\eigenvec}}_{n^2}^{h\,T} \, \mat{K}_\Gamma^1 \, \mat{\tilde{\eigenvec}}_{n^2}^{h} + \alpha^2 \, \left(1 - f^2 \right) \, \mat{\tilde{\eigenvec}}_{n^2}^{h\,T} \, \mat{K}_\Gamma^2 \, \mat{\tilde{\eigenvec}}_{n^2}^{h} = \left(\tilde{\omega}^{h\,*}_{n^2}\right)^2 - \mat{\tilde{\eigenvec}}_{n^2}^{h\,T} \, \mat{K} \, \mat{\tilde{\eigenvec}}_{n^2}^{h} \, ,
    \end{align}
\end{subequations}
where 
\begin{align}
    \beta^1 = f^1 \frac{1}{\left(\omega^{h\,*}_{n^1} \right)^2} \alpha^1 \; , \qquad \beta^2 = f^2 \frac{1}{\left(\omega^{h\,*}_{n^2} \right)^2} \alpha^2 \, .
\end{align}
Solving this system of equations, we obtain $\alpha^l$ and $\beta^l$ at each iteration. 
\changed{We find that only four or five iterations are required to obtain sufficiently converged values for $\alpha^l$, $l=1,2$.} 

In general, the resulting system of equations consists of $(p-1)$ equations and $(p-1)$ unknown parameters $\alpha^l$. Given a choice of $(p-1)$ factors $f^l$, we obtain $\beta^l$ in terms of $f^l$, $\alpha^l$, and the target maximum frequency, see \eqref{eq:iterative_beta}. 
We note that the proposed iterative scheme requires $(p-1)$ target maximum frequencies and an identification of $(p-1)$ outlier modes corresponding to the continuity constraint \eqref{eq:continuity_constraint} at each iteration.

\subsection{Spectral analysis of a second-order problem}
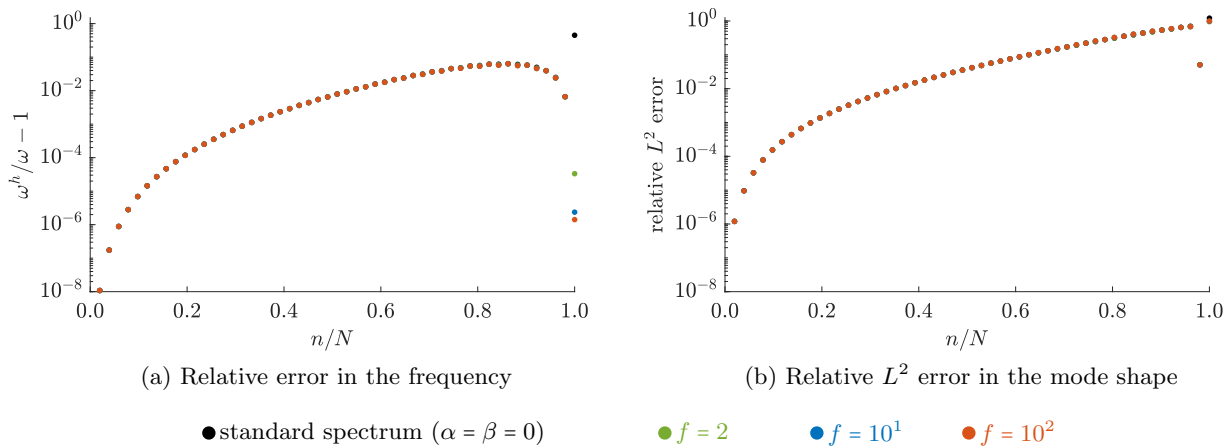
\begin{figure}[h!]
    \centering
    \subfloat[Relative error in the frequency]{{\includegraphics[width=0.5\textwidth]{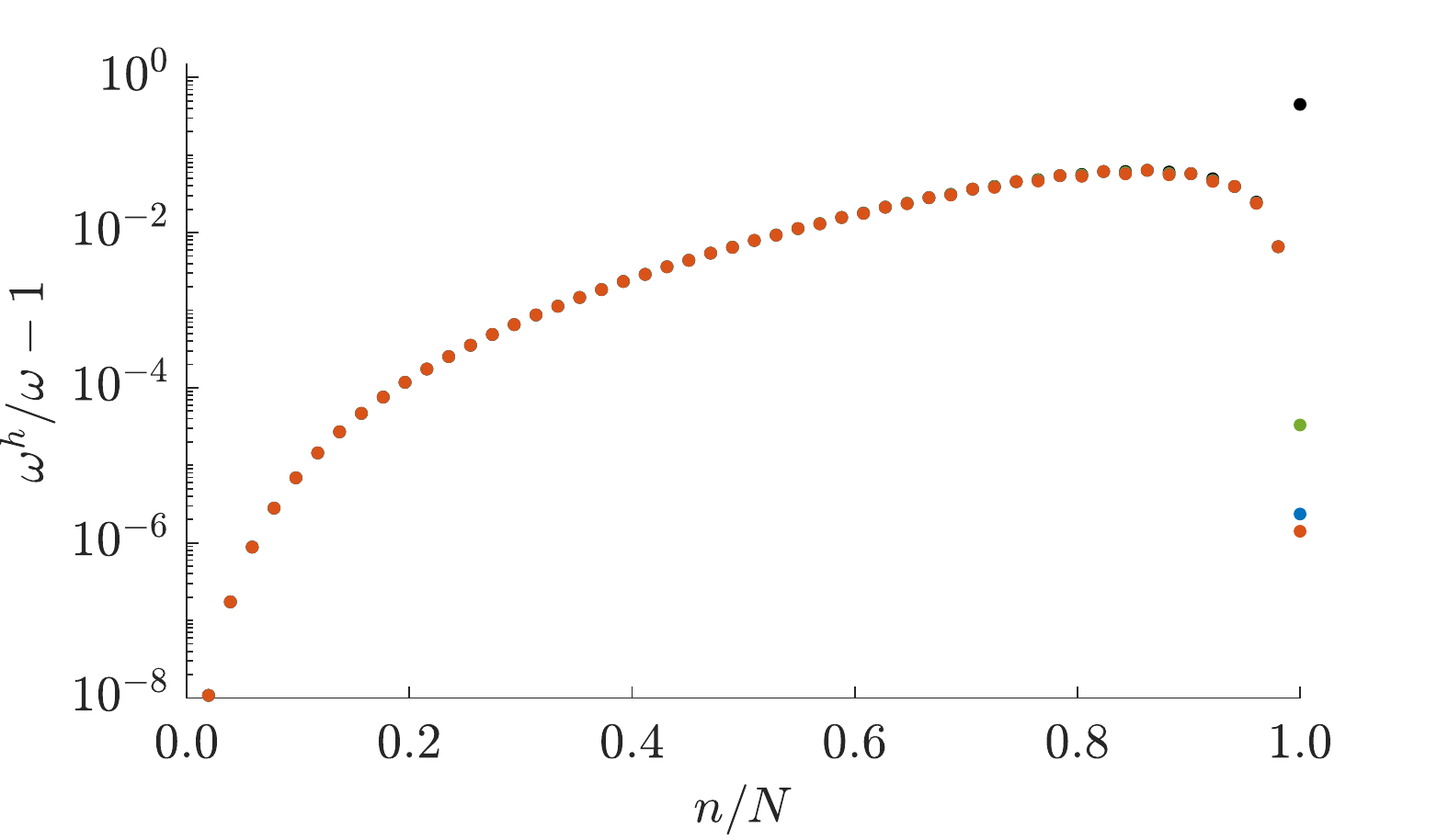} }}
    \subfloat[Relative $L^2$ error in the mode shape]{{\includegraphics[width=0.5\textwidth]{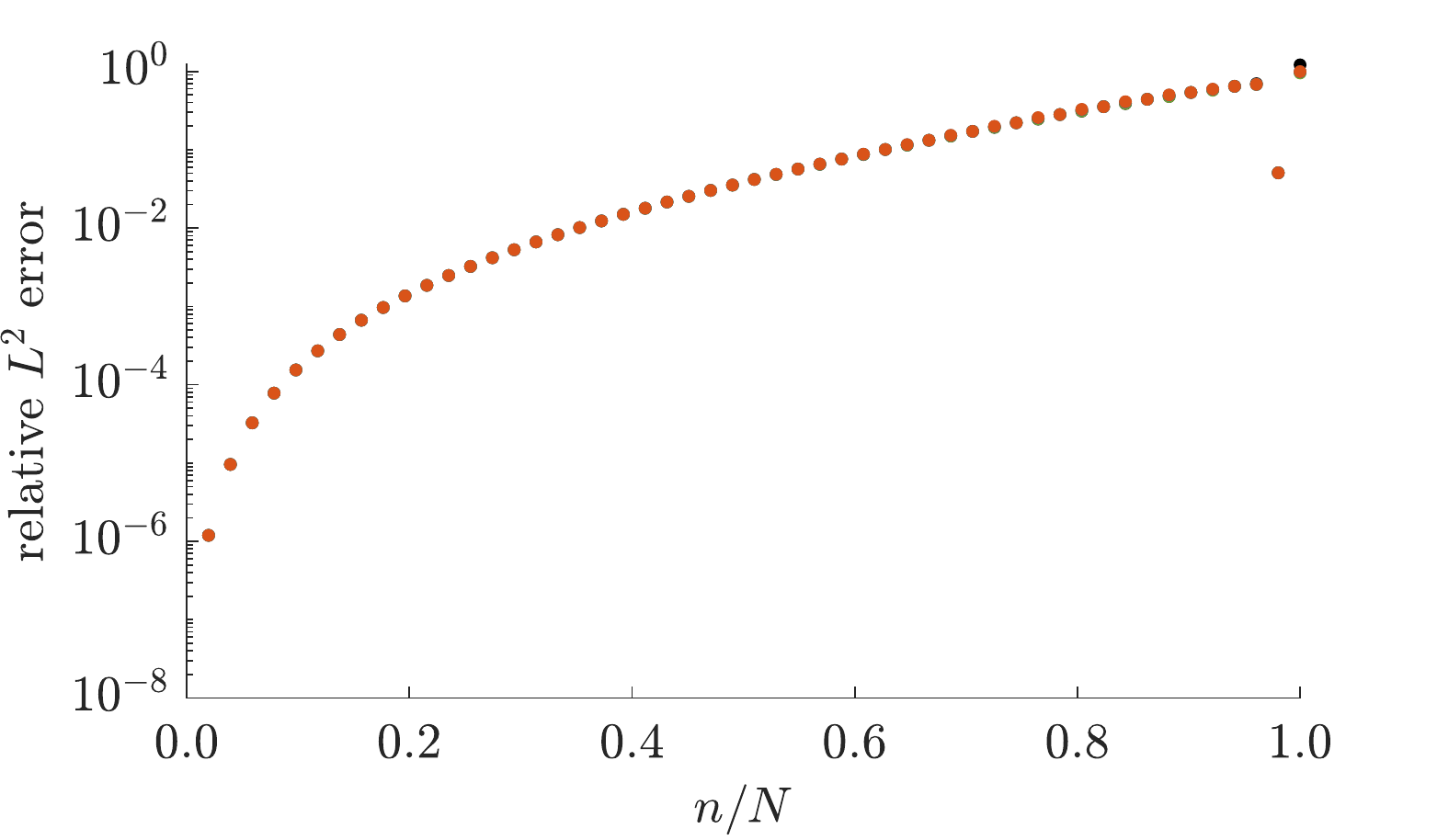} }}
    \vspace{0.2cm}
    \begin{tikzpicture}
		\filldraw[black,line width=1pt] (0,0) circle (2pt);
		\filldraw[black,line width=1pt] (0,0) node[right]{\footnotesize standard spectrum ($\alpha = \beta = 0$)};
        \filldraw[green1,line width=1pt] (6,0) circle (2pt);
		\filldraw[green1,line width=1pt] (6,0) node[right]{\footnotesize $f=2$};	
		\filldraw[blue1,line width=1pt] (8,0) circle (2pt);
		\filldraw[blue1,line width=1pt] (8,0) node[right]{\footnotesize $f=10^1$};	
        \filldraw[red1,line width=1pt] (10,0) circle (2pt);
		\filldraw[red1,line width=1pt] (10,0) node[right]{\footnotesize $f=10^2$};
	\end{tikzpicture}
    \caption{Relative error in frequencies and $L^2$ error in the mode shapes of a freely vibrating fixed bar, computed with different values $\boldsymbol f \boldsymbol > \boldsymbol 1$. We apply two patches of quadratic $C^1$ B-splines and discretize each patch with $25$ elements.}
    \label{fig:f_study2}
\end{figure}

The proposed iterative parameter estimation still requires a choice of the scaling factor $f > 1$ between $\alpha$ and $\beta$.
\changed{Figure \ref{fig:f_study2} illustrates the relative error in the frequency, $|\omega^h_\noMode / \omega_\noMode - 1|$, and the relative $L^2$ error in the mode, $\|\eigenvec^h_\noMode - \eigenvec_\noMode\|_{L^2} \, / \, \|\eigenvec_\noMode\|_{L^2}$, of a fixed bar discretized with two patches of quadratic $C^1$ B-splines. We order these results in the same way as described in Remark \ref{rm_ordering}.} \changedV{We observe that optimal accuracy is preserved in the remainder of the spectrum and modes for all values $f>1$. In all numerical studies that follow we choose a factor $f=2$.}

\begin{figure}[h!]
    \centering
    \subfloat[$25$ elements per patch]{{\includegraphics[width=0.5\textwidth]{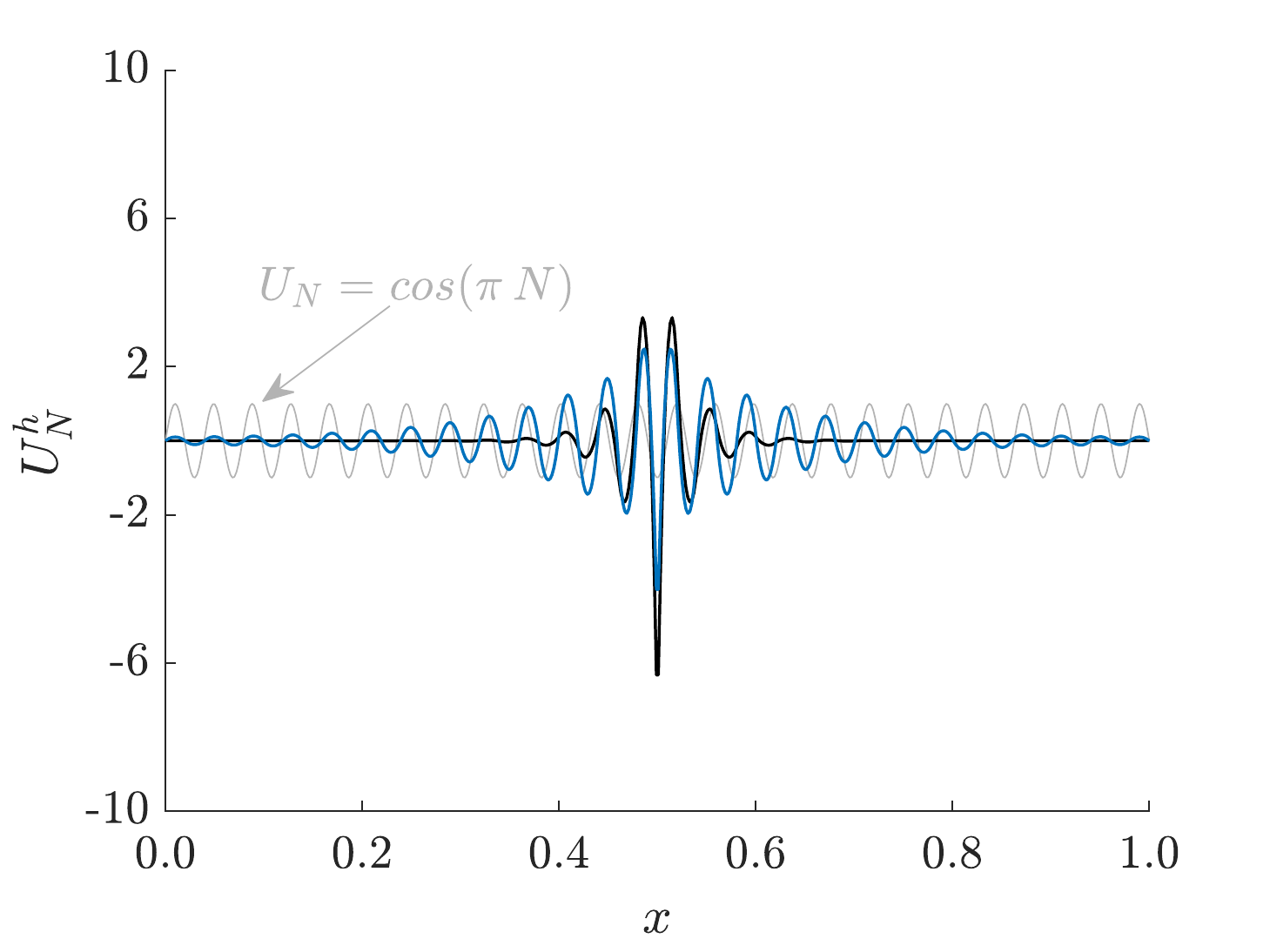}}}
    \subfloat[$50$ elements per patch]{{\includegraphics[width=0.5\textwidth]{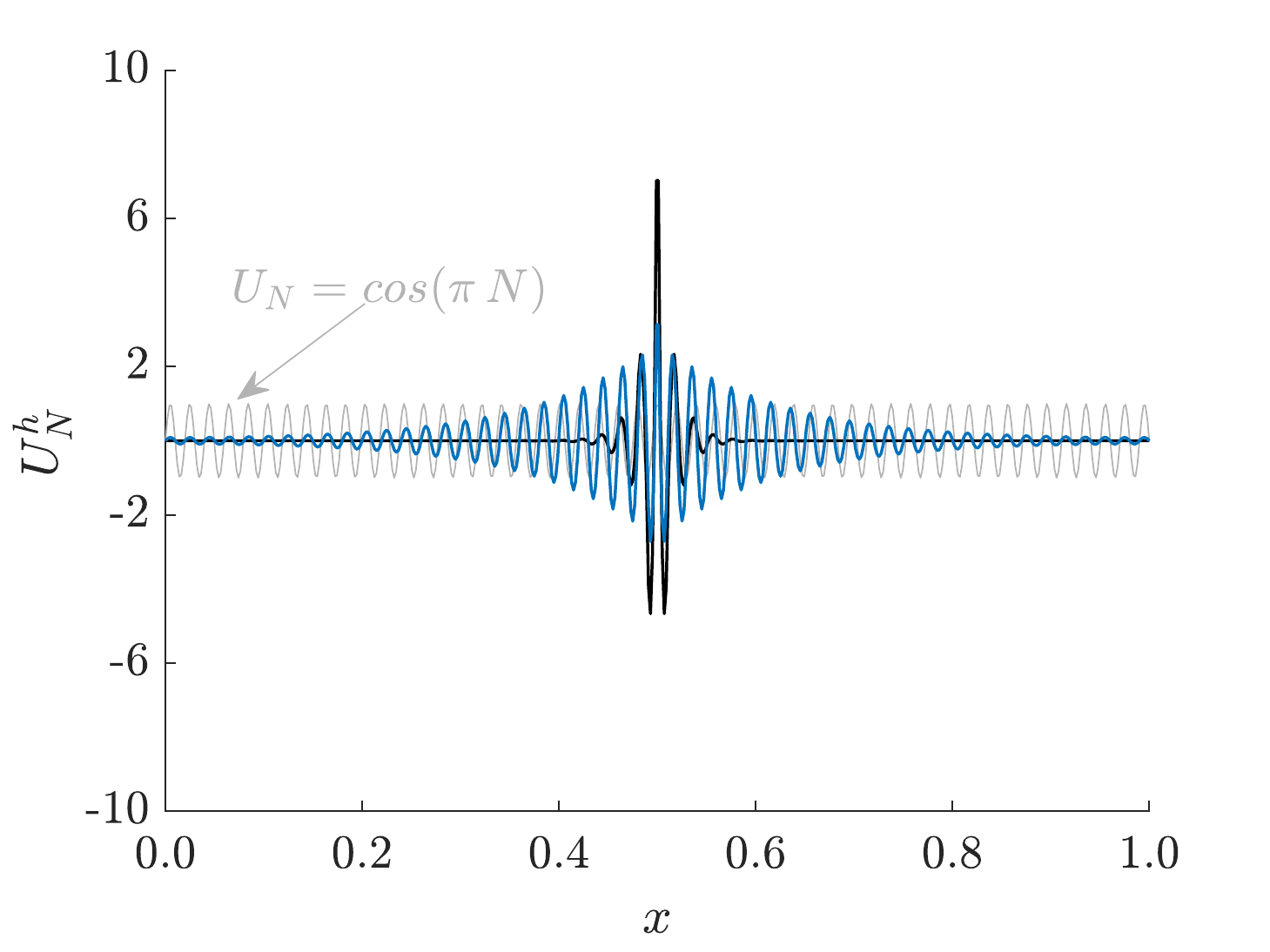}}}
    \vspace{0.2cm}
    \begin{tikzpicture}
		\filldraw[black,line width=1pt, solid] (4.0,0) -- (4.6,0);
		\filldraw[black,line width=1pt] (4.6,0) node[right]{\footnotesize standard spectrum ($\alpha = \beta = 0$)};
        \filldraw[blue1,line width=1pt] (10.0,0) -- (10.6,0);
		\filldraw[blue1,line width=1pt] (10.6,0) node[right]{\footnotesize $f=2$, $\alpha,\beta$ computed with \eqref{eq:iterative_alpha}-\eqref{eq:iterative_beta}};		
	\end{tikzpicture}
    \caption{Outlier mode of a freely vibrating fixed bar, computed with $2$ patches of quadratic $C^1$ B-splines and normalized such that $\left\|\eigenvec^h_\noMode\right\|_{L^2} = \left\|\eigenvec_\noMode\right\|_{L^2}$.}
    \label{fig:C0vsPerturbed_mode}
\end{figure}
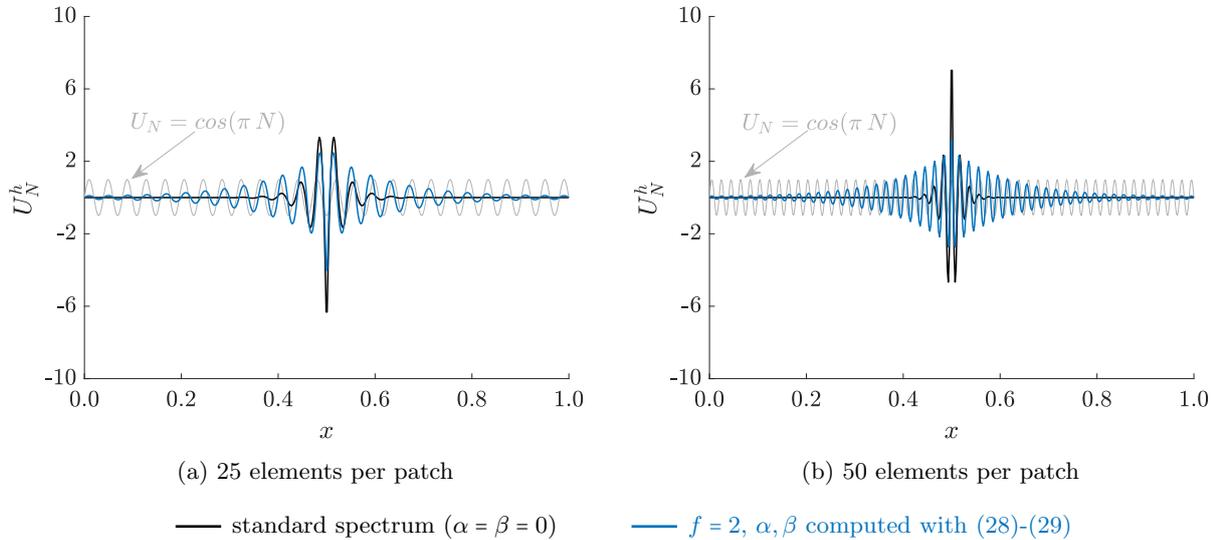

The corresponding outlier mode, as expected, remains spurious as illustrated in Figure \ref{fig:C0vsPerturbed_mode},
since the addition of perturbations improves the eigenvalue spectrum, but does not remove spurious outlier modes. 
Nevertheless, the introduced perturbation reduces the error in the outlier mode, as demonstrated in Figure \ref{fig:f_study2}b, since the perturbed outlier mode approximates the analytical solution better than the unperturbed mode, as shown in Figure \ref{fig:C0vsPerturbed_mode}. 

\begin{figure}[h!]
    \centering
    \subfloat[$p=2$]{{\includegraphics[width=0.5\textwidth]{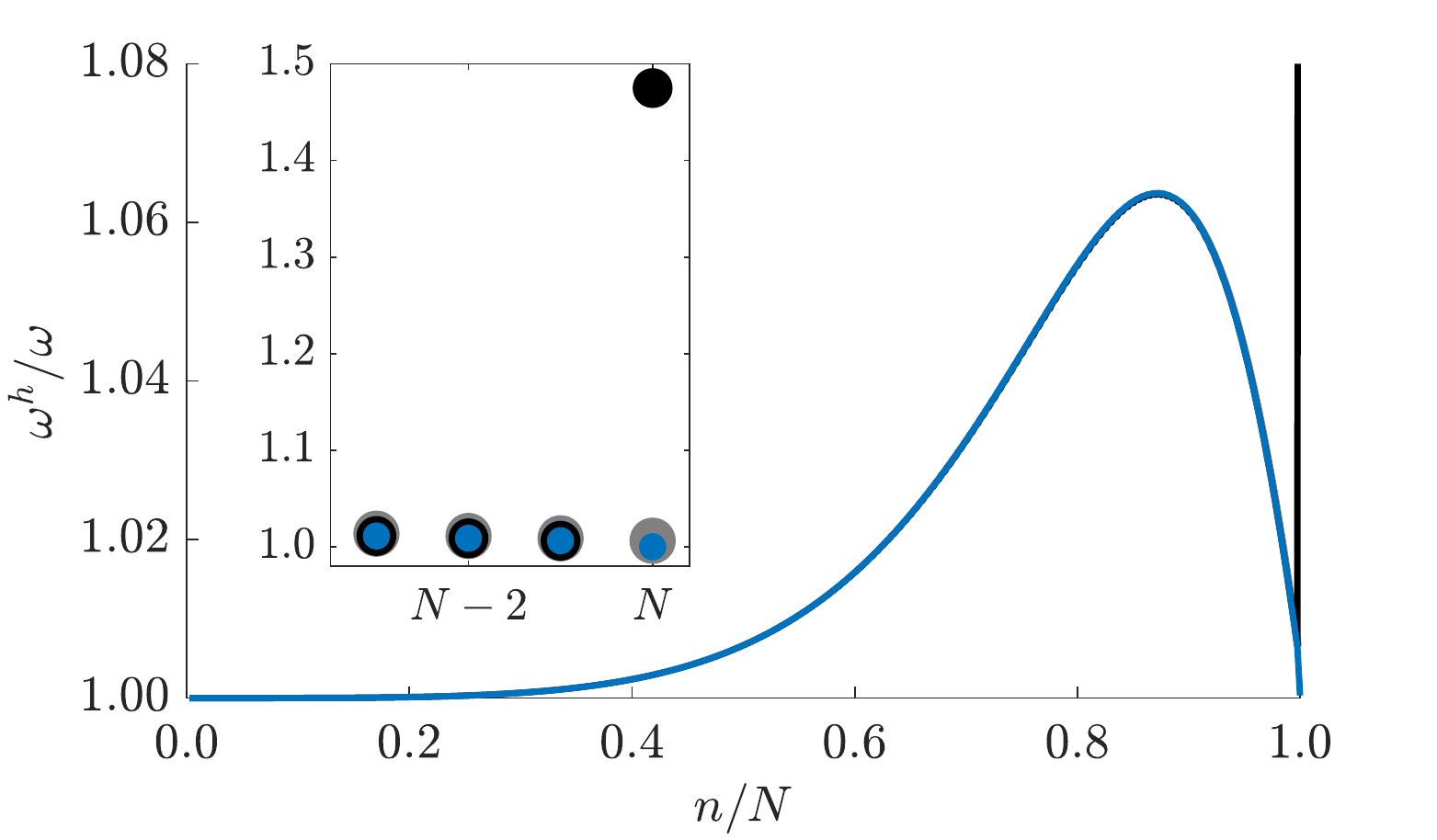} }}
    \subfloat[$p=3$]{{\includegraphics[width=0.5\textwidth]{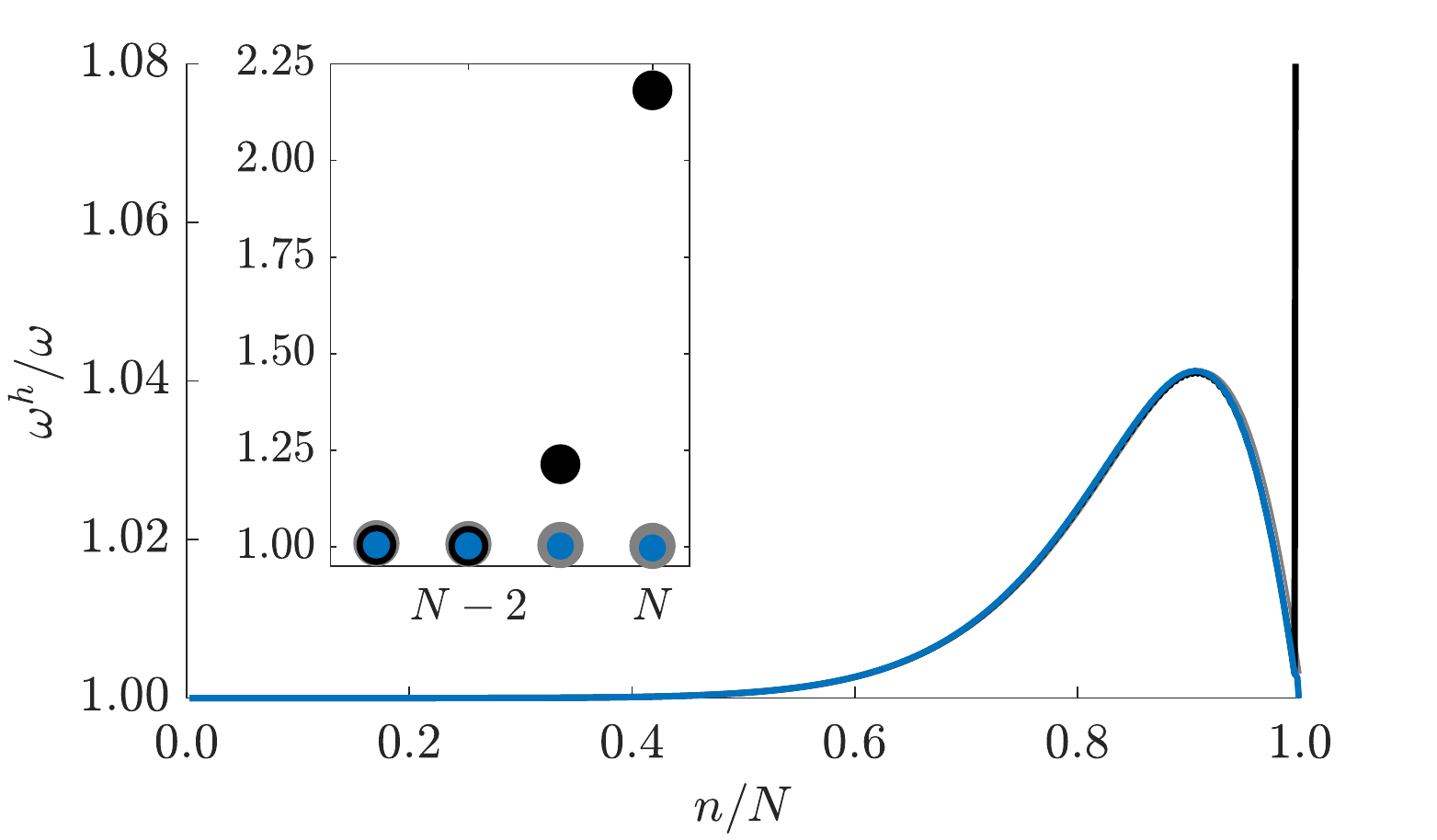} }}

    \subfloat[$p=4$]{{\includegraphics[width=0.5\textwidth]{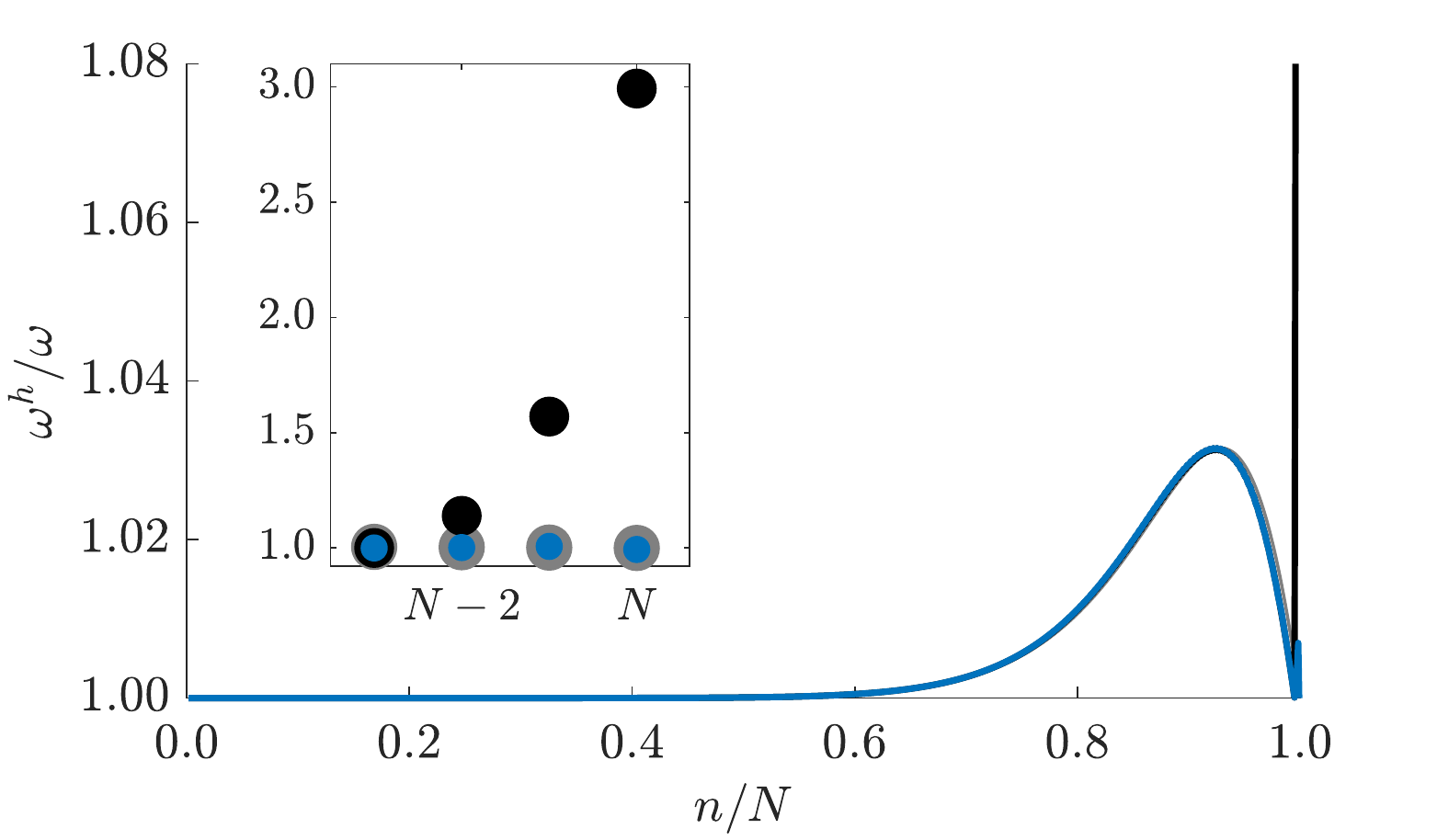} }}
    \subfloat[$p=5$]{{\includegraphics[width=0.5\textwidth]{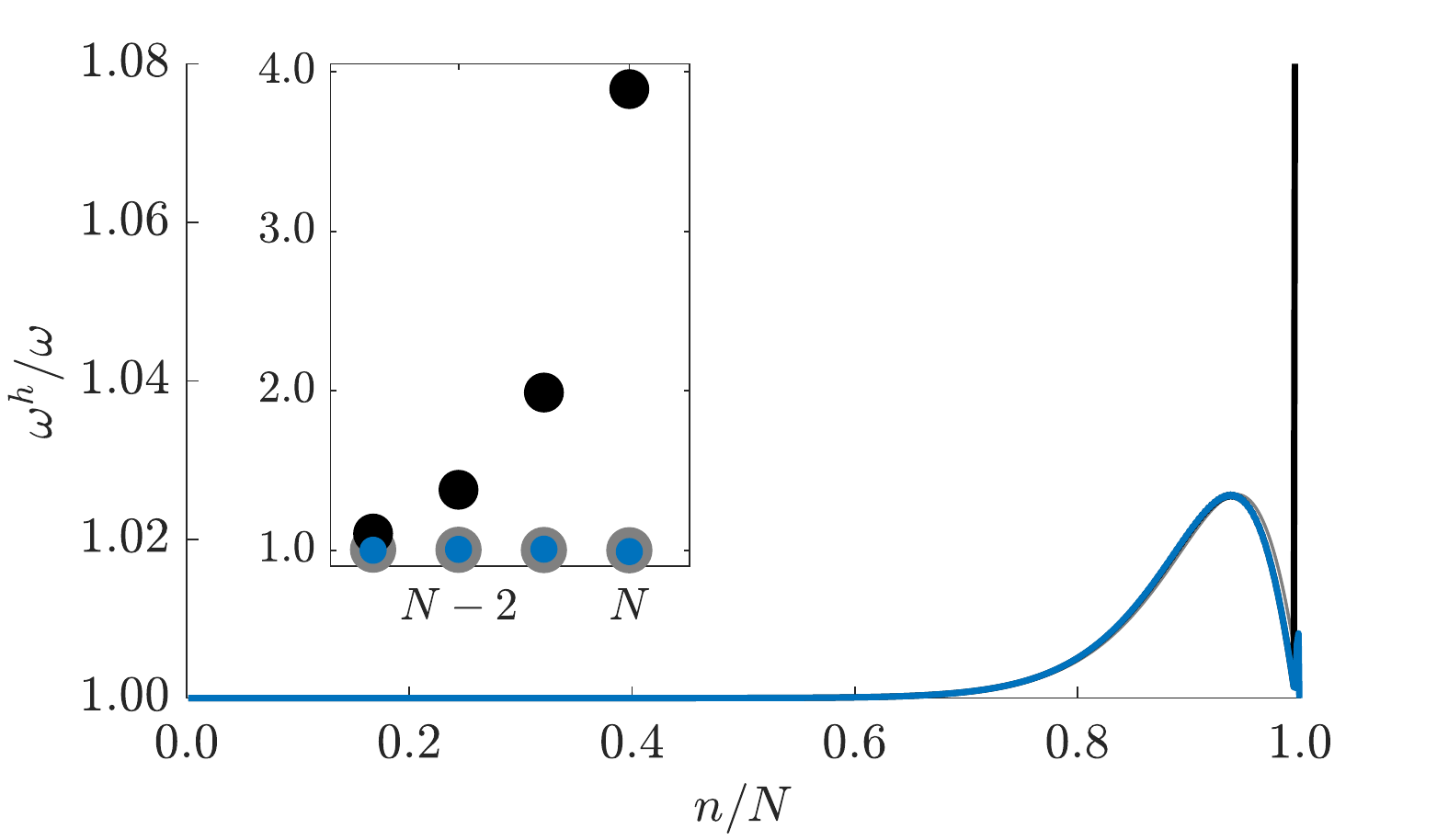} }}
    \vspace{0.2cm}
    \begin{tikzpicture}
    \filldraw[grey1,line width=1pt] (0,0) circle (2pt);
    \filldraw[grey1,line width=1pt] (0,0) node[right]{\footnotesize single-patch};
    \filldraw[black,line width=1pt] (3,0) circle (2pt);
    \filldraw[black,line width=1pt] (3,0) node[right]{\footnotesize multipatch, standard spectrum};
    \filldraw[blue1,line width=1pt] (9,0) circle (2pt);
    \filldraw[blue1,line width=1pt] (9,0) node[right]{\footnotesize multipatch, improved spectrum};
\end{tikzpicture}
    \caption{Normalized frequencies of a freely vibrating \textbf{fixed bar}, computed with \textbf{two patches} of $C^{p-1}$ B-splines and discretized with \textbf{$500$ elements}.}
    \label{fig:normalized_freq_2nd_1d}
\end{figure}

\begin{figure}[h!]
    \centering
    \subfloat[$p=2$]{{\includegraphics[width=0.47\textwidth]{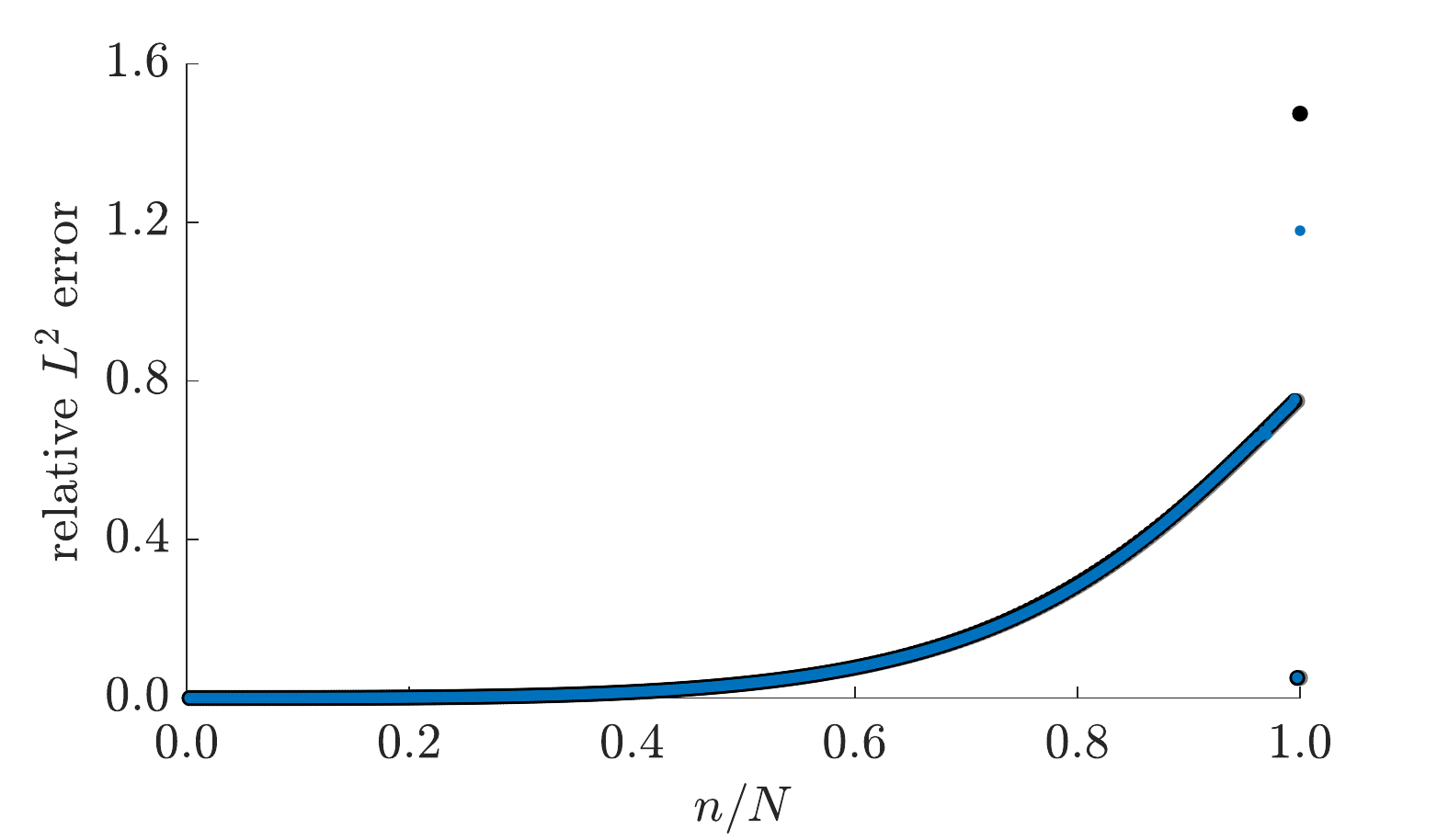} }}
    \subfloat[$p=3$]{{\includegraphics[width=0.47\textwidth]{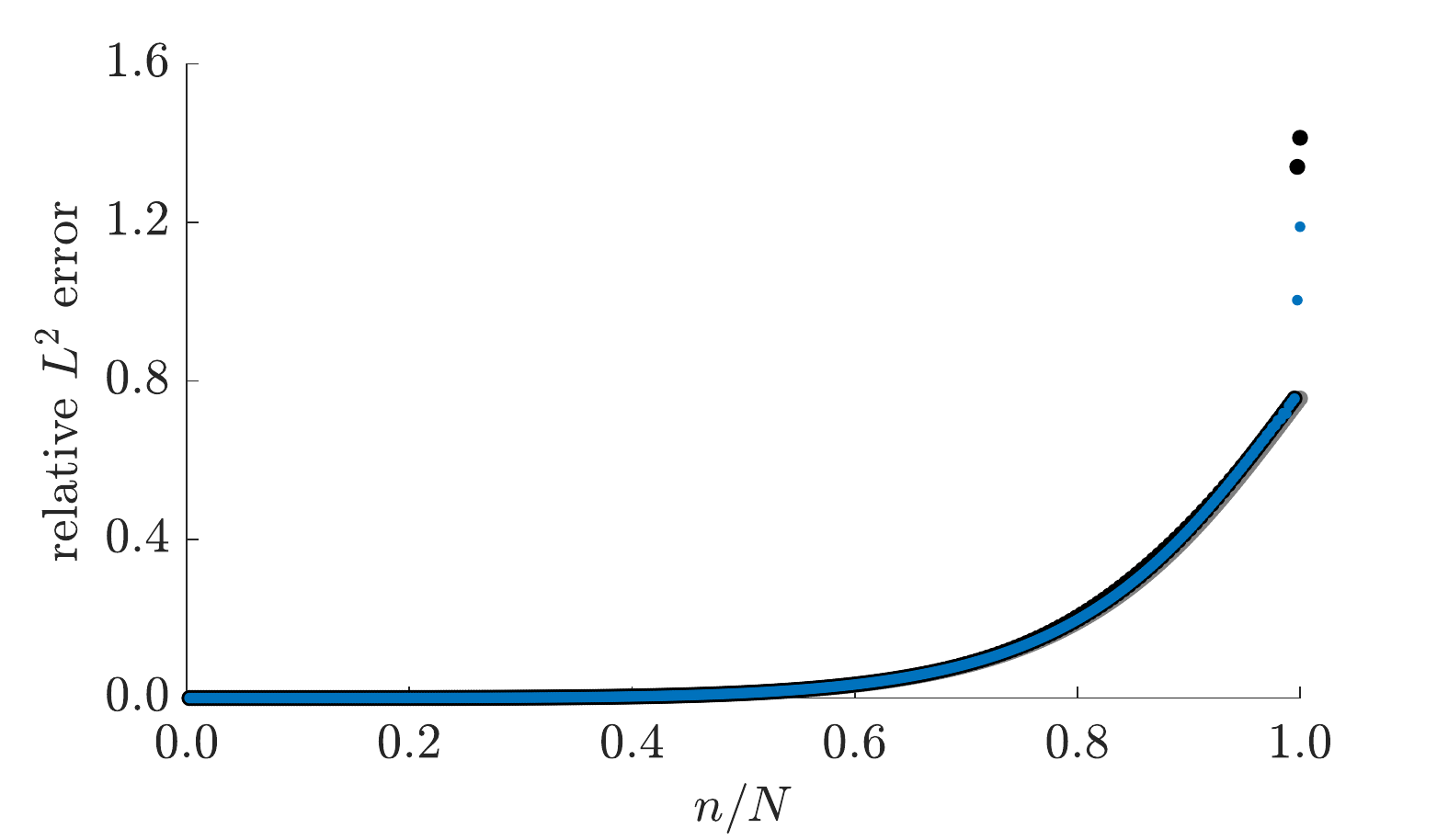} }}

    \subfloat[$p=4$]{{\includegraphics[width=0.47\textwidth]{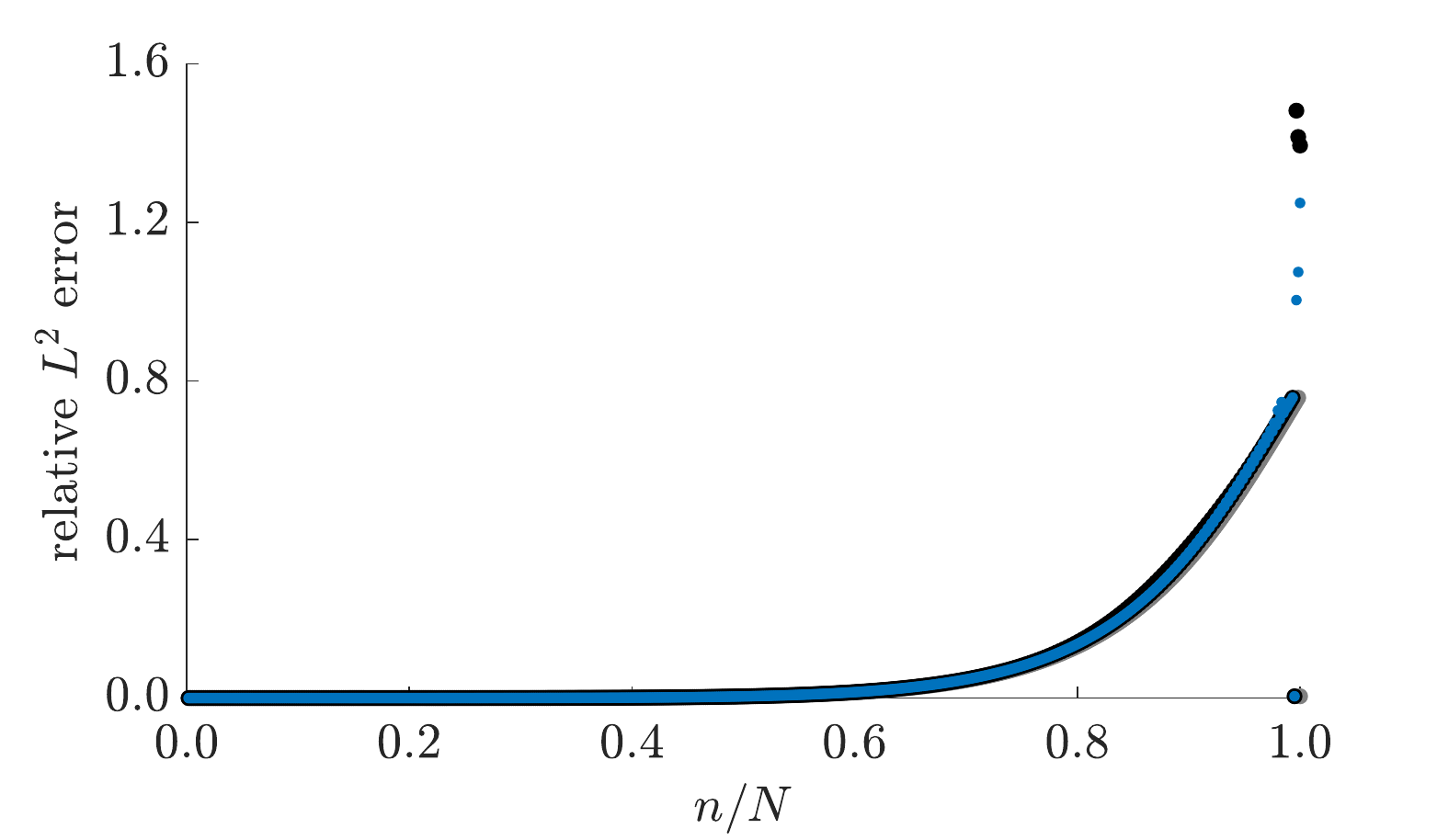} }}
    \subfloat[$p=5$]{{\includegraphics[width=0.47\textwidth]{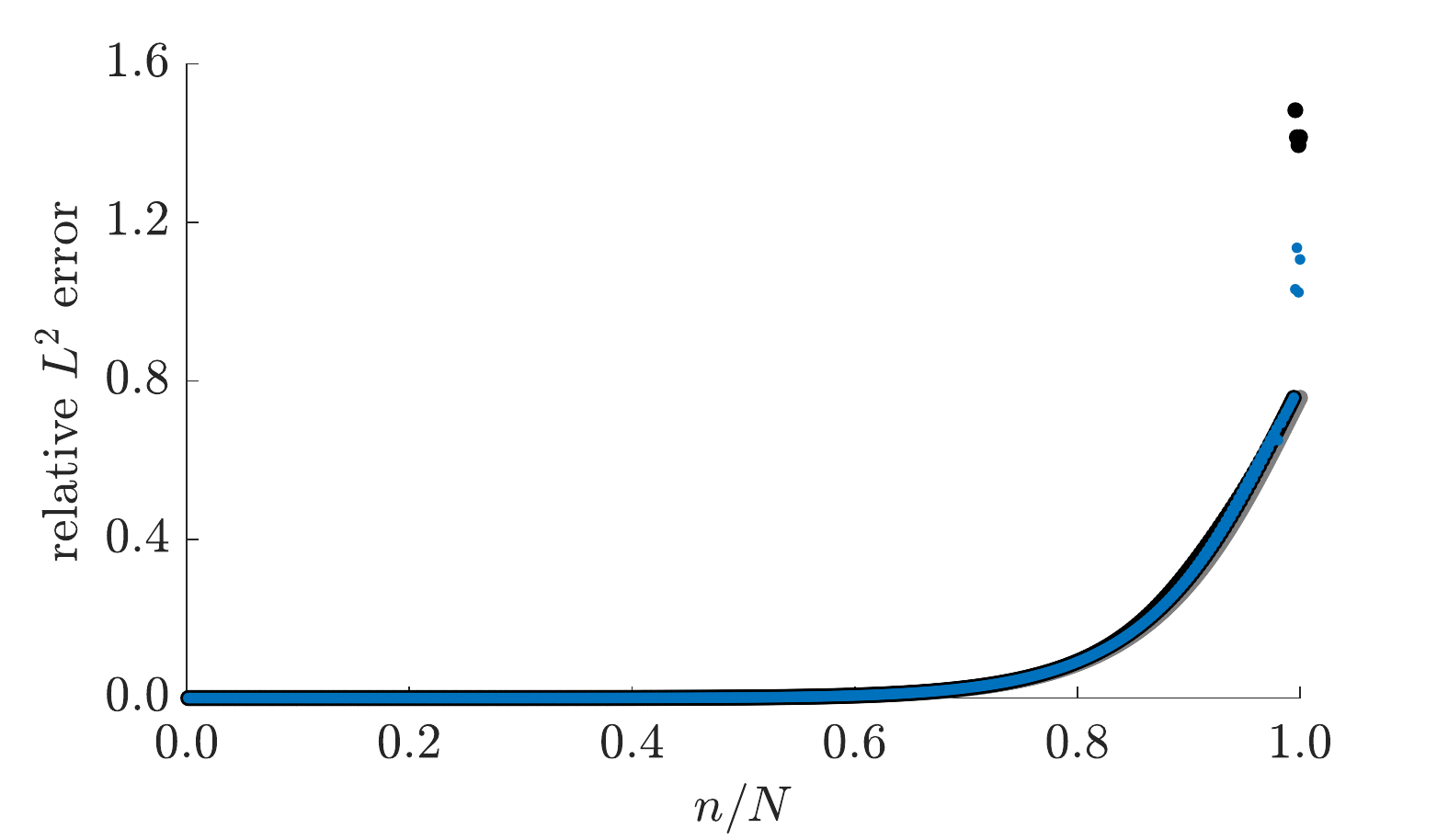} }}
    \vspace{0.2cm}
    \begin{tikzpicture}
    \filldraw[grey1,line width=1pt] (0,0) circle (2pt);
    \filldraw[grey1,line width=1pt] (0,0) node[right]{\footnotesize single-patch};
    \filldraw[black,line width=1pt] (3,0) circle (2pt);
    \filldraw[black,line width=1pt] (3,0) node[right]{\footnotesize multipatch, standard spectrum};
    \filldraw[blue1,line width=1pt] (9,0) circle (2pt);
    \filldraw[blue1,line width=1pt] (9,0) node[right]{\footnotesize multipatch, improved spectrum};
\end{tikzpicture}
    \caption{$L^2$ errors in the mode shapes of a freely vibrating \textbf{fixed bar}, computed with \textbf{two patches} of $C^{p-1}$ B-splines and discretized \textbf{with $500$ elements}.}
    \label{fig:mode_error_2nd_1d}

    \vspace{0.4cm}

    \subfloat[Normalized frequency error]{{\includegraphics[width=0.5\textwidth]{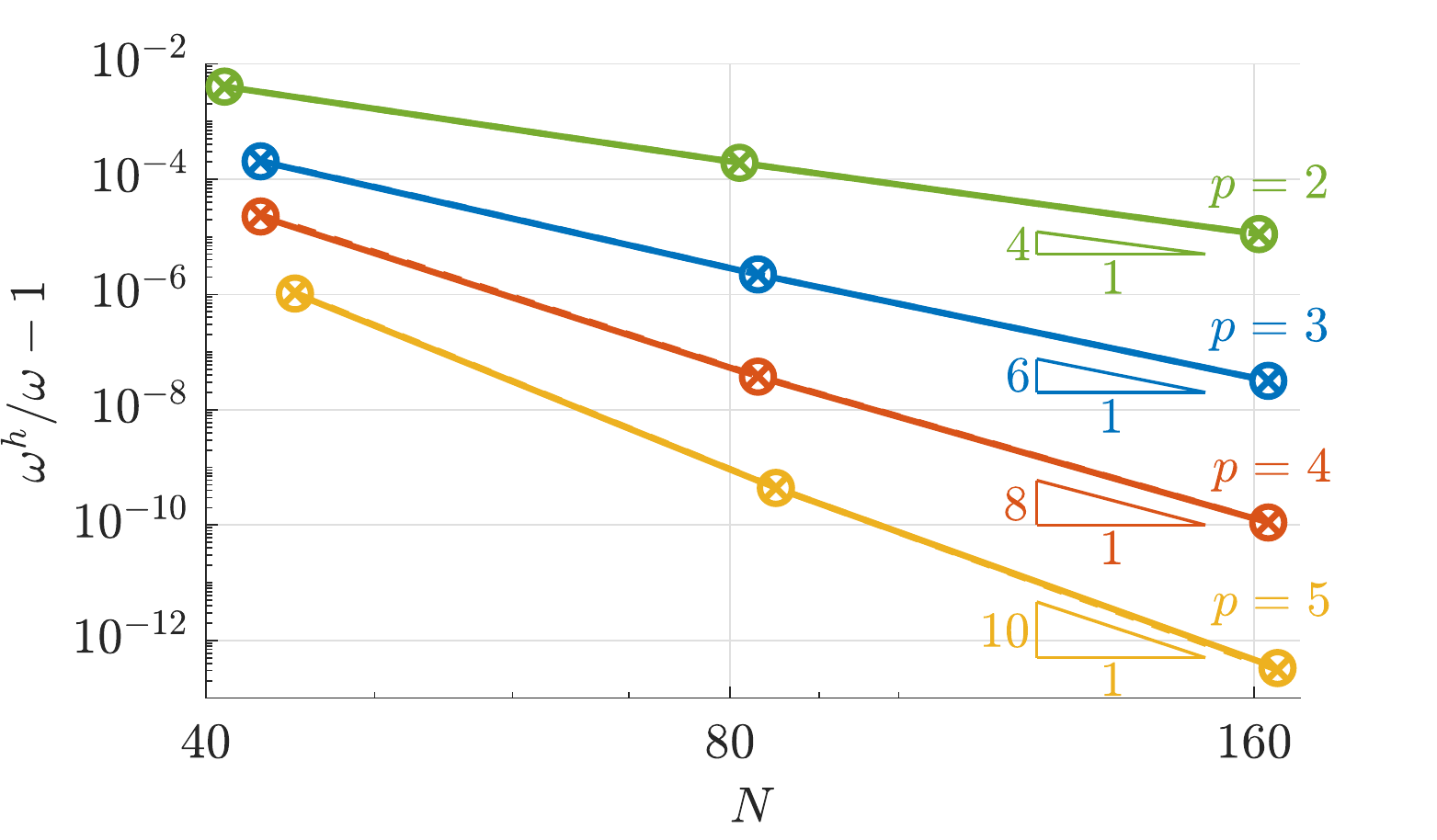} }}
    \subfloat[$L^2$ errors in the mode shapes]{{\includegraphics[width=0.5\textwidth]{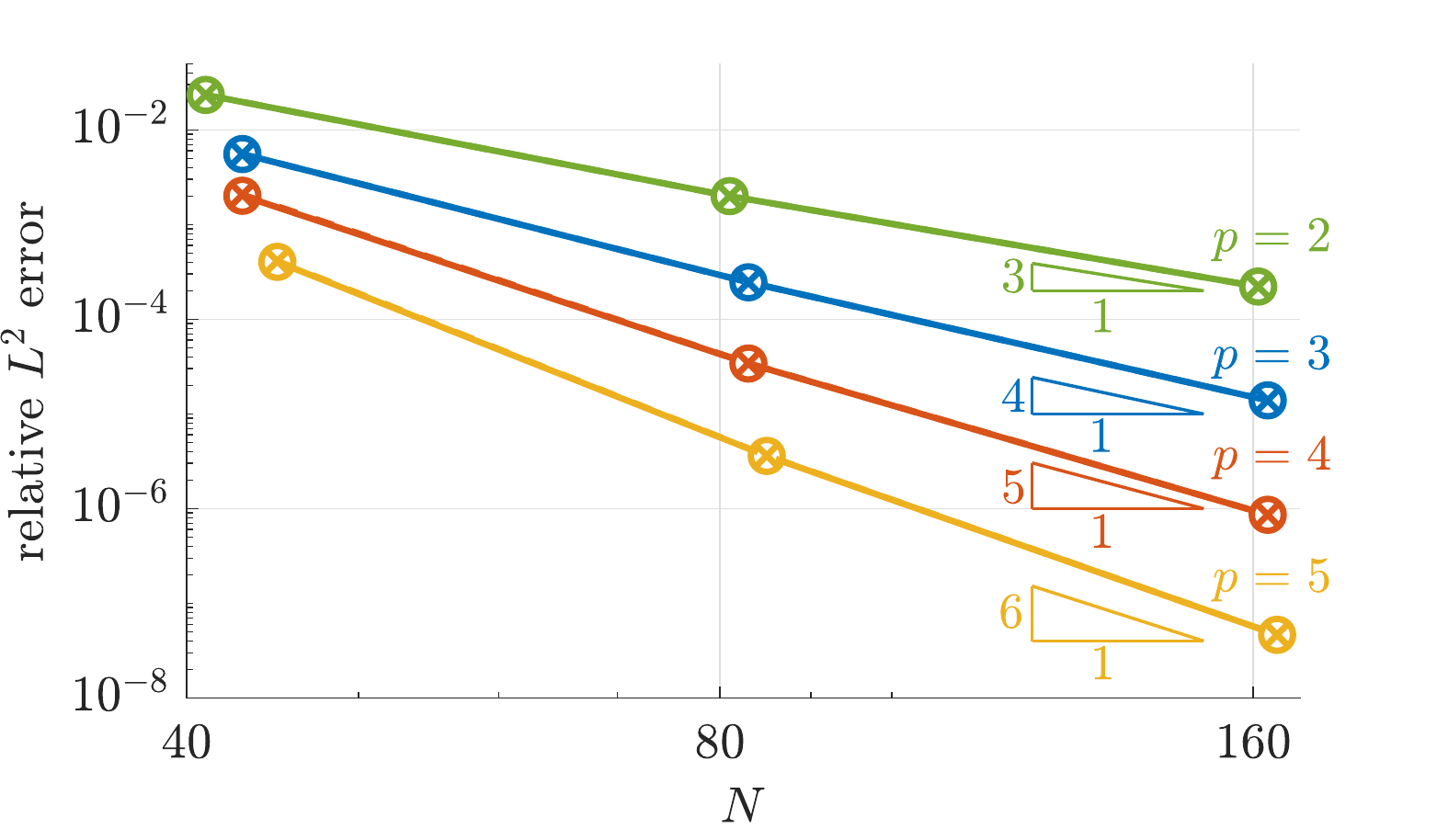} }}
    \vspace{0.2cm}
    \begin{tikzpicture}
    \filldraw[black,line width=1pt, solid] (0.0,0) -- (0.2,0);
    \filldraw[black,line width=1pt] (0.3,0) [fill=none] circle (2pt);
    \filldraw[black,line width=1pt, solid] (0.4,0) -- (0.6,0);
    \filldraw[black,line width=1pt] (0.6,0) node[right]{\footnotesize standard spectrum};
    \filldraw[black,line width=1pt, dashed] (5,0) -- (5.6,0);
    \filldraw[black,line width=1pt] (5.0,0) node[right]{\footnotesize $\boldsymbol{\bigtimes}$};
    \filldraw[black,line width=1pt] (5.6,0) node[right]{\footnotesize improved spectrum};
\end{tikzpicture}
    \caption{Convergence of the relative error in the \textbf{18$^\text{th}$} eigenfrequency and mode of \textbf{a fixed bar}, obtained with \textbf{$2$ patches} of quadratic, cubic, quartic and quintic $C^{p-1}$ B-spline basis functions.}
    \label{fig:convergence_2nd_1d}
\end{figure}

We now consider the free axial vibration of a fixed bar with unit length, unit axial stiffness and unit mass. We employ a multipatch discretization with two patches ($\npa=2$) of $C^{p-1}$ B-splines of different polynomial degrees $p=2$ through $5$ where interior outliers exist (see Table \ref{tab:number_of_outliers_1Dbar}), and $C^0$ patch continuity. 
Figures \ref{fig:normalized_freq_2nd_1d} and \ref{fig:mode_error_2nd_1d} illustrate the normalized frequency and the relative $L^2$ error in the mode shapes of the studied bar, respectively. 
We compare results obtained with multipatch discretizations based on non-perturbed and perturbed eigenvalue problems (plotted in black and blue, respectively).
We include results of the single-patch discretization in Figures \ref{fig:normalized_freq_2nd_1d} and \ref{fig:mode_error_2nd_1d} in gray as the reference solution, and keep the same number of degrees of freedom $N$ for the single- and multipatch discretizations.
 \changedV{The discrete frequencies are ordered as described in Remark \ref{rm_ordering}}.

We first focus on the normalized frequencies plotted in Figure \ref{fig:normalized_freq_2nd_1d}.
The inset figures of Figure \ref{fig:normalized_freq_2nd_1d} focus on the last four frequencies including the outlier frequencies that are present in the upper part of the spectra.
It can be observed that the entire spectrum obtained with multipatch discretizations without perturbations (in black) is accurate except for the $(p-1)$ interior outlier frequencies at the end of the spectrum (see also Table \ref{tab:number_of_outliers_1Dbar}).
These outliers are significantly reduced by our approach (in blue), while the remaining frequencies are not negatively affected.
We observe that the reduction factor increases with increasing $p$ since the outlier frequencies increase.
We note that the reduced normalized outlier frequencies are not at the same level for all cases, as illustrated via small jumps at the end of the blue spectra for $p=4$ and $5$ in Figures \ref{fig:normalized_freq_2nd_1d}c and d.
The mode errors are plotted in Figure \ref{fig:mode_error_2nd_1d}.
We observe that our approach results in smaller errors in the outlier modes (in blue), as discussed in the previous subsection, without affecting the remaining modes.
We conclude that our approach improves the spectral properties of univariate multipatch discretizations without affecting the accuracy of the remaining frequencies and modes.

We then verify that the proposed approach does not negatively affect the accuracy nor the optimal convergence behavior of the lower frequencies and modes. For second-order problems, the optimal convergence rate of the frequency error and the $L^2$ error in the mode is $\mathcal{O}(2p)$ and $\mathcal{O}(p+1)$, respectively \cite{hughes_finite_2003,cottrell_isogeometric_2006}. Figure \ref{fig:convergence_2nd_1d} illustrates the convergence of the relative error in the 18$^\text{th}$ frequency (left) and the $L^2$ errors in the corresponding mode (right) of the bar. We plot these errors vs. the mode number $N$ for degrees $p=2$ through $5$. 
We observe that our approach preserves the optimal accuracy of the frequency and mode. 
We can also see in Figure \ref{fig:convergence_2nd_1d} that the error converges optimally in all cases.

\subsection{Spectral analysis of a fourth-order problem}

\begin{figure}[h!]
    \centering
    \subfloat[$p=3$]{{\includegraphics[width=0.5\textwidth]{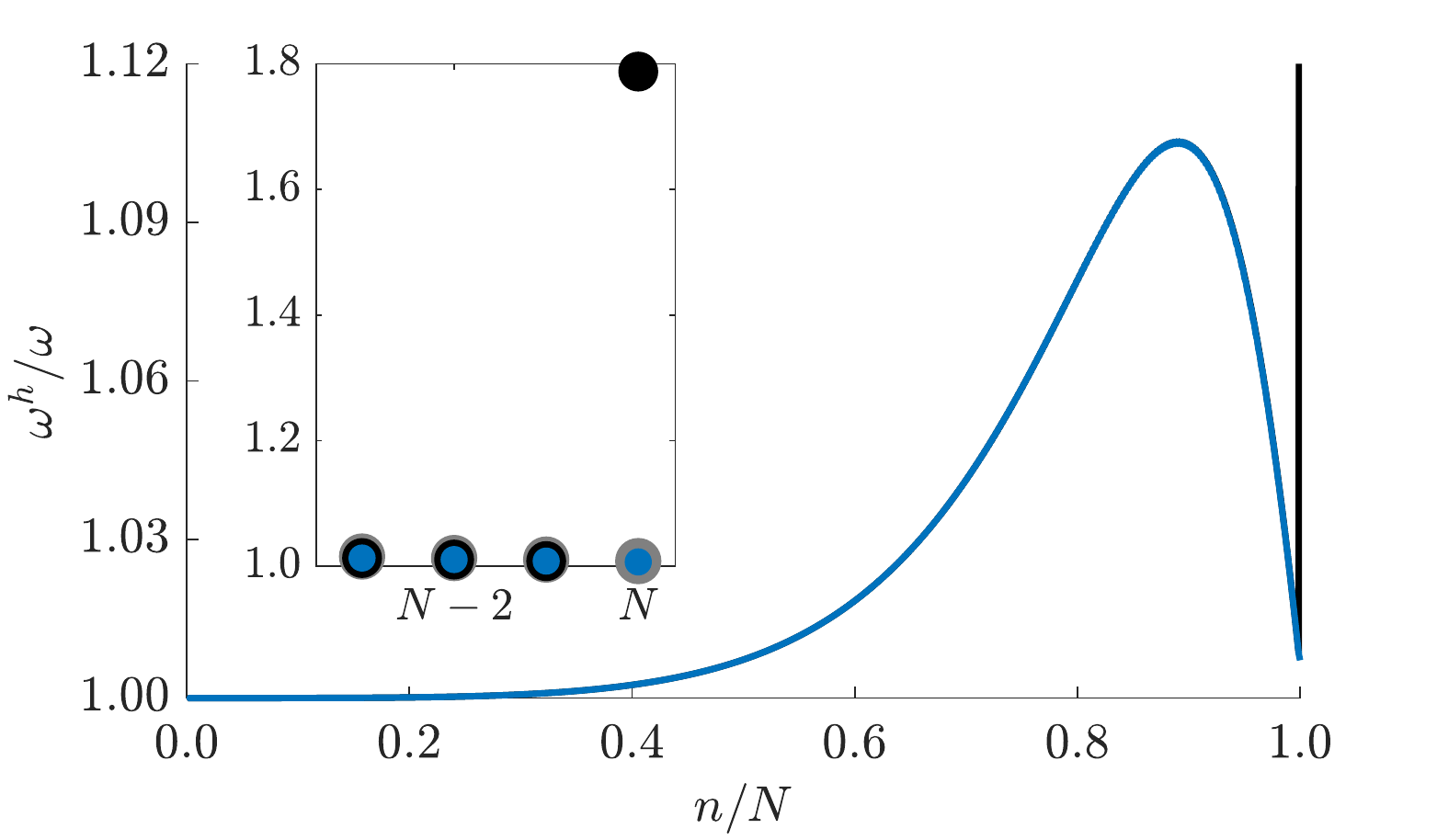} }}
    \subfloat[$p=4$]{{\includegraphics[width=0.5\textwidth]{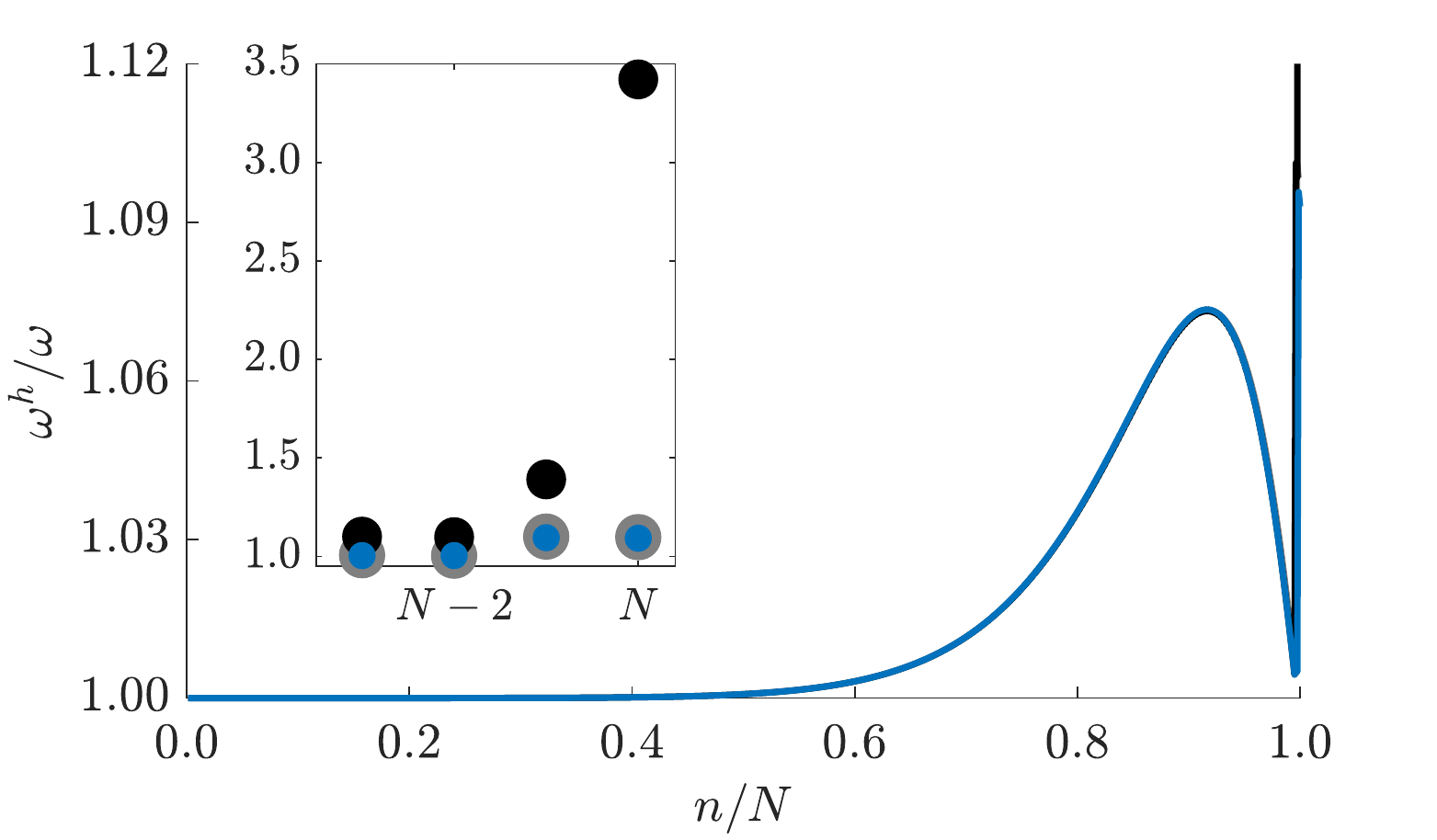} }}

    \subfloat[$p=5$]{{\includegraphics[width=0.5\textwidth]{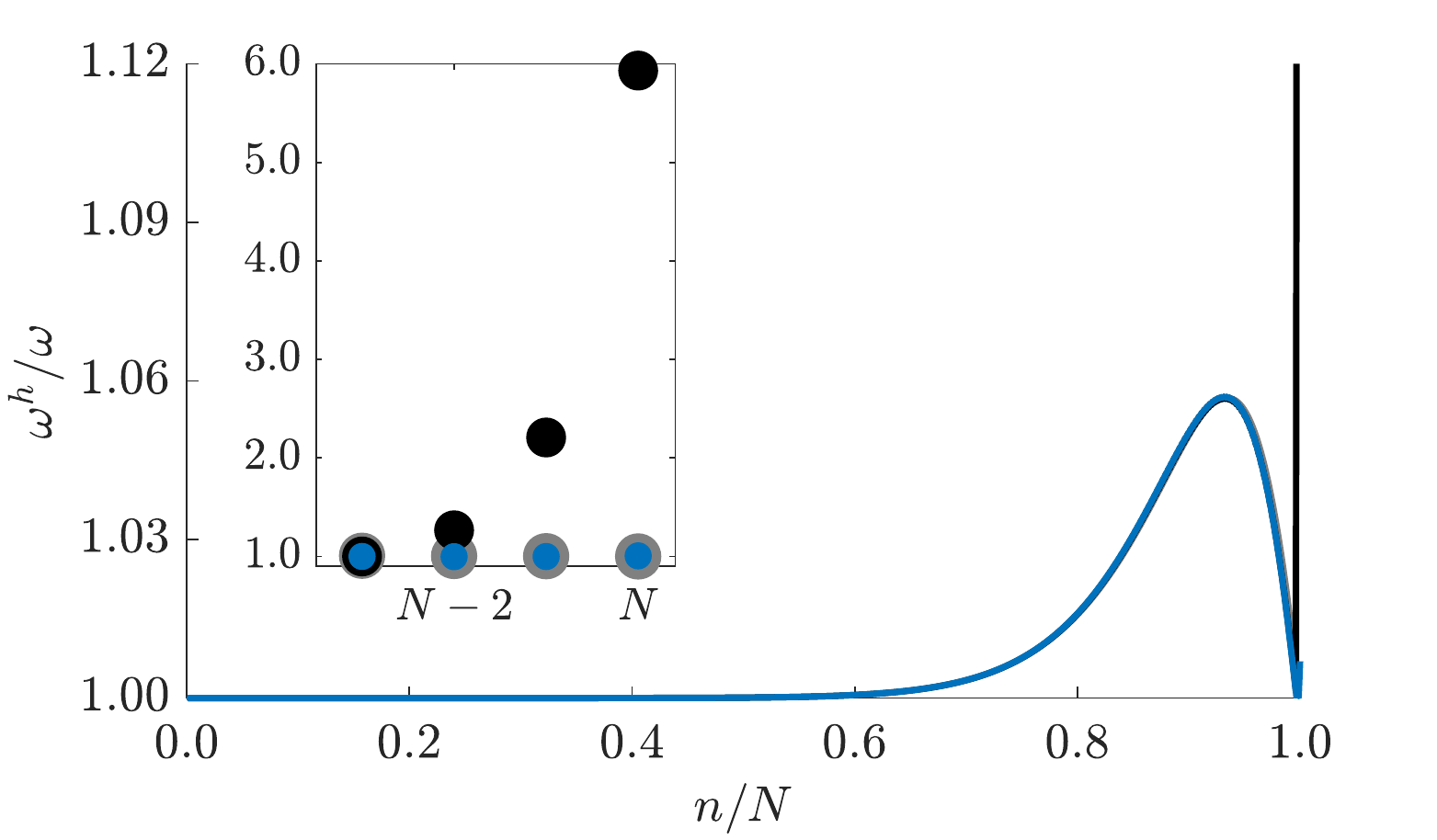} }}
    \subfloat[$p=6$]{{\includegraphics[width=0.5\textwidth]{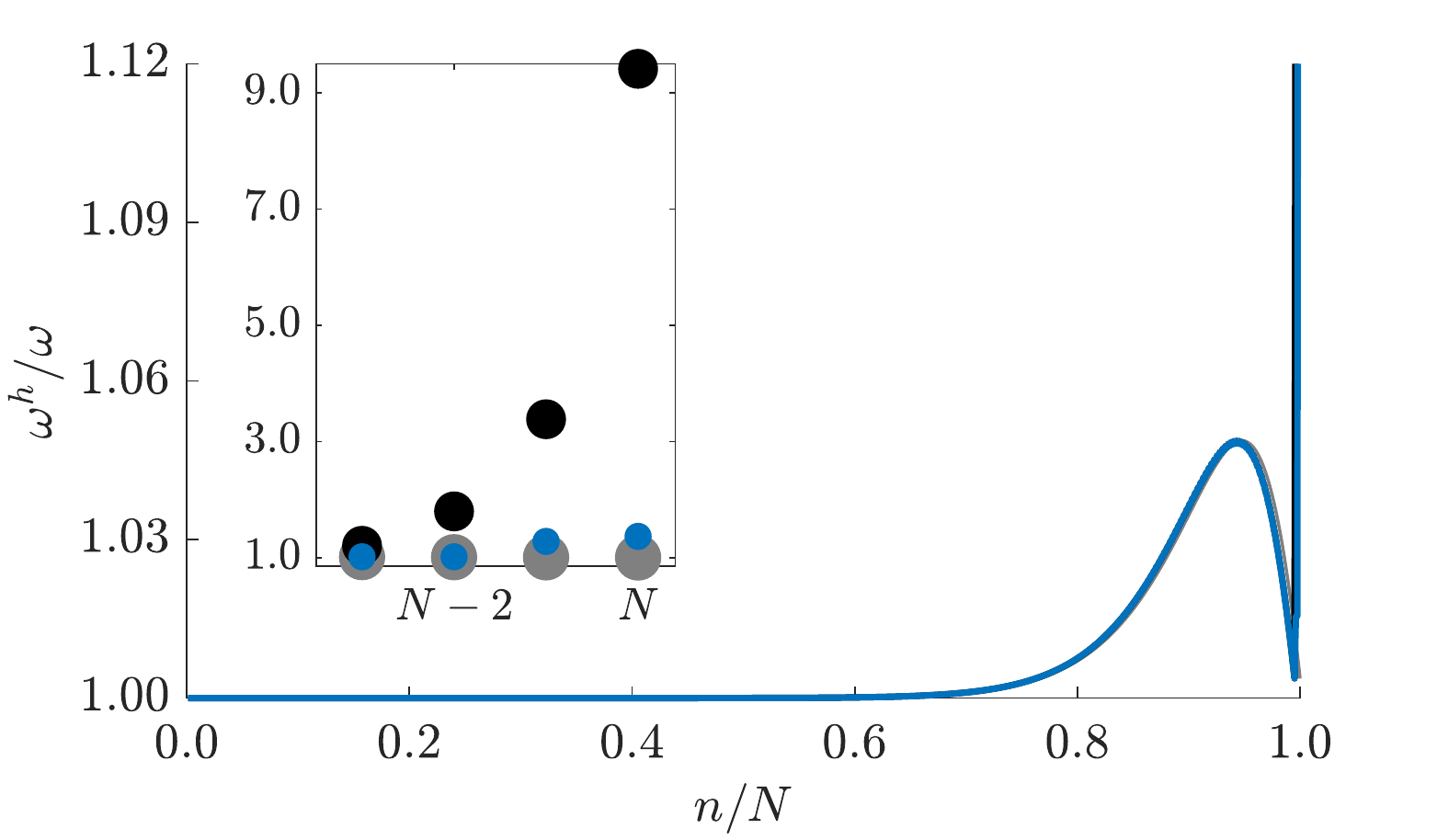} }}
    \vspace{0.2cm}
    \begin{tikzpicture}
    \filldraw[grey1,line width=1pt] (0,0) circle (2pt);
    \filldraw[grey1,line width=1pt] (0,0) node[right]{\footnotesize single-patch};
    \filldraw[black,line width=1pt] (3,0) circle (2pt);
    \filldraw[black,line width=1pt] (3,0) node[right]{\footnotesize multipatch, standard spectrum};
    \filldraw[blue1,line width=1pt] (9,0) circle (2pt);
    \filldraw[blue1,line width=1pt] (9,0) node[right]{\footnotesize multipatch, improved spectrum};
\end{tikzpicture}
    \caption{Normalized frequencies of a freely vibrating \textbf{fixed beam}, computed with \textbf{two patches} of $C^{p-1}$ B-splines and discretized \textbf{with $500$ elements}.}
    \label{fig:normalized_freq_4th_1d}
\end{figure}

\begin{figure}[h!]
    \centering
    \subfloat[$p=3$]{{\includegraphics[width=0.47\textwidth]{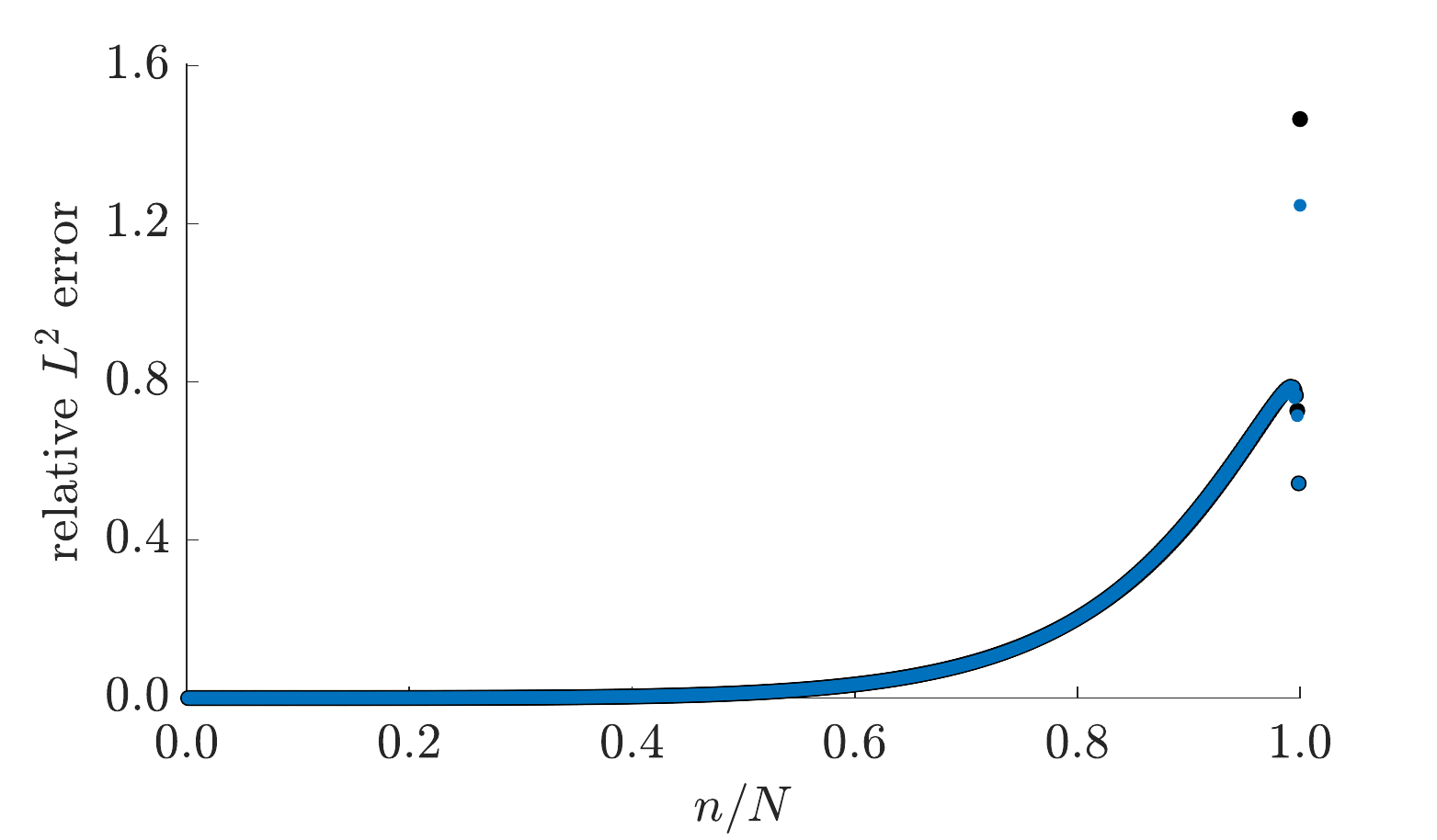} }}
    \subfloat[$p=4$]{{\includegraphics[width=0.47\textwidth]{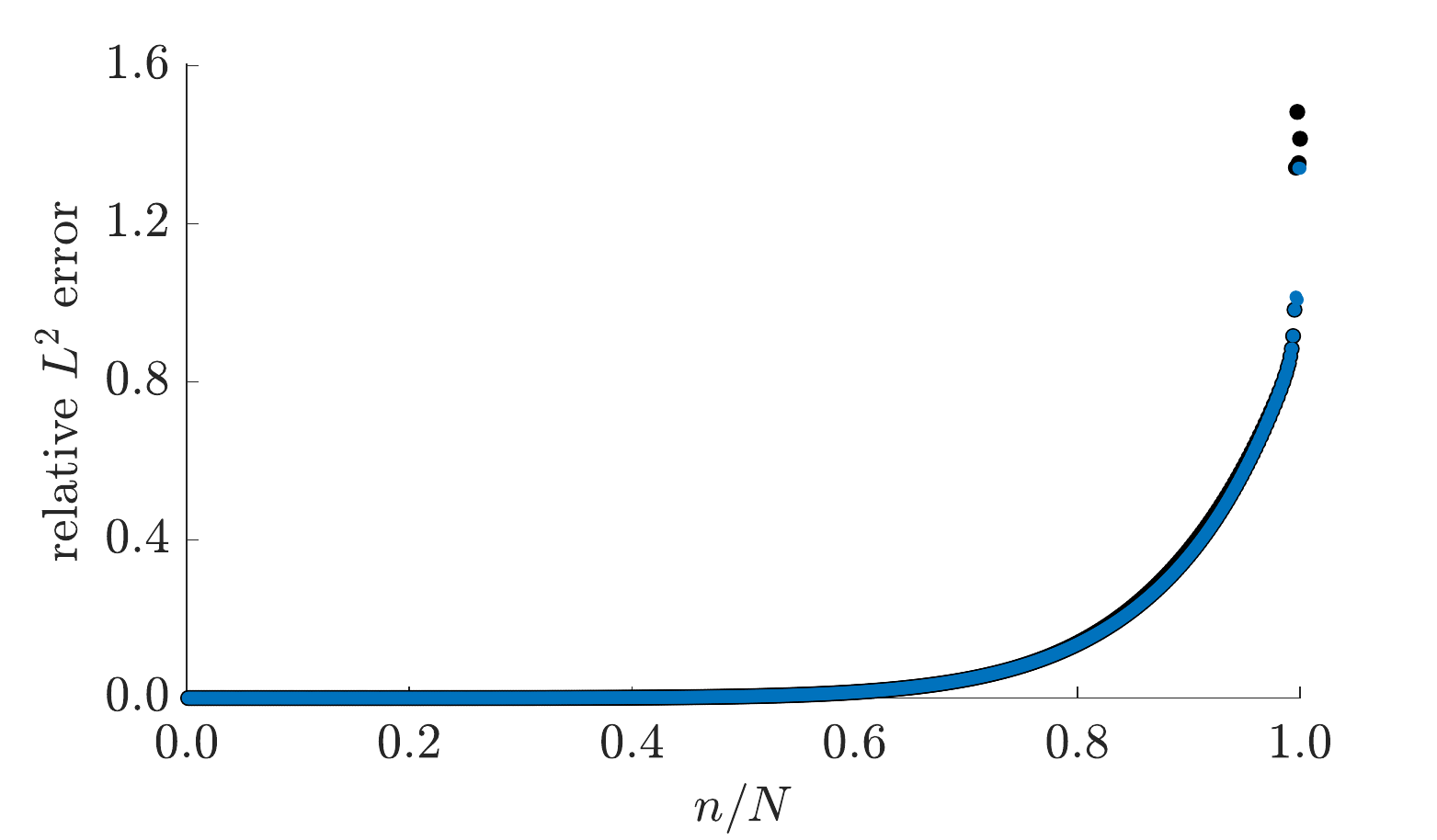} }}

    \subfloat[$p=5$]{{\includegraphics[width=0.47\textwidth]{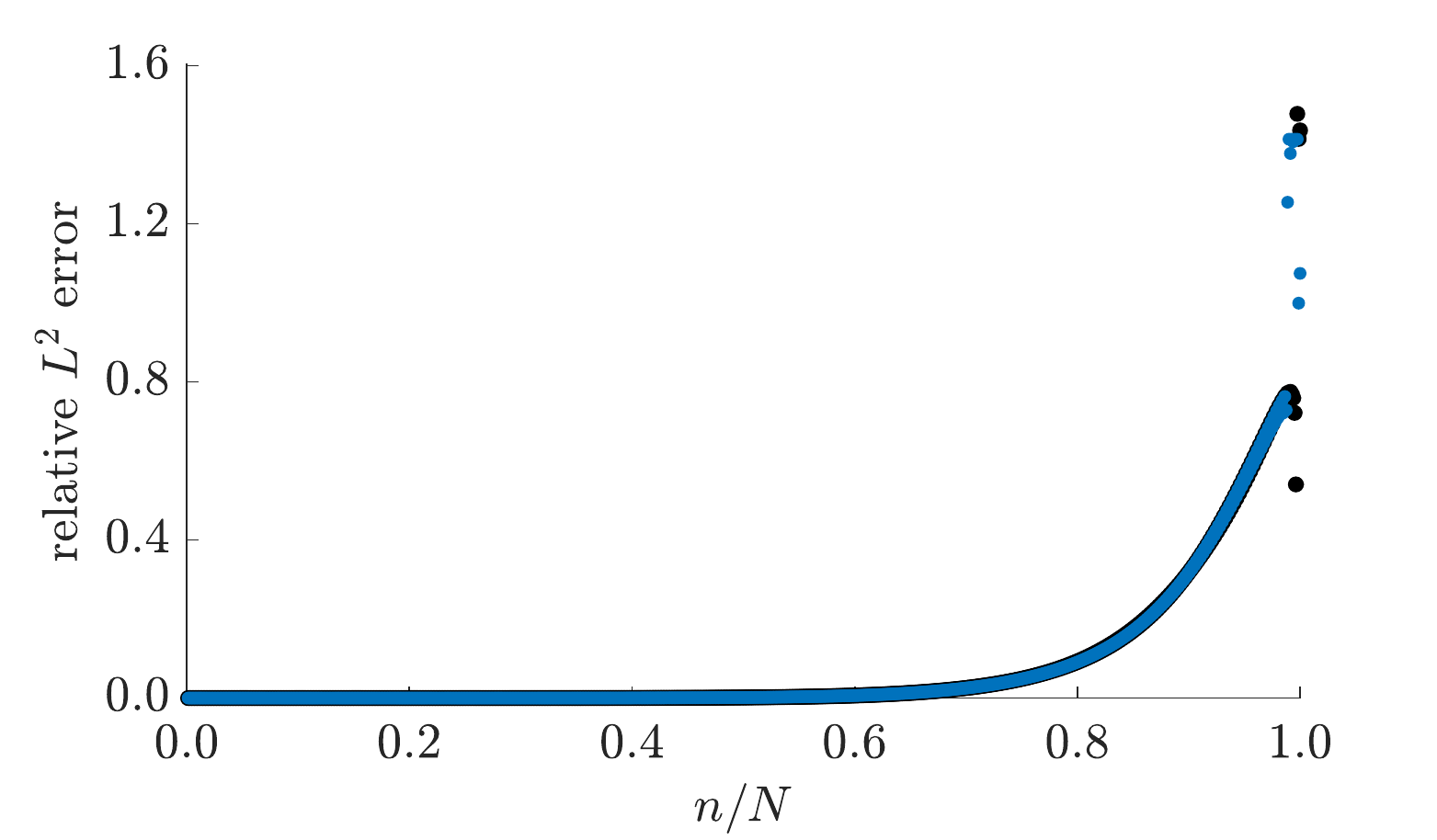} }}
    \subfloat[$p=6$]{{\includegraphics[width=0.47\textwidth]{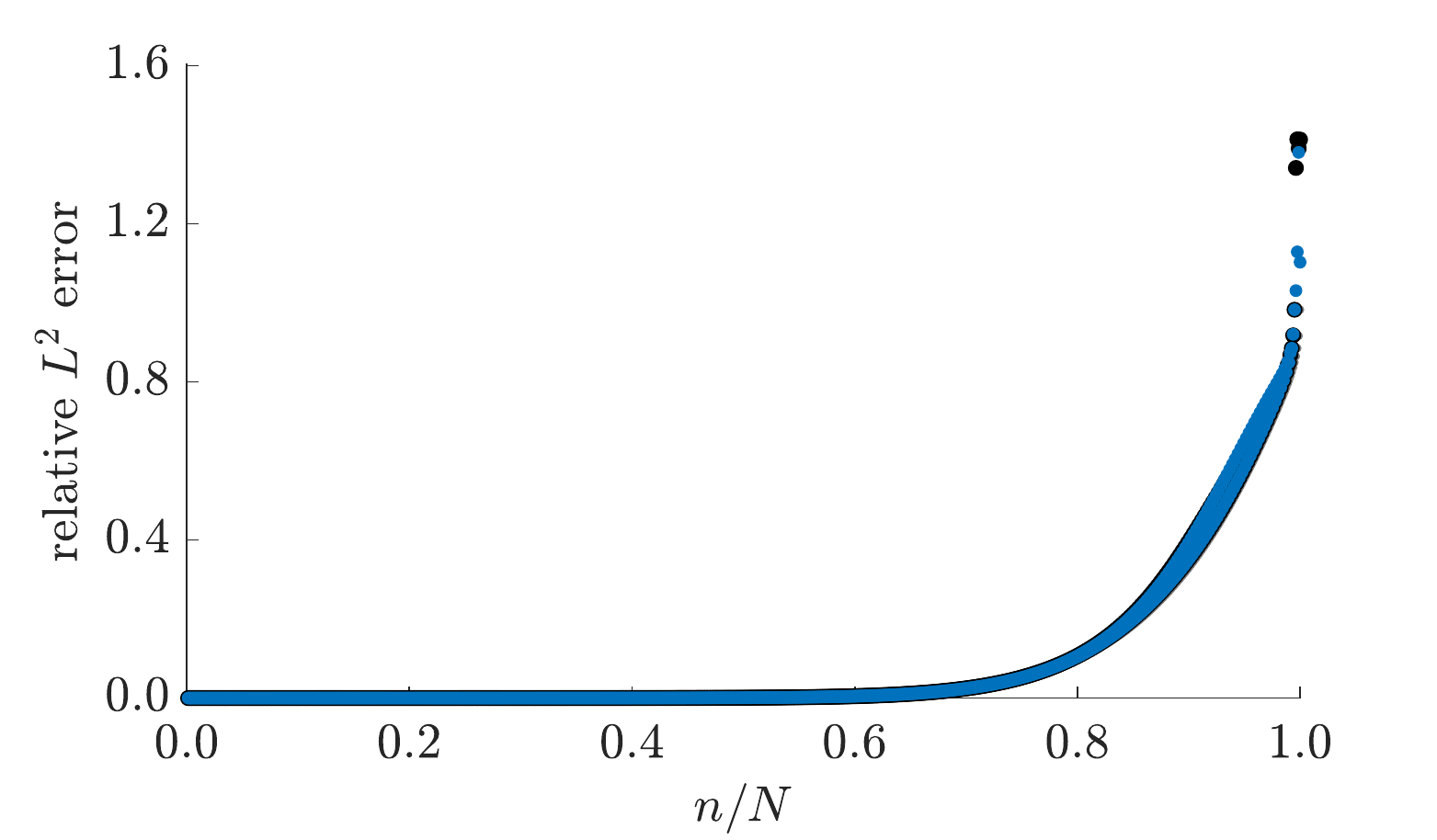} }}
    \vspace{0.2cm}
    \begin{tikzpicture}
    \filldraw[grey1,line width=1pt] (0,0) circle (2pt);
    \filldraw[grey1,line width=1pt] (0,0) node[right]{\footnotesize single-patch};
    \filldraw[black,line width=1pt] (3,0) circle (2pt);
    \filldraw[black,line width=1pt] (3,0) node[right]{\footnotesize multipatch, standard spectrum};
    \filldraw[blue1,line width=1pt] (9,0) circle (2pt);
    \filldraw[blue1,line width=1pt] (9,0) node[right]{\footnotesize multipatch, improved spectrum};
\end{tikzpicture}
    \caption{$L^2$ errors in the mode shapes of a freely vibrating \textbf{fixed beam}, computed with \textbf{two patches} of $C^{p-1}$ B-splines and discretized \textbf{with $500$ elements}.}
    \label{fig:mode_error_4th_1d}

    \vspace{0.4cm}

    \subfloat[Normalized frequency error]{{\includegraphics[width=0.5\textwidth]{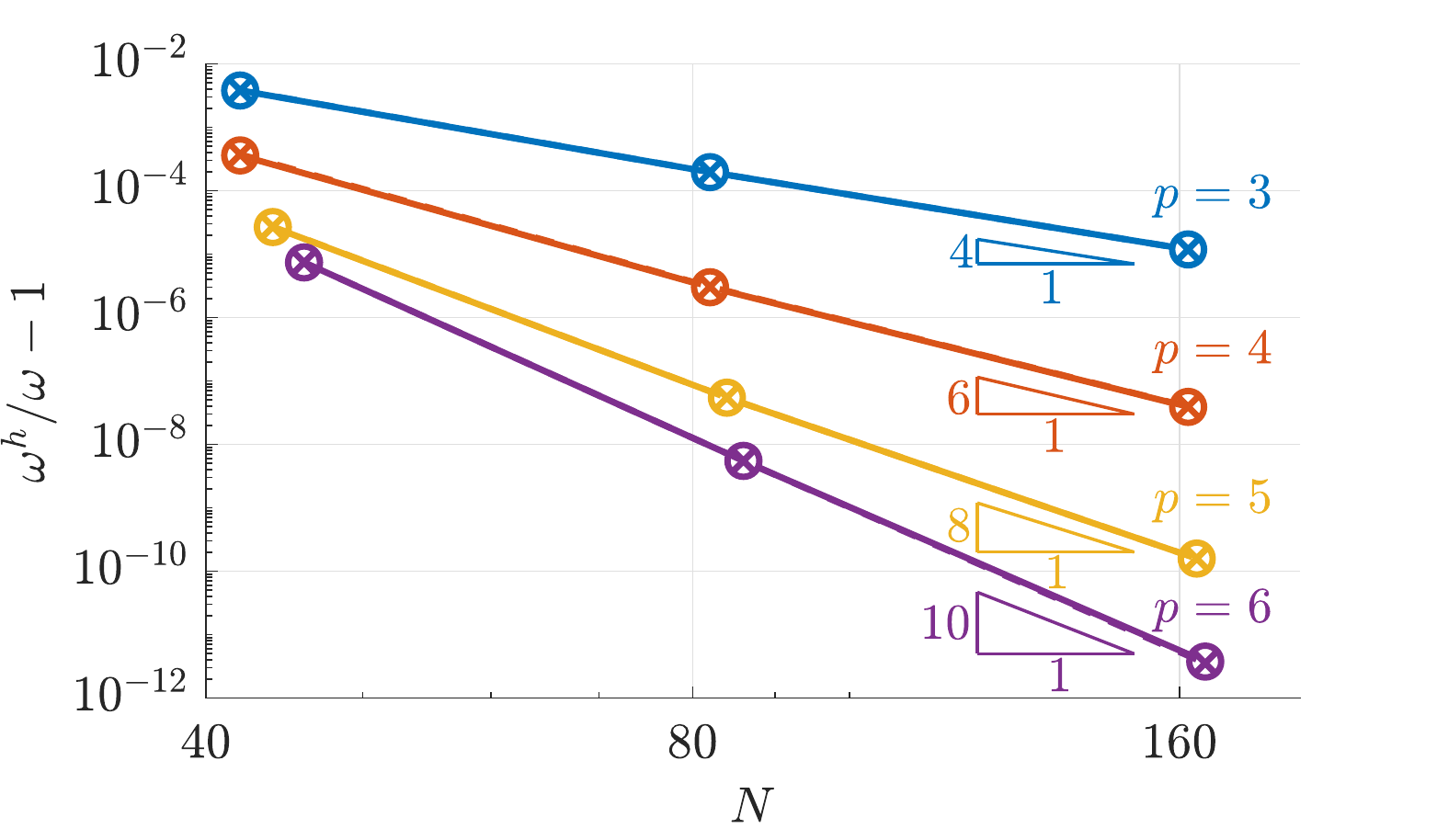} }}
    \subfloat[$L^2$ errors in the mode shapes]{{\includegraphics[width=0.5\textwidth]{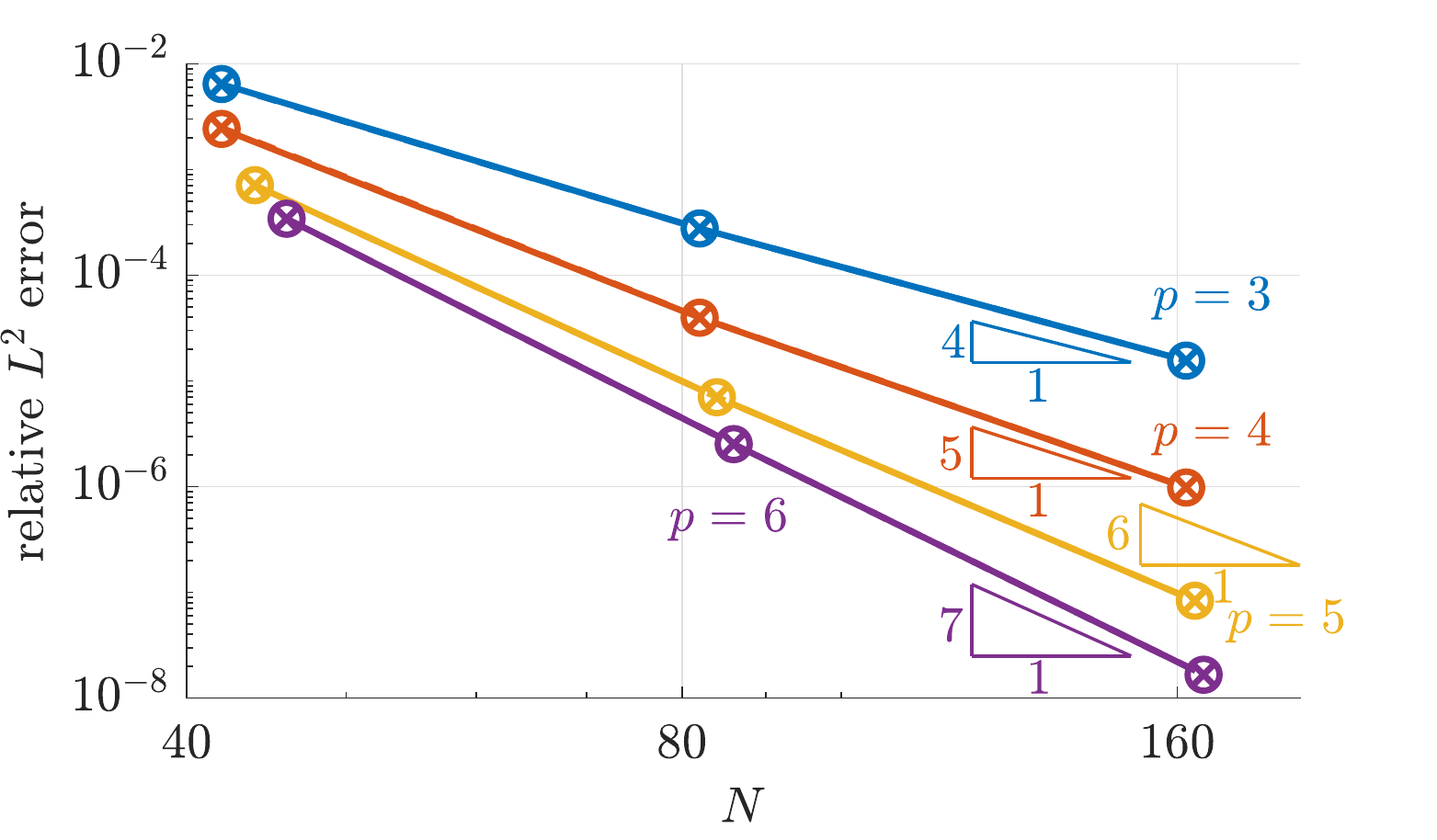} }}
    \vspace{0.2cm}
    \begin{tikzpicture}
    \filldraw[black,line width=1pt, solid] (0.0,0) -- (0.2,0);
    \filldraw[black,line width=1pt] (0.3,0) [fill=none] circle (2pt);
    \filldraw[black,line width=1pt, solid] (0.4,0) -- (0.6,0);
    \filldraw[black,line width=1pt] (0.6,0) node[right]{\footnotesize standard spectrum};
    \filldraw[black,line width=1pt, dashed] (5,0) -- (5.6,0);
    \filldraw[black,line width=1pt] (5.0,0) node[right]{\footnotesize $\boldsymbol{\bigtimes}$};
    \filldraw[black,line width=1pt] (5.6,0) node[right]{\footnotesize improved spectrum};
\end{tikzpicture}
    \caption{Convergence of the relative error in the \textbf{18$^\text{th}$} eigenfrequency and mode of \textbf{a fixed beam}, obtained with \textbf{$2$ patches} of cubic, quartic, quintic and sextic $C^{p-1}$ B-spline basis functions.}
    \label{fig:convergence_4th_1d}
\end{figure}

We then consider the free transverse vibration of a fixed beam with unit length, unit bending stiffness and unit mass. We employ a multipatch discretization with two patches ($\npa=2$) of $C^{p-1}$ B-splines of different polynomial degrees $p=3$ through $6$ where interior outliers exist (see Table \ref{tab:number_of_outliers_1Dbeam}), and with sufficient regularity, i.e.\ $C^1$ patch continuity. 
Figures \ref{fig:normalized_freq_4th_1d} and \ref{fig:mode_error_4th_1d} illustrate the normalized frequency and the $L^2$ error in the mode shapes corresponding to the studied beam, respectively.
We order and \changed{present} the results in the same way as those of the fixed bar in the previous subsection.
We observe a similar effect of the proposed approach on the frequencies and modes as in the case of the fixed bar.
Thus, we 
\changed{conclude} that the proposed approach improves the spectral properties of univariate multipatch discretizations for both second- and fourth-order problems.

We also verify that the proposed approach preserves the optimal accuracy and optimal convergence behavior of the lower frequencies and modes. For fourth-order problems, the optimal convergence rate of the frequency error and the $L^2$ error in the mode is $\mathcal{O}(2(p-1))$ and $\mathcal{O}(p+1)$, respectively \cite{hughes_finite_2003,cottrell_isogeometric_2006}. 
Figure \ref{fig:convergence_4th_1d} demonstrates the convergence of the relative error in the 18$^\text{th}$ frequency (left) and the $L^2$ errors in the corresponding mode (right) of the beam, as functions of the mode number $N$ with polynomial degrees $p=3$ through $6$. 
\changed{We observe that our approach preserves the optimal accuracy and convergence behavior in all cases.}

\section{Generalization to multidimensional discretizations}\label{sec:2d-study}

\changed{The iterative scheme proposed in the previous section is shown to work effectively for univariate multipatch discretizations.}  \changedV{There are two major issues that need to be addressed in order to generalize this scheme to multidimensional discretizations. First}, the scheme requires $(p-1)$ target maximum frequencies, which are typically unknown in practical applications. Second, \changedV{for each continuity constraint \eqref{eq:continuity_constraint} a corresponding outlier mode needs to be identified that maximizes \eqref{eq:maximum}. This would necessarily imply that all modes need to be precomputed.} In this section, we tackle both these problems by further simplifying the current approach \eqref{eq:vdgep_perturbed}, giving rise to a new iterative scheme for \changed{estimating the scaling parameters of the perturbation terms}. We then demonstrate via spectral analysis of two-dimensional second- and fourth-order problems that the simplified approach is able to effectively improve the spectra of multipatch discretizations.

\subsection{A pragmatic approach \changedV{to} parameter estimation}\label{sec:pragmatic_approach}
\changedV{We assume a uniform mesh with mesh size $h$. The variational formulation \eqref{eq:vdgep_perturbed} can then be simplified by choosing the $l^{\text{th}}$ parameter as: $\alpha^l = \alpha \; h^{2l-2}$ and $\beta^l = \beta \; h^{2l-2}$}. The scaling factor $h^{2l-2}$ is based on the dimensional consistency of the inner-products of the $l^\text{th}$ derivatives in the perturbation \eqref{eq:perturbation_bilinear_form}. The proposed variational formulation of the perturbed eigenvalue problem \eqref{eq:vdgep_perturbed} becomes: find $(\tilde{\eigenvec}^h_\noMode, \tilde{\omega}^h_\noMode) \in \mathcal{V}^h \times \mathbb{R}^+$, for $n = 1,2,\ldots,N$, such that $\forall \, v^h \in \mathcal{V}^h \subset \mathcal{V}$:
\begin{align}
    a(\tilde{\eigenvec}_\noMode^h, v^h) + \alpha \, \sum_{l=1}^{p-1} \, h^{2l-2} \, c^l(\tilde{\eigenvec}_\noMode^h, v^h) = \left(\tilde{\omega}_\noMode^h\right)^2 \, \left[b(\tilde{\eigenvec}_\noMode^h, v^h) + \beta \, \sum_{l=1}^{p-1} \, h^{2l-2} \, c^l(\tilde{\eigenvec}_\noMode^h, v^h) \right] \, . \label{eq:vdgep_perturbed2}
\end{align}
The corresponding matrix equation \eqref{eq:dgep_perturbed} is then simplified to:
\begin{align}
    \left( \, \mat{K} + \alpha \, \underbrace{\sum_{l=1}^{p-1} \, h^{2l-2} \, \mat{K}_{\Gamma}^l}_{\mat{K}_\Gamma}  \right) \, \mat{\tilde{\eigenvec}}_\noMode^h \; = \; \left(\tilde{\omega}_\noMode^h \right)^2 \, \left( \, \mat{M} + \beta \, \underbrace{\sum_{l=1}^{p-1} \, h^{2l-2} \, \mat{K}_{\Gamma}^l}_{\mat{K}_\Gamma} \right) \, \mat{\tilde{\eigenvec}}_\noMode^h \, , \label{eq:dgep_perturbed2}
\end{align}

\begin{figure}[h!]
    \centering
    \begin{center}
        \captionsetup{labelfont={bf,color=black}}
        \captionof{algorithm}{\changed{Iterative estimation of the parameters involved in  \eqref{eq:vdgep_perturbed2}. The algorithm terminates when the maximum frequency does not correspond to an outlier anymore (see Remark \ref{rm_alg1}).}}
        \label{fig:parameter_estimation}
        \includegraphics[width=0.62\textwidth]{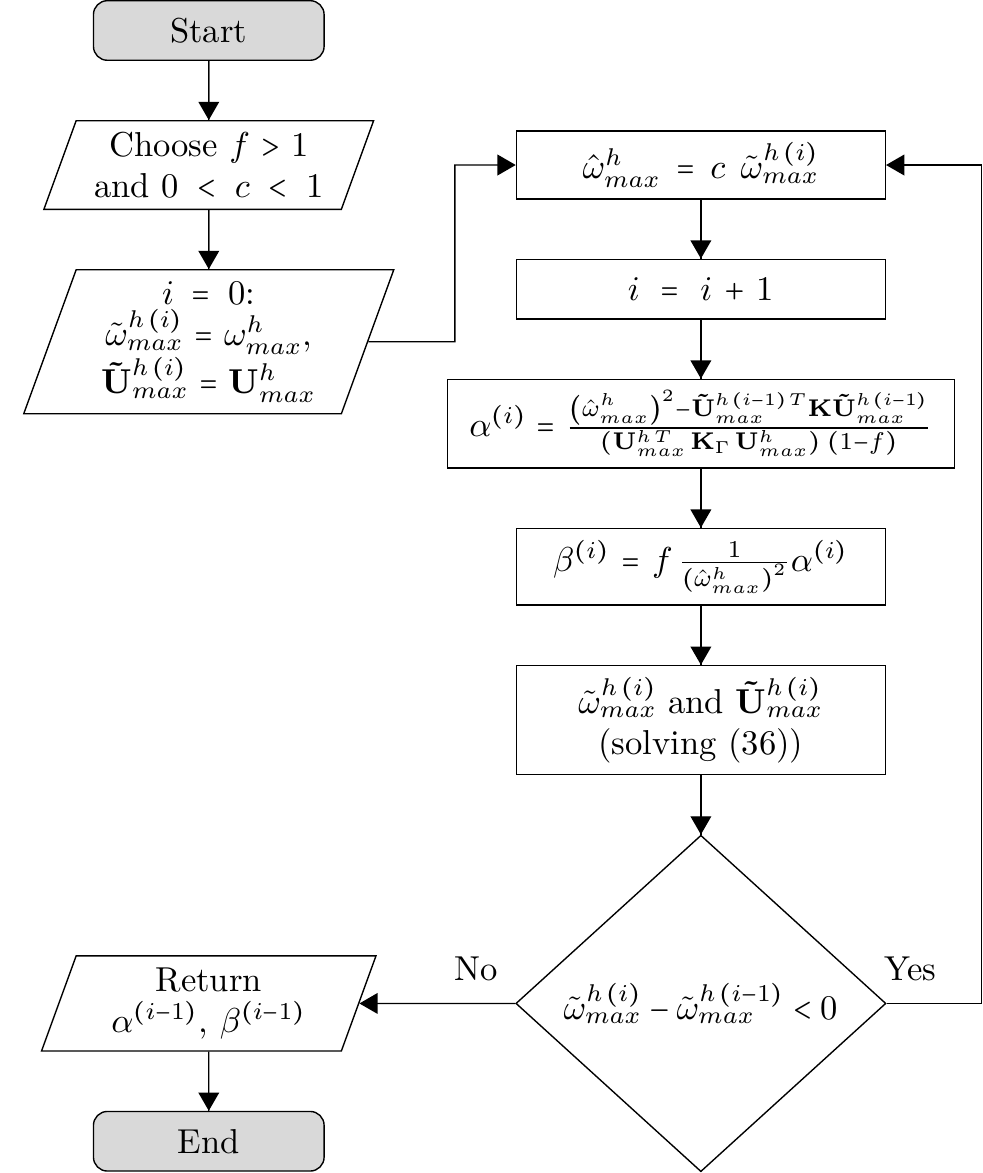}      
    \end{center}

    \vspace{0.4cm}
    \subfloat[Discrete frequency spectrum]{{\includegraphics[width=0.47\textwidth]{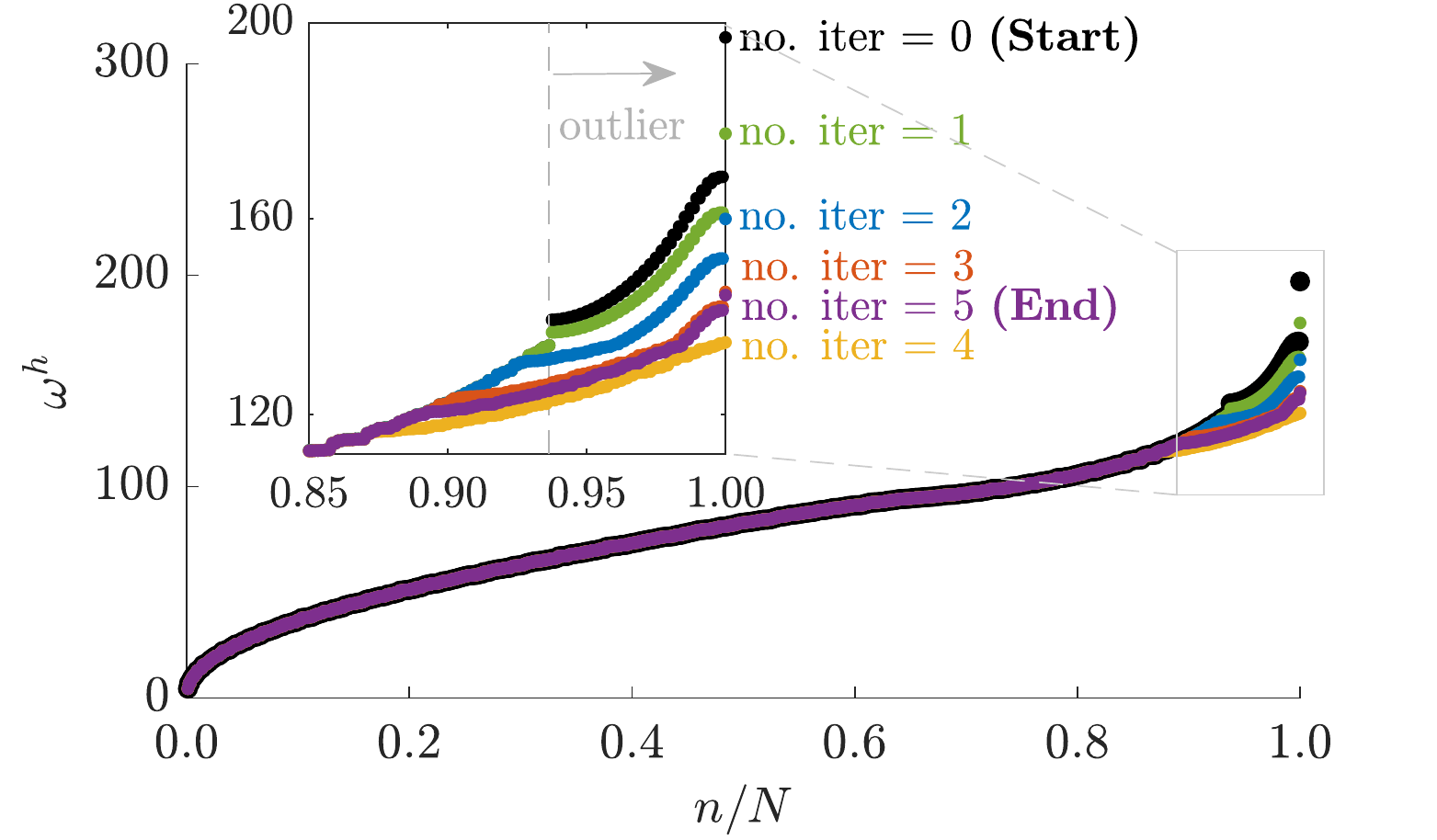}}}
    \subfloat[Evolution of $\alpha$ and $\beta$]{{\includegraphics[width=0.47\textwidth]{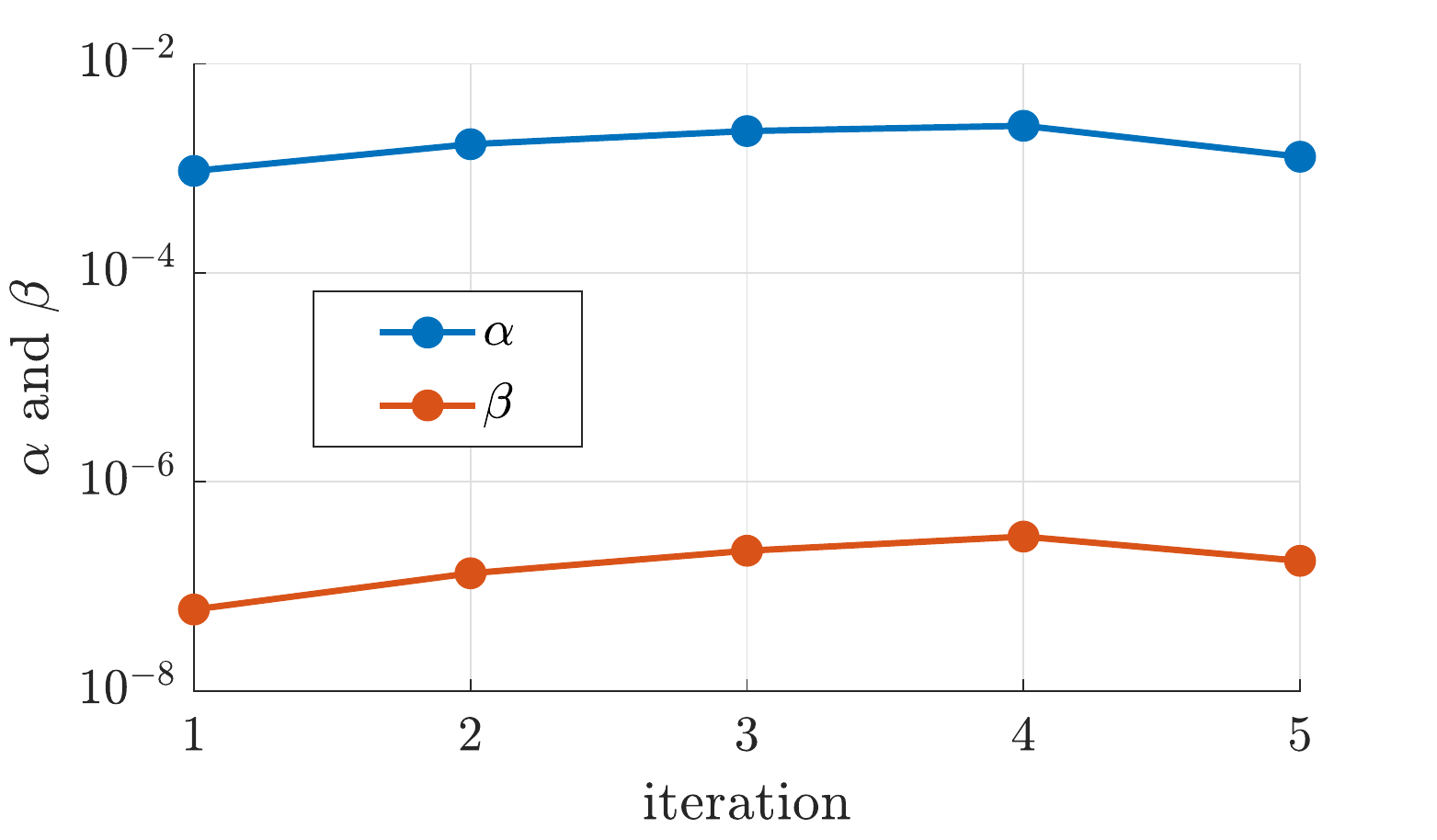}}}
    \caption{Discrete frequencies and parameters, $\alpha$, $\beta$ in each iteration, corresponding to a freely vibrating \textbf{square membrane with fixed boundary conditions}, computed with $2 \times 2$ patches and $15 \times 15$ elements of \textbf{quadratic B-splines} ($p=2$) via Algorithm \ref{fig:parameter_estimation} with $\boldsymbol{f=2,\,c=0.9}$.}
    \label{fig:iteration}
\end{figure}

On this basis, we design \changed{Algorithm \ref{fig:parameter_estimation} for estimating the parameters $\alpha$ and $\beta$ in \eqref{eq:vdgep_perturbed2} (or \eqref{eq:dgep_perturbed2}). 
Here, we iteratively estimate $\alpha$ and $\beta$ such that the maximum frequency is reduced by a target factor $c \in (0,1)$ in each iteration, and thus do not require a choice of the target maximum frequency or any analytical value. In Algorithm \ref{fig:parameter_estimation}, the target maximum frequency in each iteration, $\hat{\omega}^h_{max}$, is a fraction of the maximum frequency obtained from the previous iteration}. 
To obtain $\alpha$ and $\beta$ in each iteration, instead of the outlier, we focus on the maximum frequency and corresponding mode, such that no identification of outlier modes is required. 
The equation for estimating these parameters in \changed{Algorithm \ref{fig:parameter_estimation}} is based on the iterative scheme \eqref{eq:iterative_alpha}-\eqref{eq:iterative_beta} described in the previous section. 
We note that to ensure stability, we employ the unperturbed outlier mode $\eigenvec^h_{max}$ instead of the perturbed one $\tilde{\eigenvec}^{h}_{max}$ in the denominator for computing $\alpha$, since $\tilde{\eigenvec}^{h}_{max}$ is not necessarily an outlier mode in all iterations and thus could result in infinitesimal values in the denominator. 

Figure \ref{fig:iteration} illustrates how the Algorithm \ref{fig:parameter_estimation} affects the discrete frequency (left) and changes the parameters $\alpha$ and $\beta$ (right) in each iteration for a fixed square membrane (see next subsection for details of the example). 
We observe that the maximum outlier frequency is reduced in each iteration, from the first through the fourth iteration (see Figure \ref{fig:iteration}a), and the parameters $\alpha$ and $\beta$ increase in the same iteration (see Figure \ref{fig:iteration}b). 
We can also see in the inset figure of Figure \ref{fig:iteration}a that our approach reduces not only the outlier frequencies \changed{(data points on the right of the dashed gray line)} but also the non-outlier ones \changed{(data points on the left of the dashed gray line)}. 
In the last (fifth) iteration, the maximum frequency increases (purple curve in Figure \ref{fig:iteration}a), which then meets the stopping criterion defined in \changed{Algorithm \ref{fig:parameter_estimation}}. This increase in the maximum frequency corresponds to decreasing parameters $\alpha$ and $\beta$ (see Figure \ref{fig:iteration}b).  
Then, the output of \changed{Algorithm \ref{fig:parameter_estimation}} are the parameters $\alpha$ and $\beta$ of the previous (fourth) iteration.

\begin{remark}\label{rm_alg1} 
Algorithm \ref{fig:parameter_estimation} stops when the perturbation reduces the outlier frequencies by such a factor that they become smaller than the highest non-outlier frequency, see Figure \ref{fig:iteration}a.  
\end{remark}

\begin{remark}\label{rm_c_factor1} 
The factor $c \in (0,1)$ determines the reduction in each iteration of Algorithm \ref{fig:parameter_estimation}. Due to the chosen stopping criterion (Remark \ref{rm_alg1}), the resulting maximum frequency lies within a range of $100 \cdot(1-c) \%$ of the minimum value that can be achieved. For example, a choice of $c=0.9$ results in a reduced maximum frequency within $10$\% of the minimum value.
\end{remark}

\begin{figure}[t!]
    \centering
    \subfloat[$p=2$]{{\includegraphics[width=0.5\textwidth]{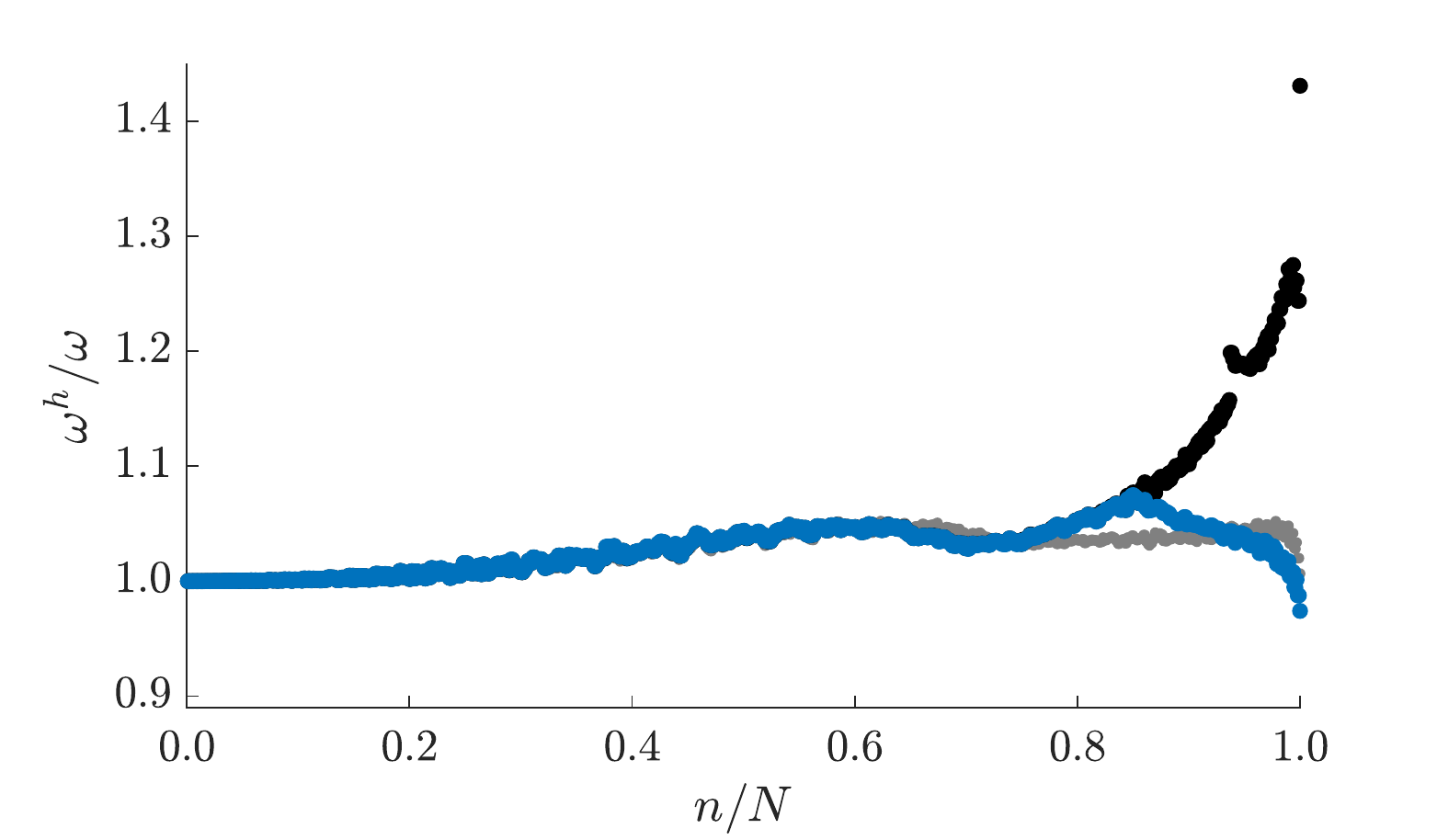} }}
    \subfloat[$p=3$]{{\includegraphics[width=0.5\textwidth]{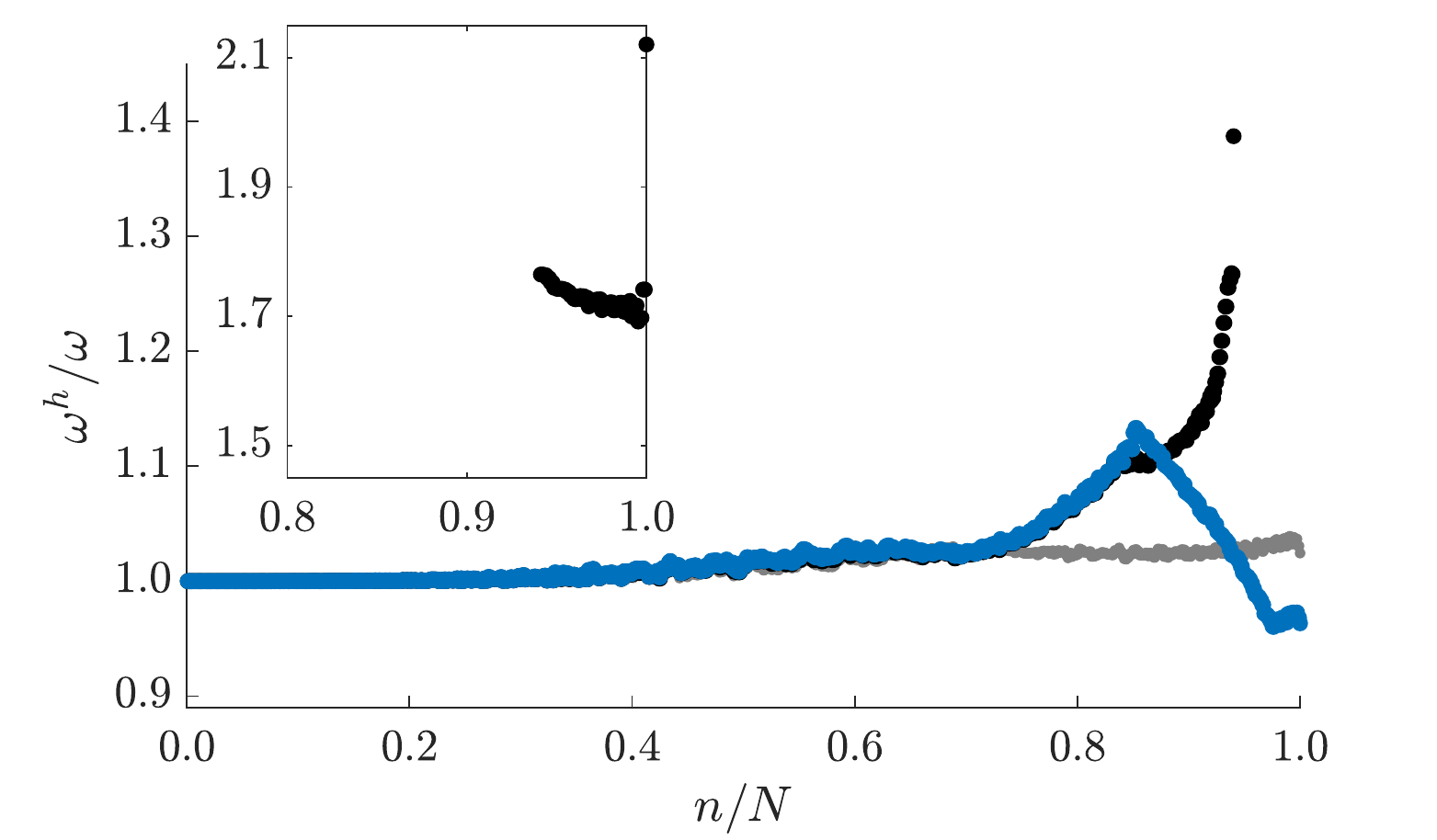} }}

    \subfloat[$p=4$]{{\includegraphics[width=0.5\textwidth]{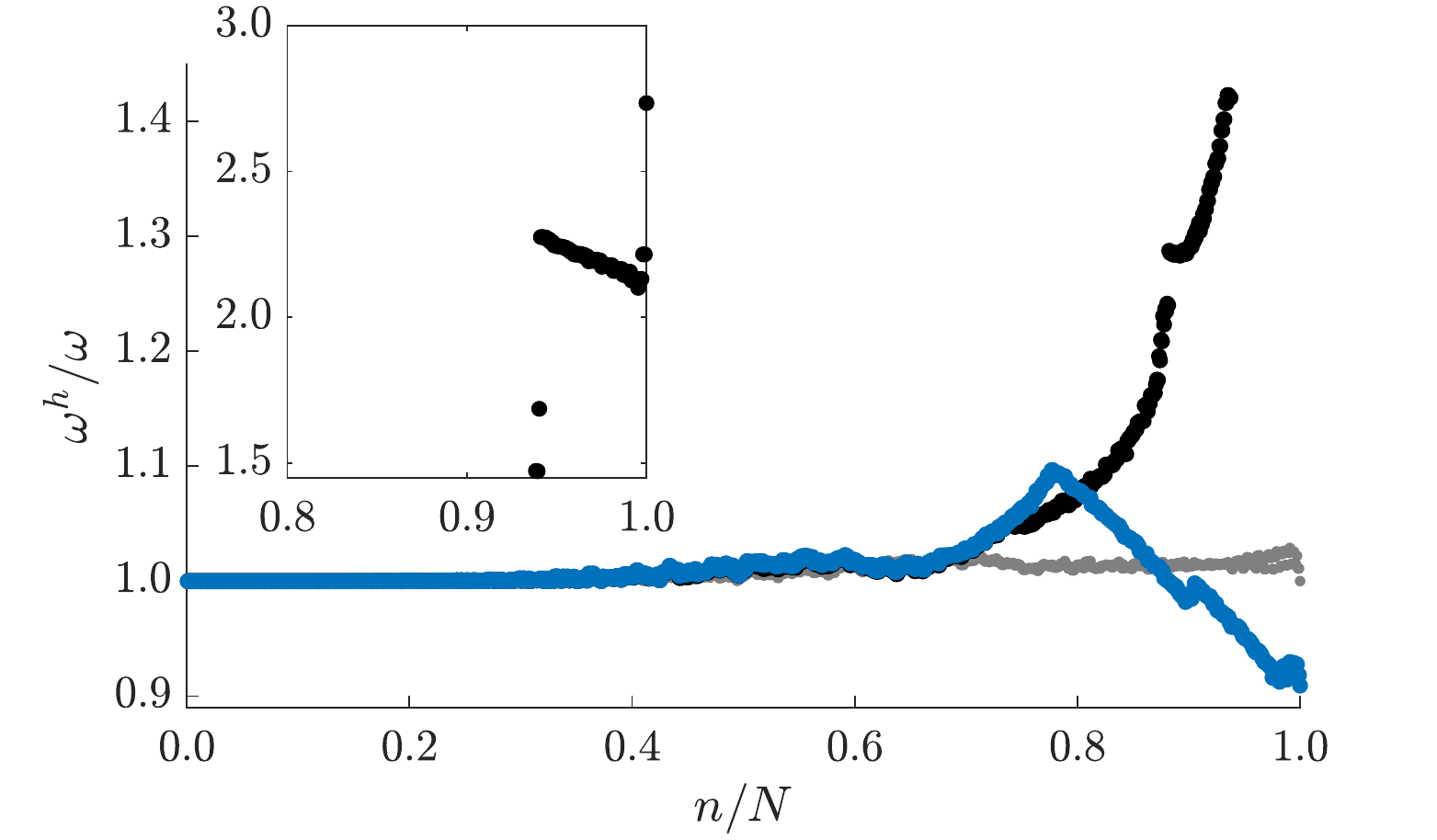} }}
    \subfloat[$p=5$]{{\includegraphics[width=0.5\textwidth]{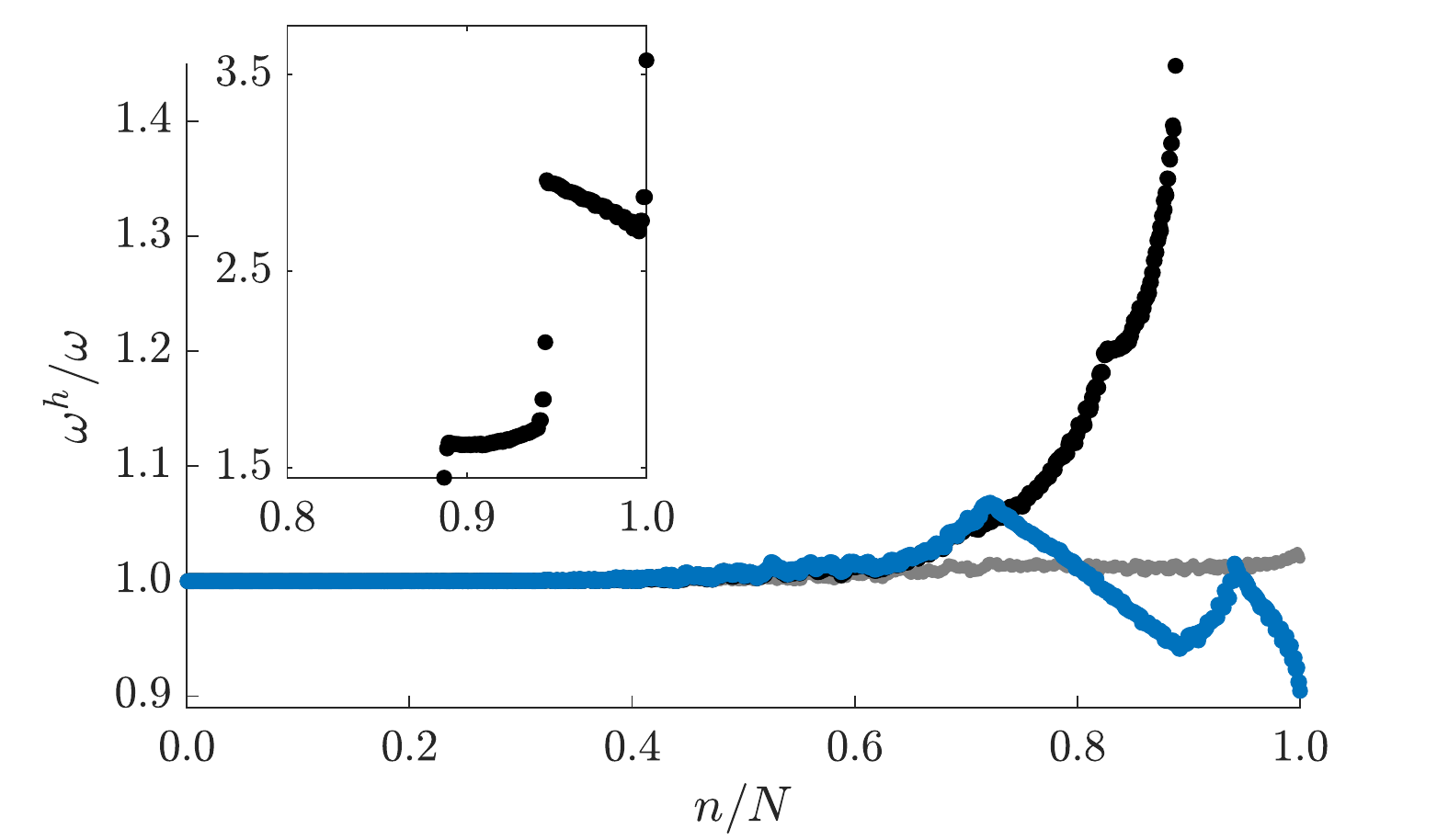} }}
    \vspace{0.2cm}
    \begin{tikzpicture}
    \filldraw[grey1,line width=1pt] (0,0) circle (2pt);
    \filldraw[grey1,line width=1pt] (0,0) node[right]{\footnotesize single-patch};
    \filldraw[black,line width=1pt] (3,0) circle (2pt);
    \filldraw[black,line width=1pt] (3,0) node[right]{\footnotesize multipatch, standard spectrum};
    \filldraw[blue1,line width=1pt] (9,0) circle (2pt);
    \filldraw[blue1,line width=1pt] (9,0) node[right]{\footnotesize multipatch, improved spectrum};
\end{tikzpicture}
    \caption{Normalized frequencies of a freely vibrating \textbf{square membrane with fixed boundary conditions}, computed with \textbf{$2 \times 2$ patches of $C^{p-1}$ B-splines}. Each patch is \textbf{discretized with $15 \times 15$ elements}.}
    \label{fig:normalized_freq_2nd_2d}
\end{figure}

\changed{Algorithm \ref{fig:parameter_estimation}} updates only the maximum frequency and the corresponding mode in each iteration. It requires two input parameters:  
a scaling factor $f > 1$ between $\alpha$ and $\beta$, which can be chosen as discussed in the previous section, and a target reduction factor $c \in (0,1)$ in  each iteration. 
A choice of a small $c$ results in a large reduction step of the maximum frequency and a small number of iterations. The reduced maximum frequency, however, could be far away from the minimum value that can be achieved (see Remark \ref{rm_c_factor1}), i.e.\  the maximum frequency is reduced ineffectively. Moreover, a very small $c$ leads to a very small value in the denominator of the equation for $\beta$ (see Algorithm \ref{fig:parameter_estimation}), negatively affecting the convergence of $\beta$.
Large values of $c$ avoid this issue, but require more iterations. For our numerical studies in the remainder of this paper, we choose $f=2$ and $c=0.9$ for all cases. We find that this choice of factor $c$ typically requires only up to five iterations to sufficiently reduce the maximum frequency, i.e. the maximum frequency is significantly reduced within $10$\% of the lowest possible value (see also \ref{rm_c_factor1}) with a small number of iterations.

\subsection{Spectral analysis of 2D second- and fourth-order model problems}

\begin{figure}[t!]
    \centering
    \subfloat[$p=2$]{{\includegraphics[width=0.5\textwidth]{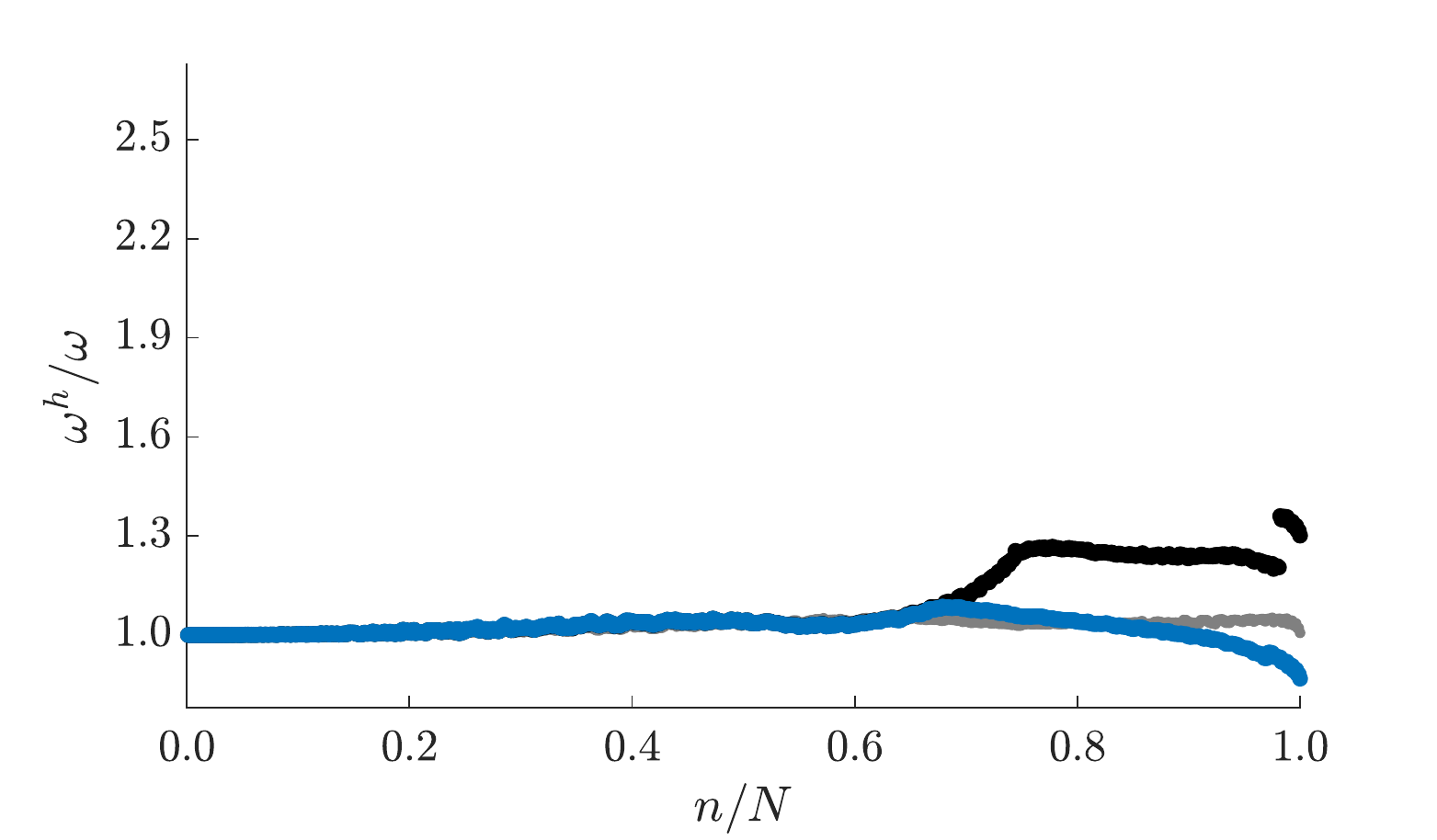} }}
    \subfloat[$p=3$]{{\includegraphics[width=0.5\textwidth]{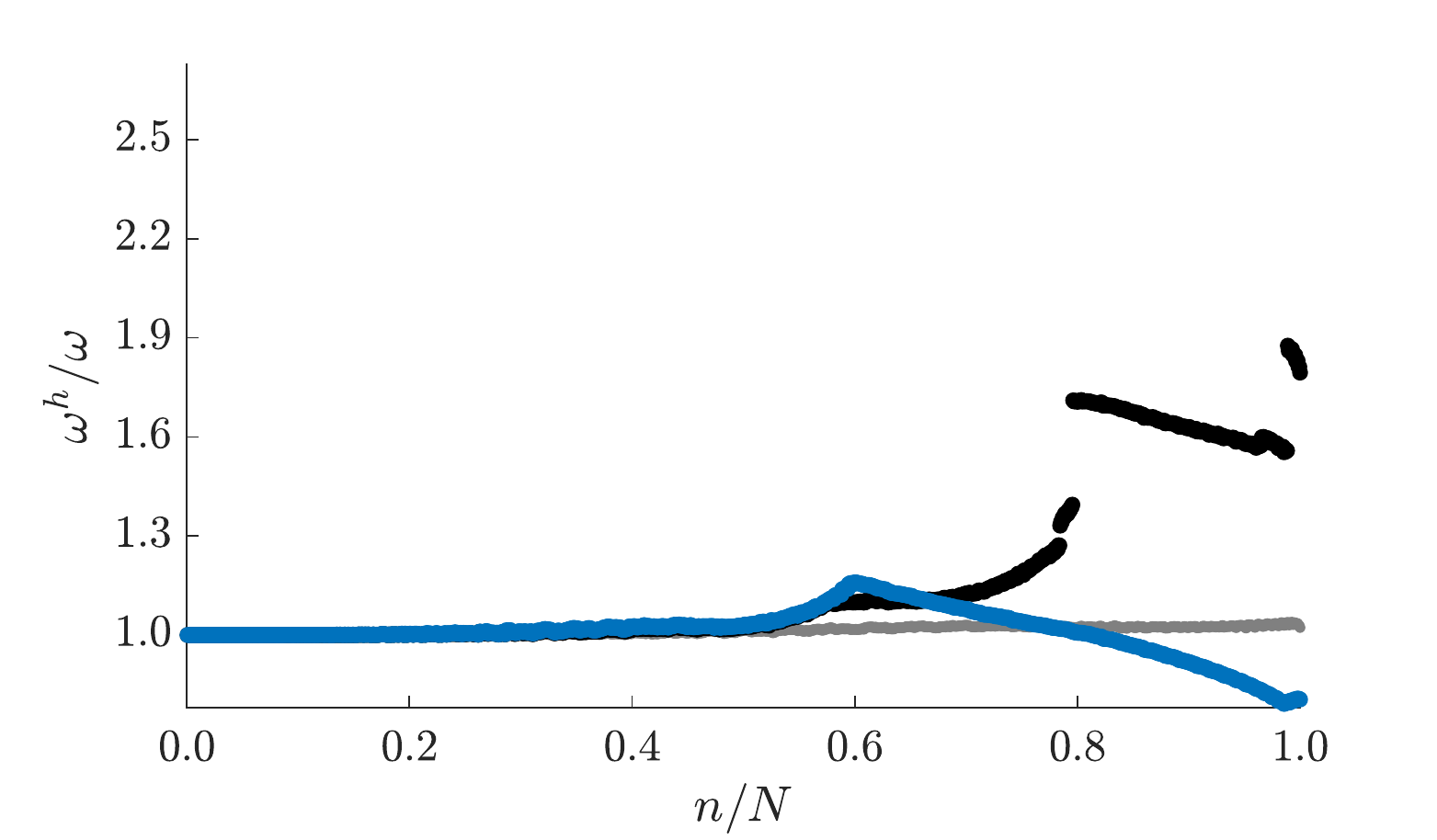} }}

    \subfloat[$p=4$]{{\includegraphics[width=0.5\textwidth]{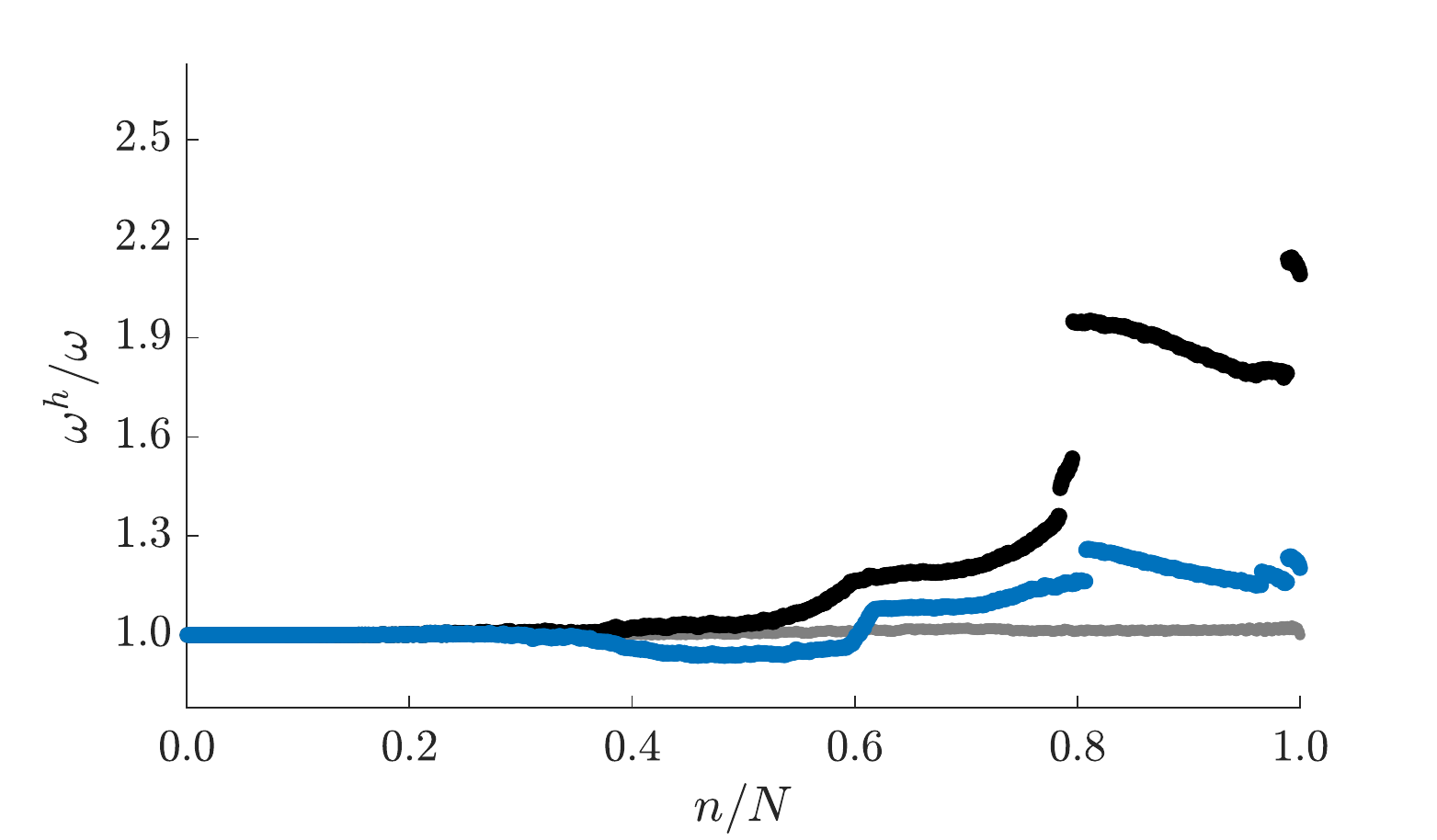} }}
    \subfloat[$p=5$]{{\includegraphics[width=0.5\textwidth]{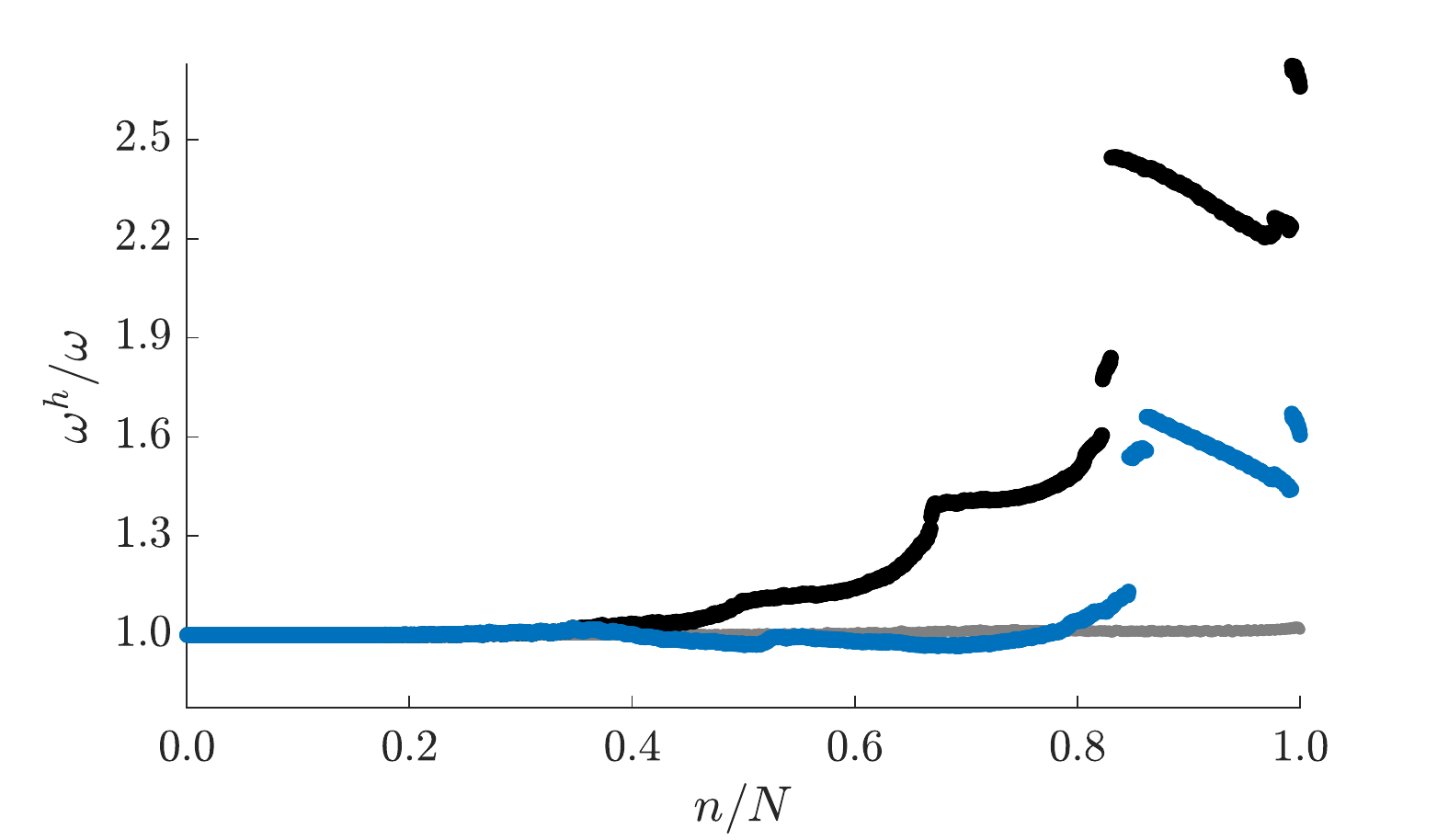} }}
    \vspace{0.2cm}
    \begin{tikzpicture}
    \filldraw[grey1,line width=1pt] (0,0) circle (2pt);
    \filldraw[grey1,line width=1pt] (0,0) node[right]{\footnotesize single-patch};
    \filldraw[black,line width=1pt] (3,0) circle (2pt);
    \filldraw[black,line width=1pt] (3,0) node[right]{\footnotesize multipatch, standard spectrum};
    \filldraw[blue1,line width=1pt] (9,0) circle (2pt);
    \filldraw[blue1,line width=1pt] (9,0) node[right]{\footnotesize multipatch, improved spectrum};
\end{tikzpicture}
    \caption{Normalized frequencies of a freely vibrating \textbf{square membrane with fixed boundary conditions}, computed with \textbf{$5 \times 5$ patches of $C^{p-1}$ B-splines}. Each patch is \textbf{discretized with $5 \times 5$ elements}.}
    \label{fig:normalized_freq2_2nd_2d}
\end{figure}

\begin{figure}[t!]
    \centering
    \subfloat[$p=2$]{{\includegraphics[width=0.5\textwidth]{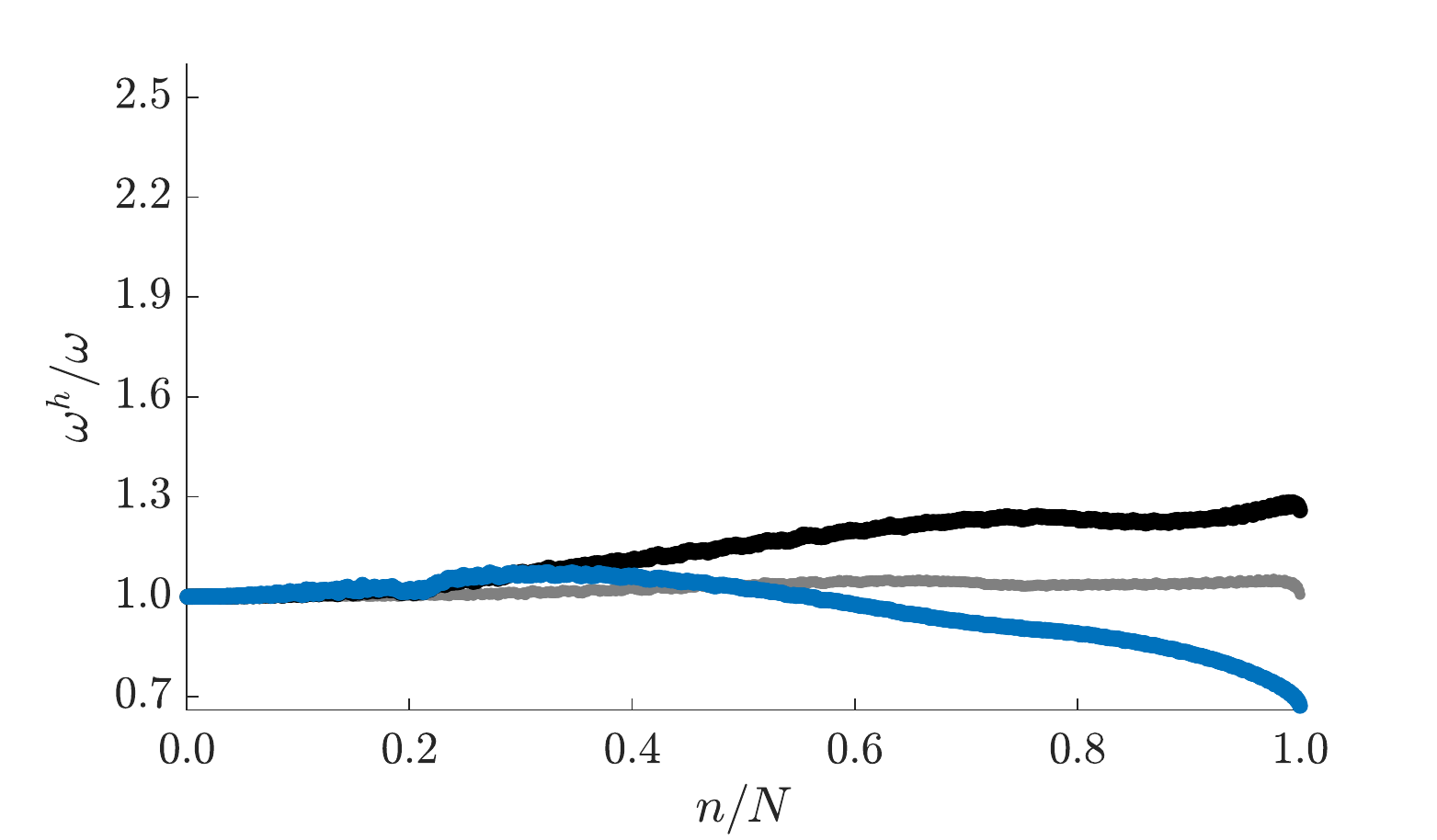} }}
    \subfloat[$p=3$]{{\includegraphics[width=0.5\textwidth]{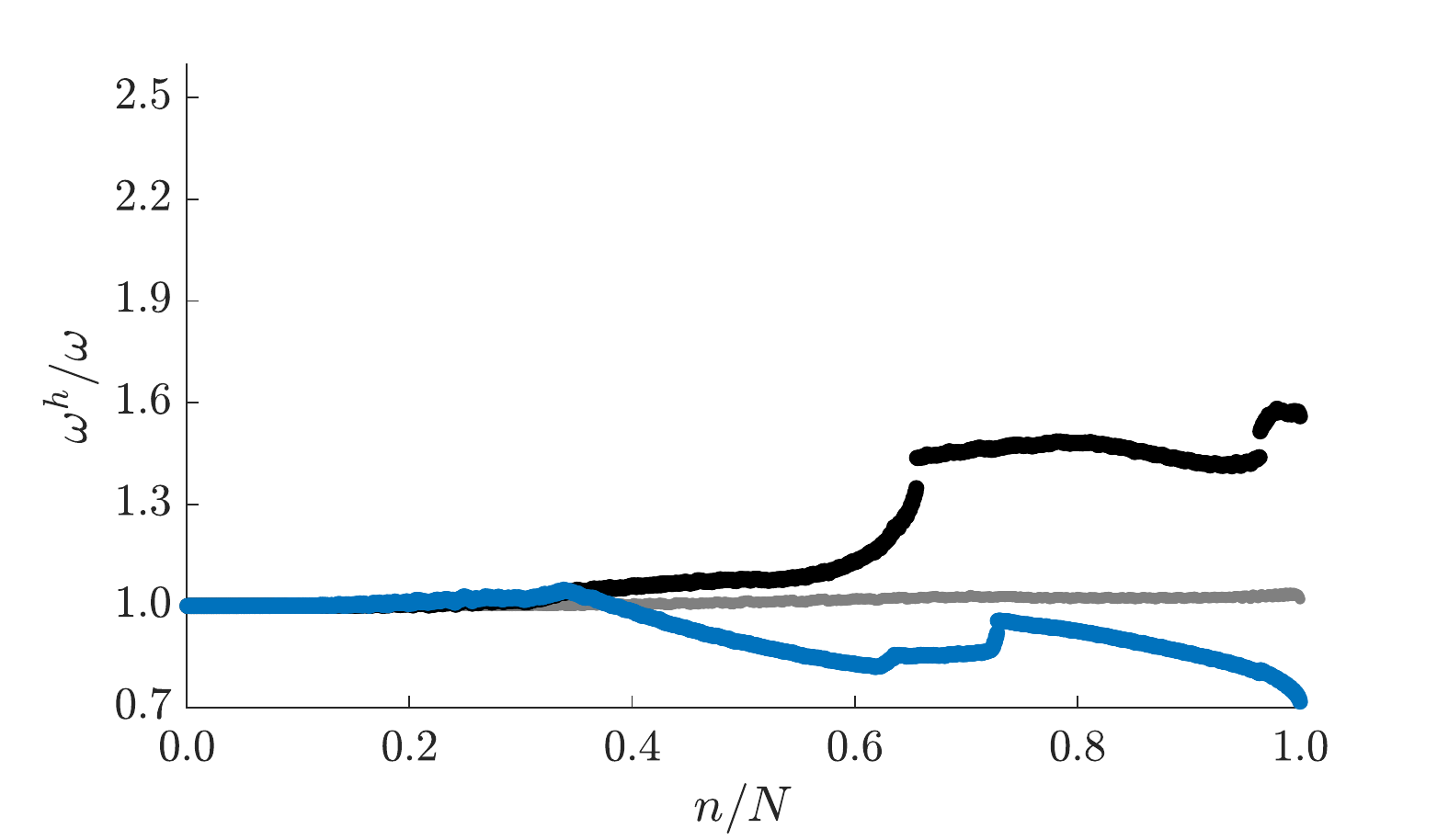} }}

    \subfloat[$p=4$]{{\includegraphics[width=0.5\textwidth]{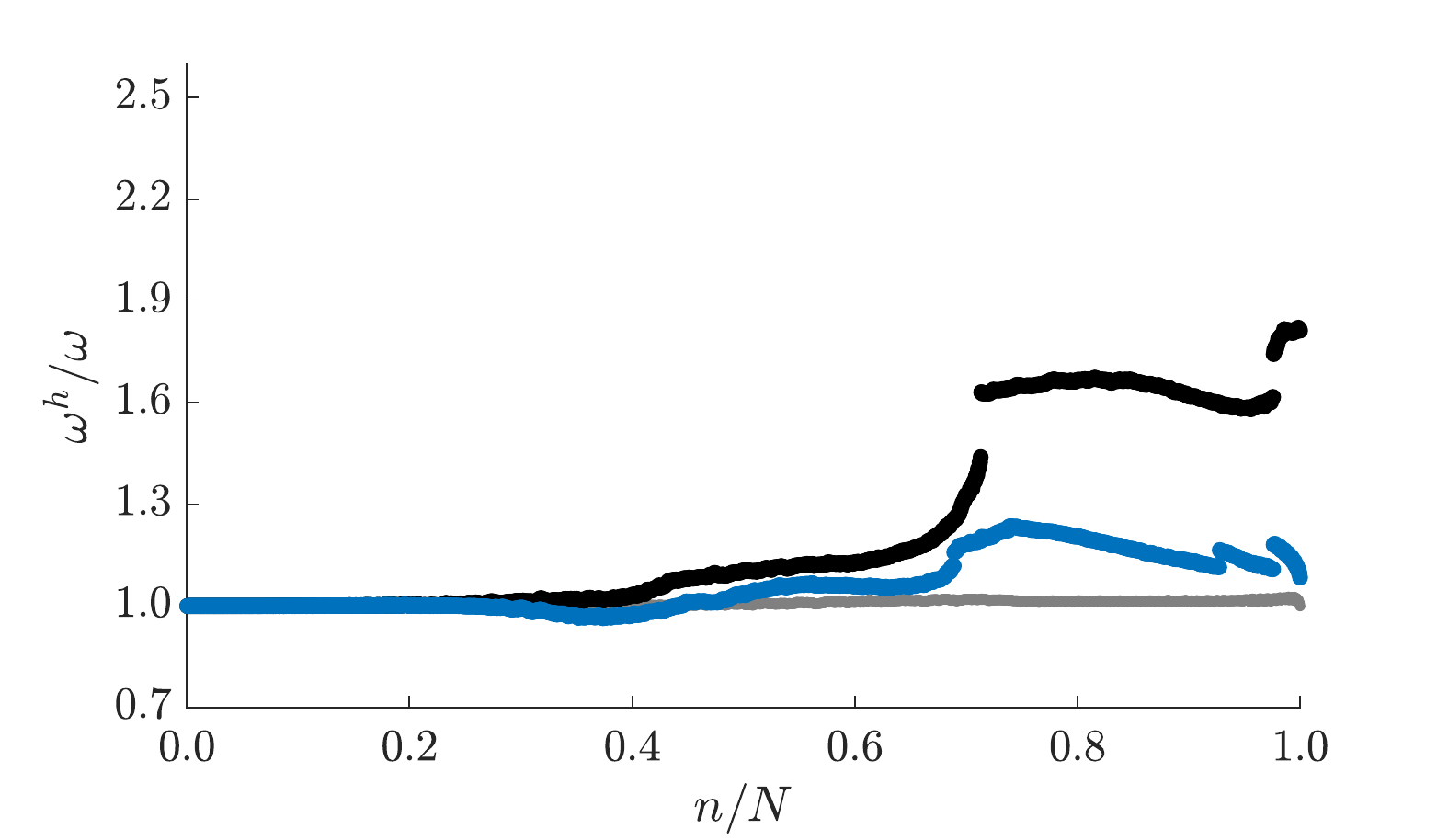} }}
    \subfloat[$p=5$]{{\includegraphics[width=0.5\textwidth]{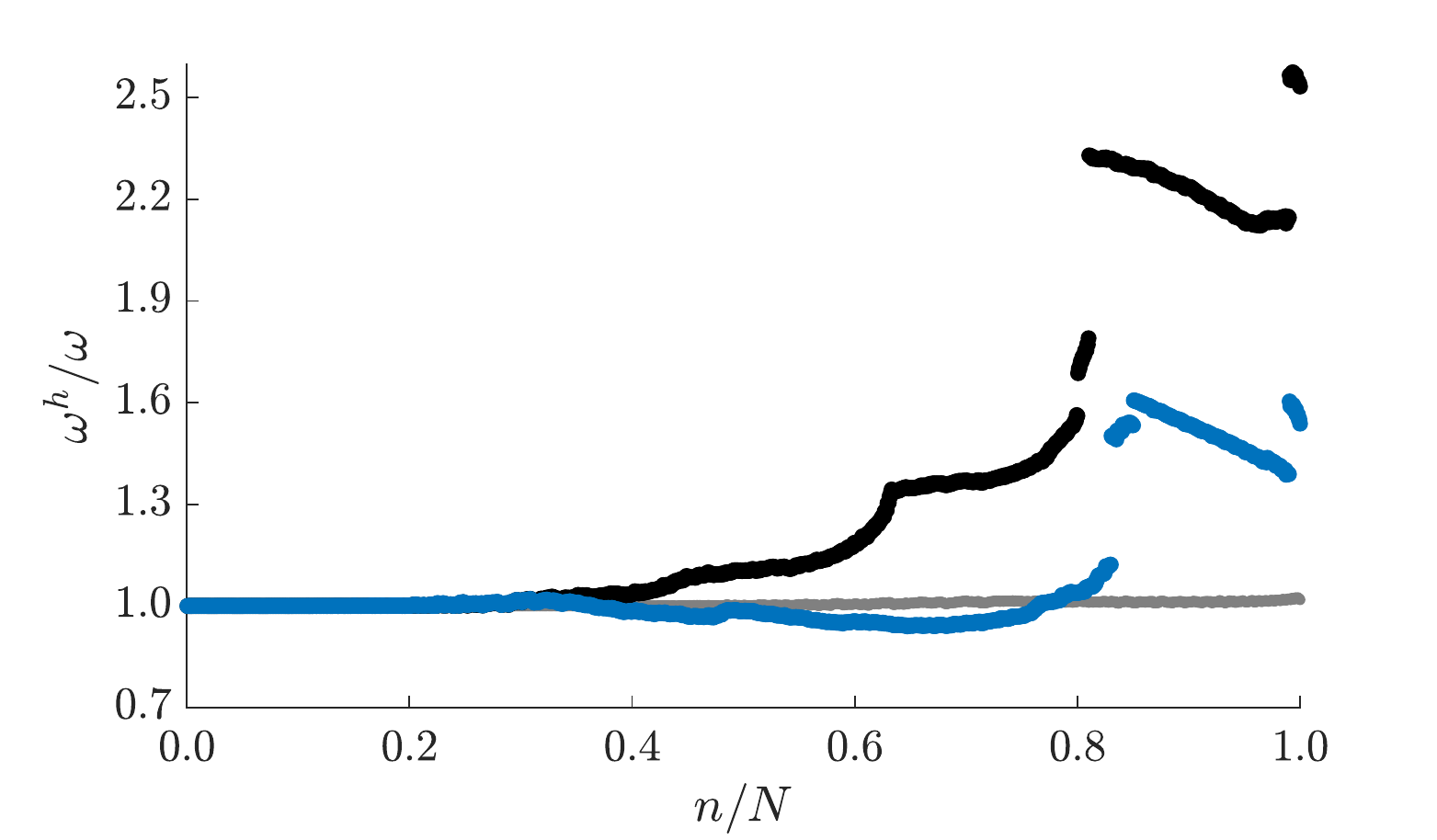} }}
    \vspace{0.2cm}
    \begin{tikzpicture}
    \filldraw[grey1,line width=1pt] (0,0) circle (2pt);
    \filldraw[grey1,line width=1pt] (0,0) node[right]{\footnotesize single-patch};
    \filldraw[black,line width=1pt] (3,0) circle (2pt);
    \filldraw[black,line width=1pt] (3,0) node[right]{\footnotesize multipatch, standard spectrum};
    \filldraw[blue1,line width=1pt] (9,0) circle (2pt);
    \filldraw[blue1,line width=1pt] (9,0) node[right]{\footnotesize multipatch, improved spectrum};
\end{tikzpicture}
    \caption{Normalized frequencies of a freely vibrating \textbf{square membrane with fixed boundary conditions}, computed for the limit case of \textbf{$C^0$ B\'ezier elements} ($\npa = \nele = 15 \times 15$).}
    \label{fig:normalized_freqC0_2nd_2d}
\end{figure}

We consider the free transverse vibration of a square membrane, of unit edge size with fixed boundary conditions, unit membrane stiffness and unit mass. 
We study $C^{p-1}$ B-splines of polynomial degrees $p=2$ through $5$ and $C^0$ patch continuity. We employ multipatch discretizations of $2 \times 2$, $5 \times 5$ patches, and the limit case of $C^0$ B\'ezier elements (one element per patch, $\npa = \nele = 15 \times 15$).
Figures \ref{fig:normalized_freq_2nd_2d}, \ref{fig:normalized_freq2_2nd_2d}, and \ref{fig:normalized_freqC0_2nd_2d} present the normalized frequencies $\omega^h_\noMode / \omega_\noMode$ corresponding to the square membrane discretized with $2 \times 2$, $5 \times 5$, and $15 \times 15$ patches, respectively. 
The inset figures of \ref{fig:normalized_freq_2nd_2d} focus on the upper last twenty percent of the spectra.
Due to their tensor product structure, the spectra of multivariate discretizations exhibit a higher number of outliers than univariate discretizations.

\begin{figure}[h!]
    \centering
    \subfloat[$p=3$]{{\includegraphics[width=0.5\textwidth]{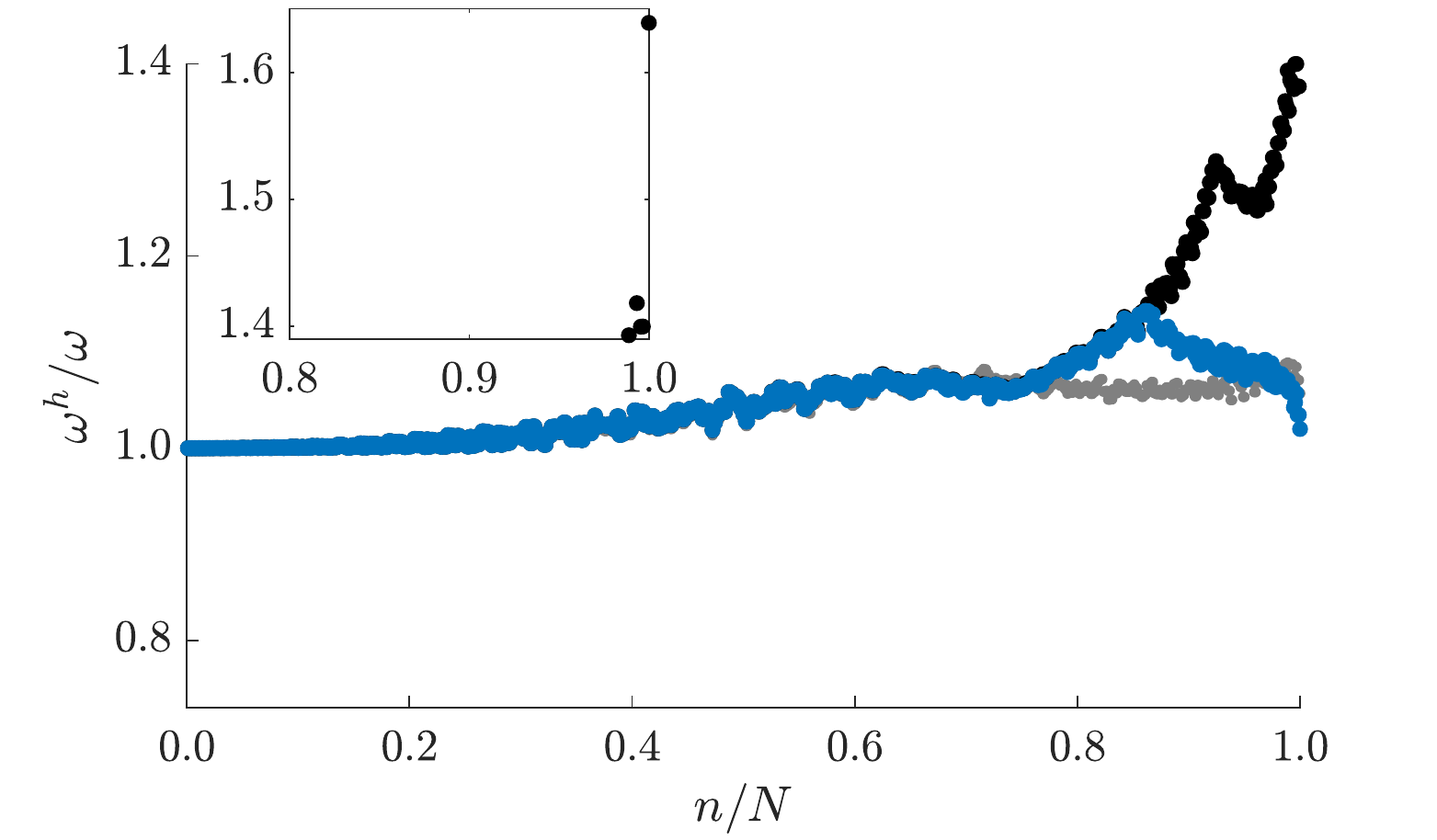} }}
    \subfloat[$p=4$]{{\includegraphics[width=0.5\textwidth]{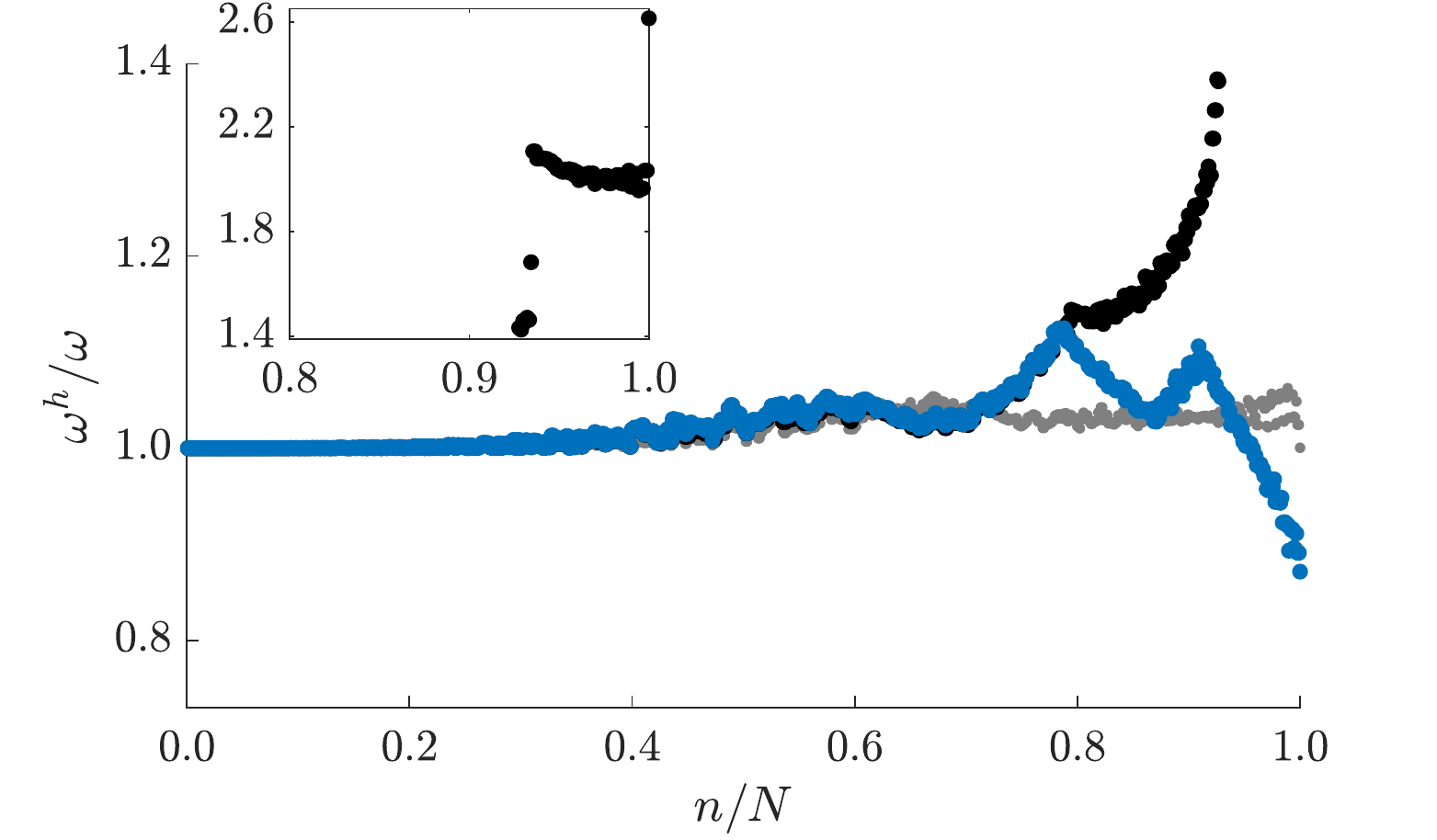} }}

    \subfloat[$p=5$]{{\includegraphics[width=0.5\textwidth]{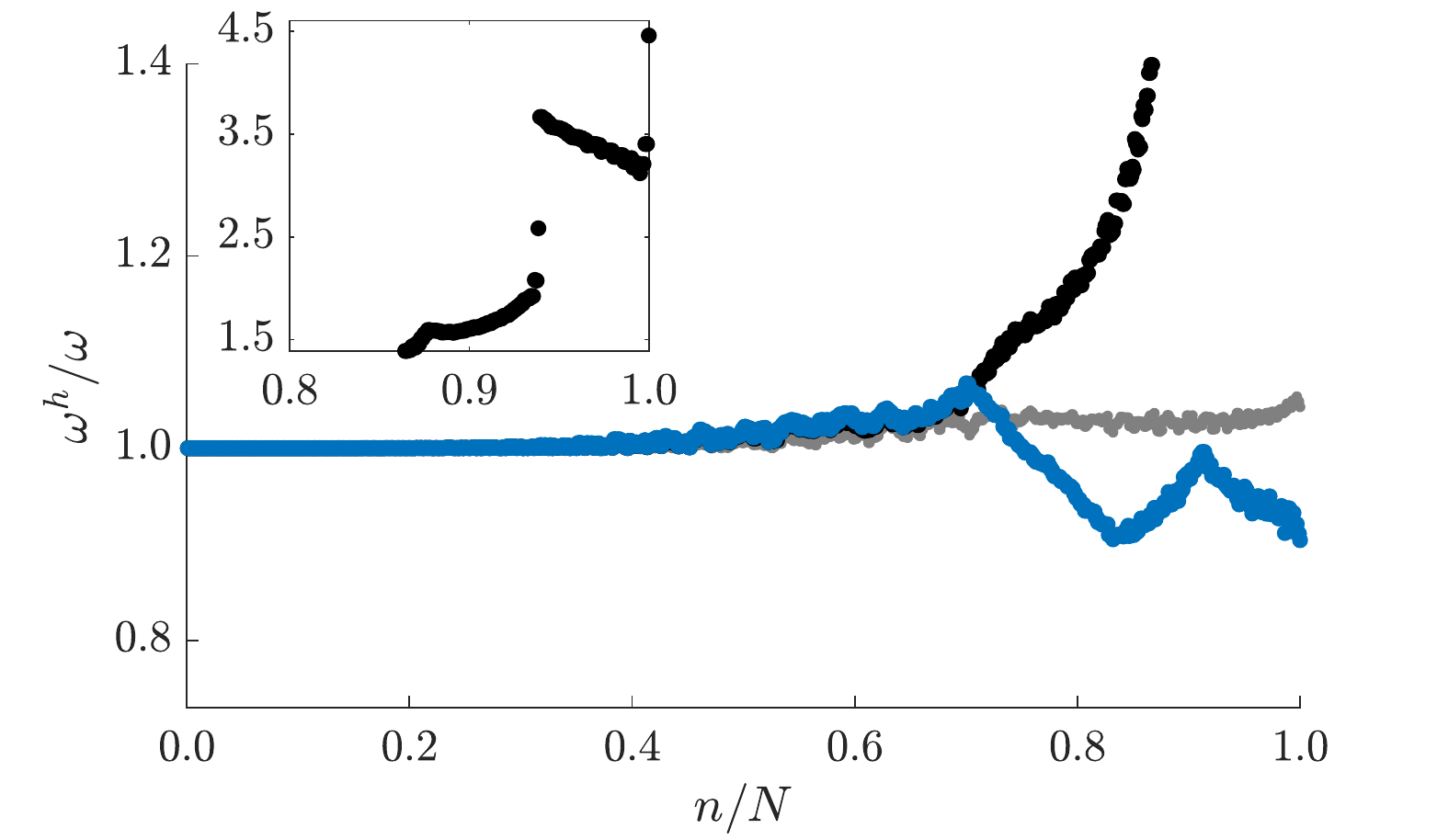} }}
    \subfloat[$p=6$]{{\includegraphics[width=0.5\textwidth]{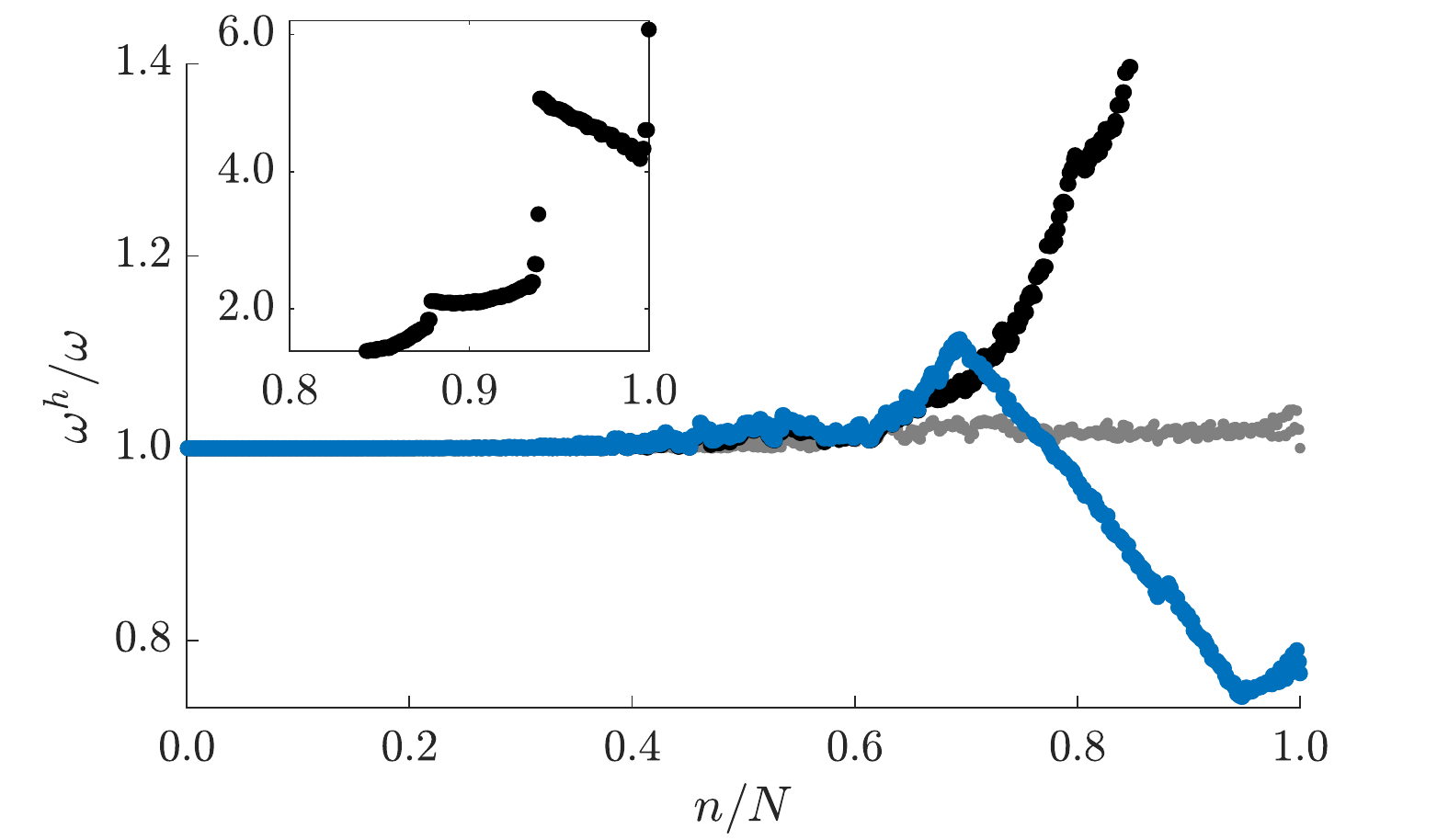} }}
    \vspace{0.2cm}
    \begin{tikzpicture}
    \filldraw[grey1,line width=1pt] (0,0) circle (2pt);
    \filldraw[grey1,line width=1pt] (0,0) node[right]{\footnotesize single-patch};
    \filldraw[black,line width=1pt] (3,0) circle (2pt);
    \filldraw[black,line width=1pt] (3,0) node[right]{\footnotesize multipatch, standard spectrum};
    \filldraw[blue1,line width=1pt] (9,0) circle (2pt);
    \filldraw[blue1,line width=1pt] (9,0) node[right]{\footnotesize multipatch, improved spectrum};
\end{tikzpicture}
    \caption{Normalized frequencies of a freely vibrating \textbf{square plate with simply supported boundary conditions}, computed with \textbf{$2 \times 2$ patches of $C^{p-1}$ B-splines}. Each patch is \textbf{discretized with $15 \times 15$ elements}.}
    \label{fig:normalized_freq_4th_2d}
\end{figure}

\begin{figure}[h!]
    \centering
    \subfloat[$p=3$]{{\includegraphics[width=0.5\textwidth]{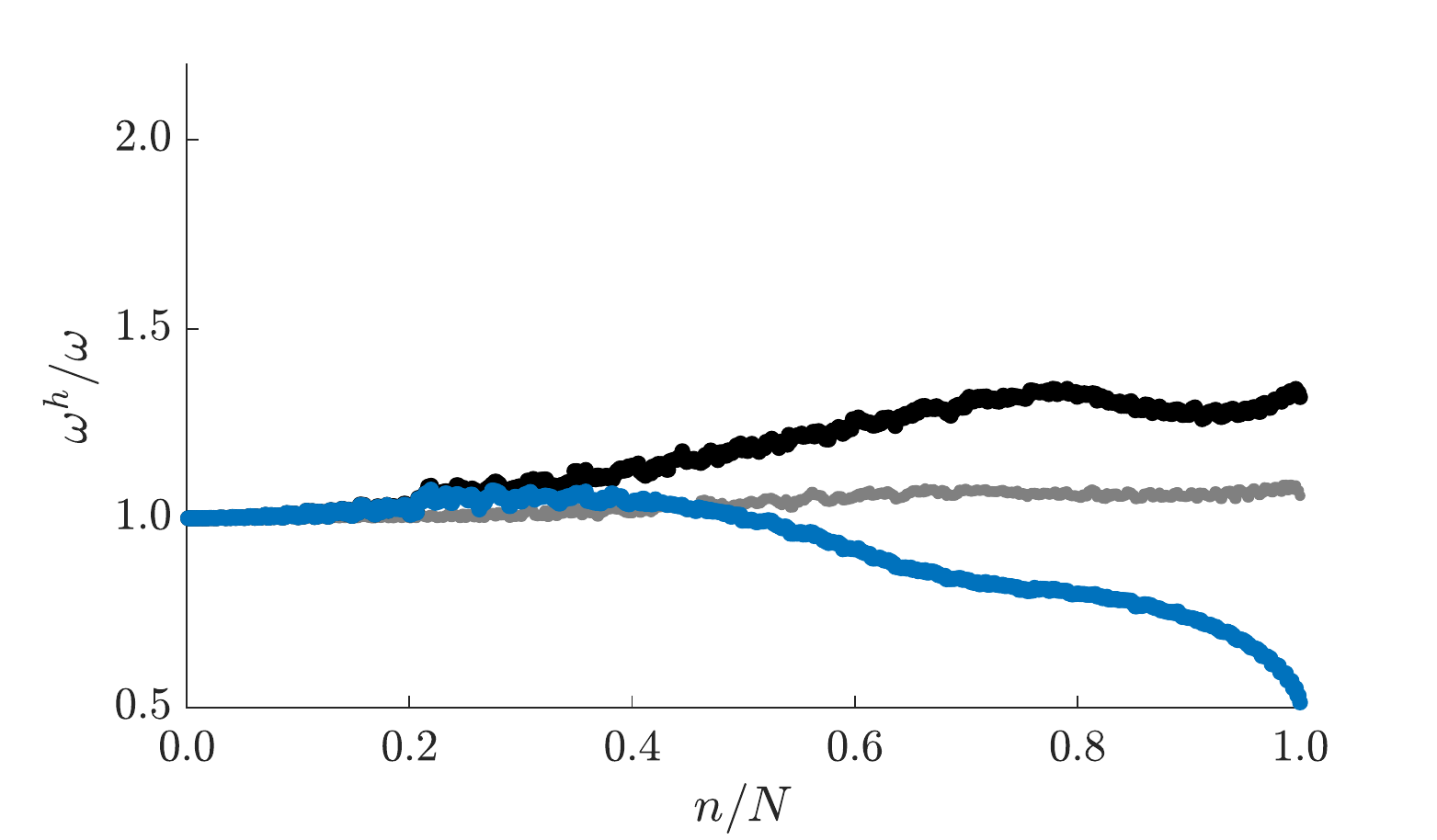} }}
    \subfloat[$p=4$]{{\includegraphics[width=0.5\textwidth]{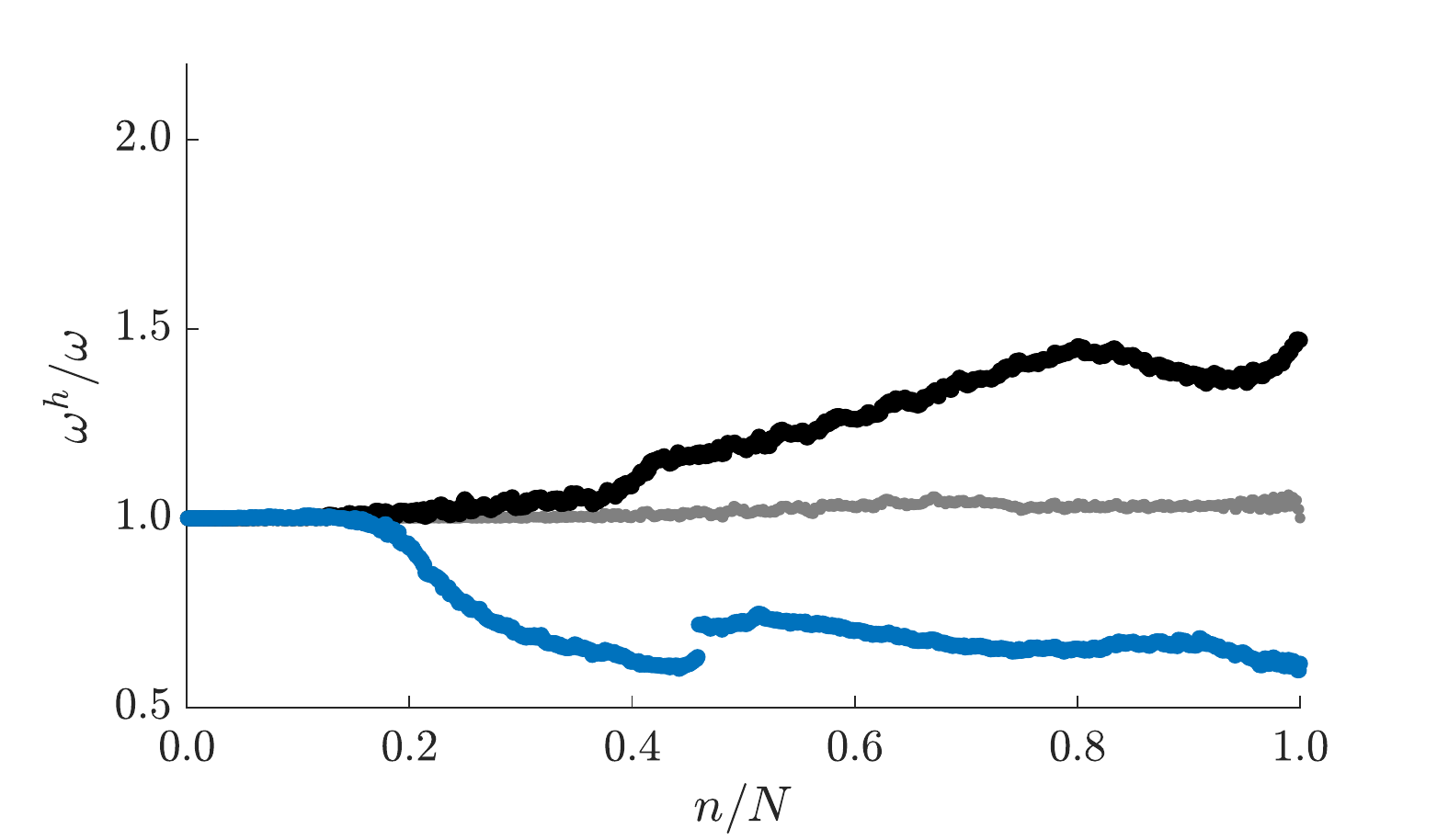} }}

    \subfloat[$p=5$]{{\includegraphics[width=0.5\textwidth]{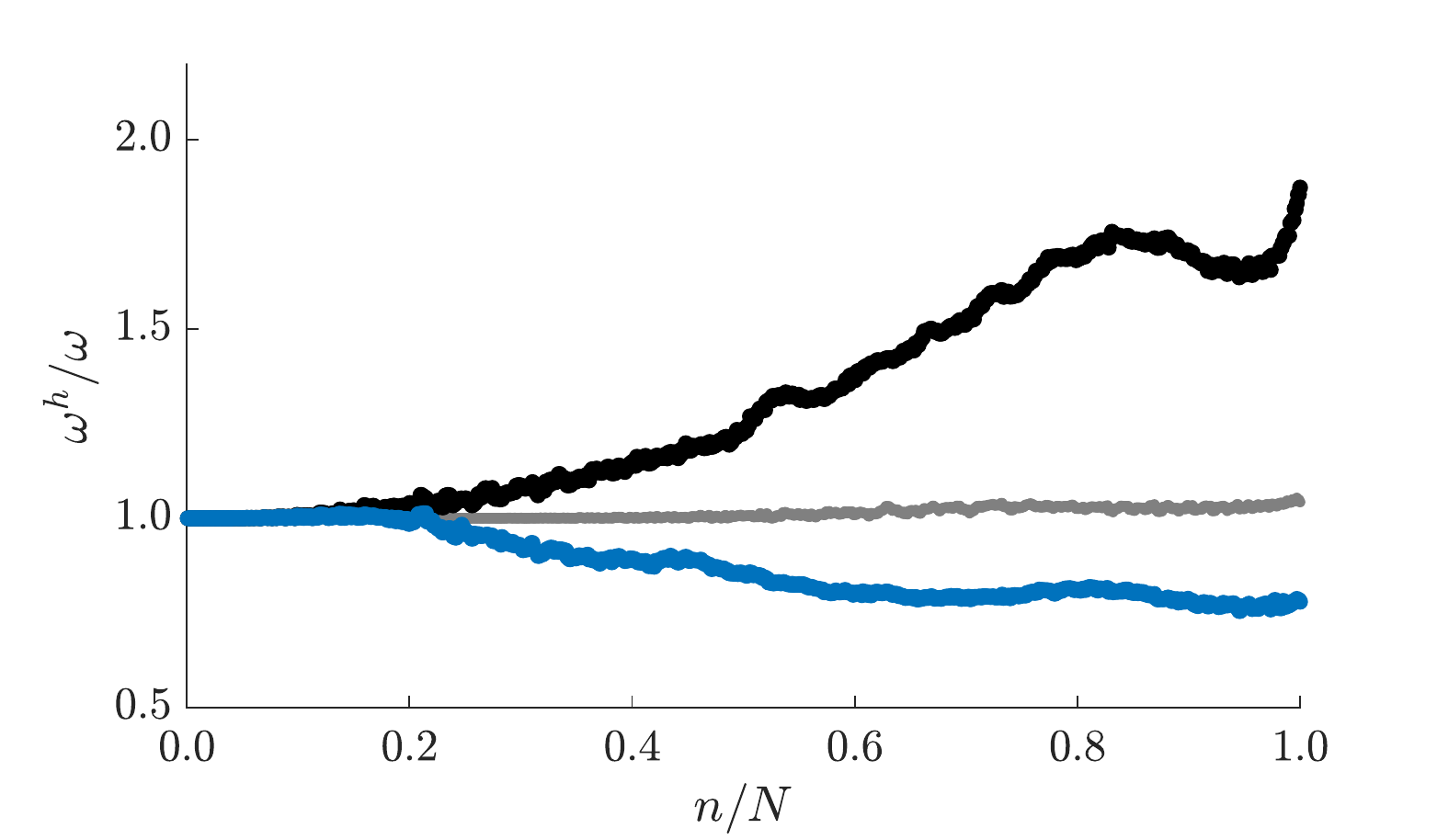} }}
    \subfloat[$p=6$]{{\includegraphics[width=0.5\textwidth]{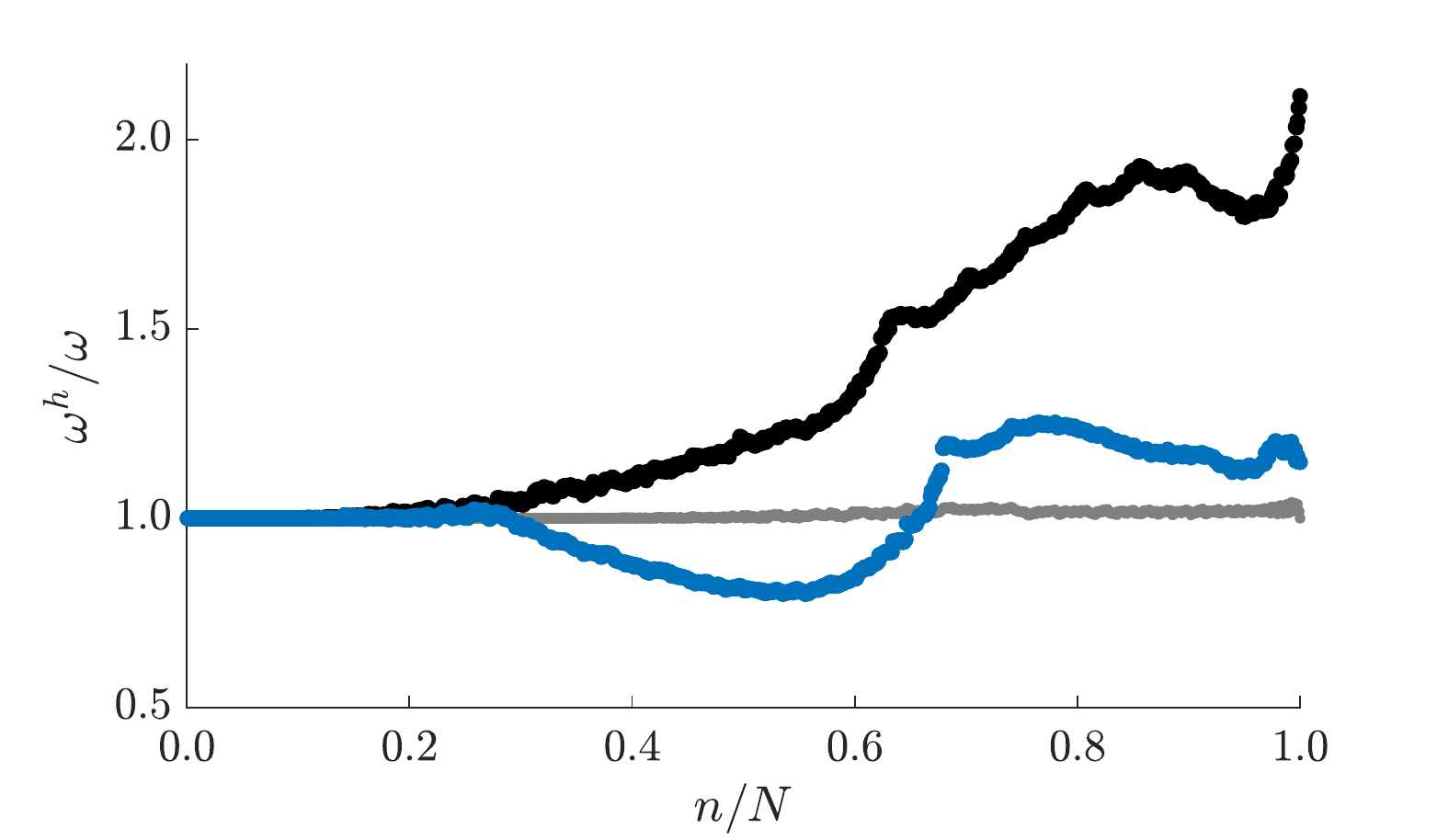} }}
    \vspace{0.2cm}
    \begin{tikzpicture}
    \filldraw[grey1,line width=1pt] (0,0) circle (2pt);
    \filldraw[grey1,line width=1pt] (0,0) node[right]{\footnotesize single-patch};
    \filldraw[black,line width=1pt] (3,0) circle (2pt);
    \filldraw[black,line width=1pt] (3,0) node[right]{\footnotesize multipatch, standard spectrum};
    \filldraw[blue1,line width=1pt] (9,0) circle (2pt);
    \filldraw[blue1,line width=1pt] (9,0) node[right]{\footnotesize multipatch, improved spectrum};
\end{tikzpicture}
    \caption{Normalized frequencies of a freely vibrating \textbf{square plate with simply supported boundary conditions}, computed for the limit case of \textbf{$C^1$ B\'ezier elements} ($\npa = \nele = 15 \times 15$).}
    \label{fig:normalized_freqC1_4th_2d}
\end{figure}

We then consider the free vibration of a square plate structure of unit edge size with simply supported boundary conditions, unit bending stiffness and unit mass. 
We use $C^{p-1}$ B-splines of polynomial degrees $p=3$ through $6$ and $C^1$ patch continuity. 
We employ multipatch discretizations of $2 \times 2$ and the limit case of $C^1$ B\'ezier elements (same number of patches and elements, $\npa = \nele = 15 \times 15$). 
Figures \ref{fig:normalized_freq_4th_2d} and \ref{fig:normalized_freqC1_4th_2d} illustrate the normalized frequencies, $\omega^h_\noMode / \omega_\noMode$, corresponding to the square plate discretized with $2 \times 2$ and $15 \times 15$ patches, respectively. 
These results confirm that our approach \eqref{eq:vdgep_perturbed2} in combination with \changed{Algorithm \ref{fig:parameter_estimation}} reduces the outlier frequencies of multipatch discretizations effectively for both second- and fourth-order problems, without negatively affecting lower frequencies, and works well for different polynomial degrees and patch configurations.

\section{Application in explicit dynamics}\label{sec:results-dyn}

The critical time-step size in explicit dynamics calculations is inversely proportional to the maximum discrete eigenfrequency \cite{hughes_finite_2003}. Significantly overestimated outlier frequencies therefore negatively affect the critical time-step size, and hence the computational cost of explicit dynamics calculations. It can be thus expected that the approaches presented in this work are able to effectively improve this issue for multipatch isogeometric discretizations, which we will illustrate in the following.

\subsection{Semidiscrete formulation}

In this section, we consider the semi-discrete form \eqref{deom} of an free-vibrating, undamped linear structural system, which can be expressed as follows:
\begin{align}
    \mat{M} \, \frac{d^2}{d\,t^2} \, \mat{u}^h(t) \; = \; \underbrace{- \, \mat{K} \, \mat{u}^h(t)}_{\mat{F}_{\text{int}}} \, , \label{deom2}
\end{align}
where $\mat{F}_{\text{int}}$ is the vector of internal forces. 
Using the proposed approach in Section \ref{sec:pragmatic_approach} leads to the following semi-discrete form:
\begin{align}
    \left(\, \mat{M} + \alpha \, \mat{K}_\Gamma \, \right) \, \frac{d^2}{d\,t^2} \, \mat{u}^h(t) \; = \; \underbrace{- \, \left( \, \mat{K} + \beta \, \mat{K}_\Gamma \, \right) \, \mat{u}^h(t)}_{\tilde{\mat{F}}_{\text{int}}} \, , \label{deom3}
\end{align}
where $\tilde{\mat{F}}_{\text{int}}$ is the perturbed vector of internal forces. Both the left- and right-hand sides of the standard formulation \eqref{deom2} are modified using our approach. This is not the case when using the mass scaling approach where only the mass matrix on the left-hand side is modified by adding artificial mass terms \cite{olovsson_mass_scaling_2005,tkachuk_mass_scaling_2014, schaeuble_mass_scaling_2017,gonzalez_mass_tailoring_2020}. 
The perturbation matrix $\mat{K}_\Gamma$ is computed once and hence no reassembly of the mass matrix is needed, so that this does not increase the computational cost of explicit dynamics calculations.

\begin{remark}
    Explicit dynamics applications typically involve the use of a lumped mass matrix in combination with an explicit time integration scheme. There are currently no widely accepted mass lumping techniques that maintain higher-order spatial accuracy. To demonstrate that our methodology maintains higher-order spatial accuracy, we employ the consistent mass matrix in all subsequent computations.
\end{remark}

\subsection{Optimum spatial accuracy}

\begin{figure}[h!]
    \centering
    \subfloat[Coarsest \Bezier mesh with boundary and interface constraints]{{\includegraphics[width=0.42\textwidth]{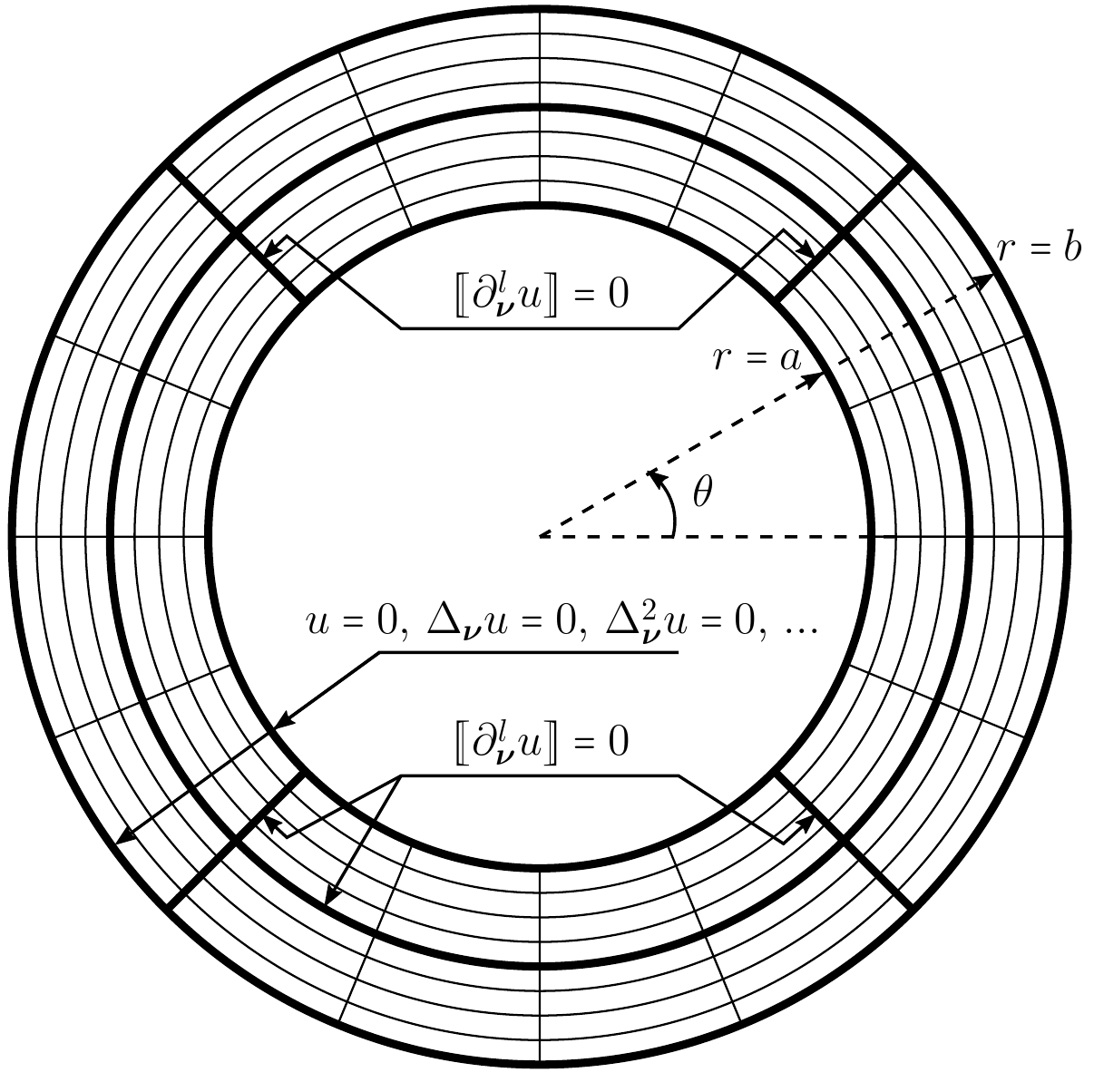} }}\hspace{0.5cm}
    \subfloat[Initial displacement field $u(r,\theta,0)$]{{\includegraphics[width=0.42\textwidth]{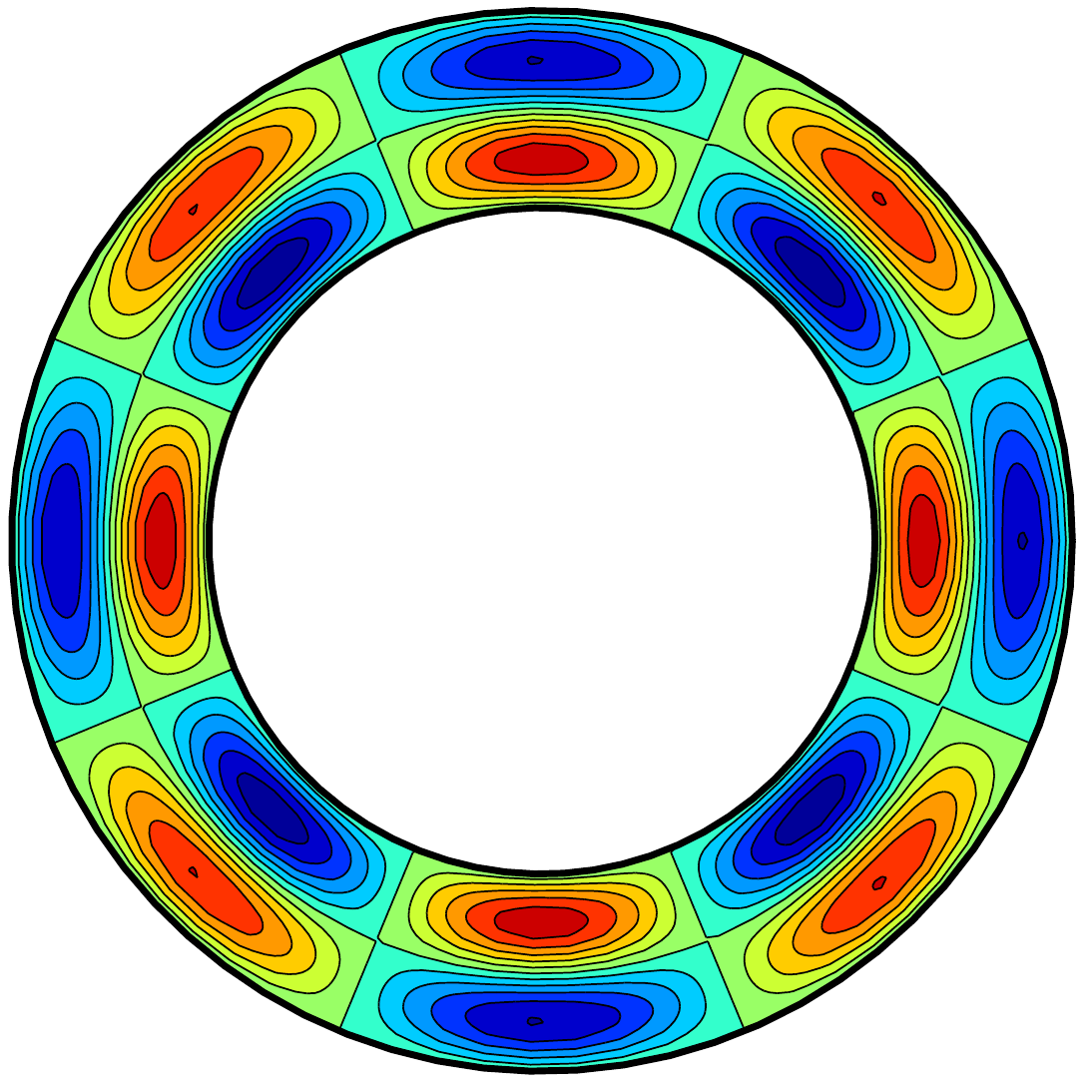} }}
    \caption{Transient model problem on an annulus.}
    \label{fig:annulus_initial}
\end{figure}

We consider the free vibration of the annular membrane problem illustrated in Figure \ref{fig:annulus_initial}a, see also \cite{hiemstra_outlier_2021}, where boundary conditions are fixed \changed{along the inner radius $a$ and outer radius $b$}.
We choose the following displacement solution $u$ that satisfies the differential equation \eqref{eom}:
\begin{align}
    u(r,\theta,t) = J_4(r) \; \cos(\lambda_2\,t) \; \cos(4\,\theta) \, , \label{eq:annulus_analytic}
\end{align}
with radial coordinate $r$, angular coordinate $\theta$ and time $t$. 
Here, $J_4(r)$ denotes the $4^{\text{th}}$ Bessel function of the first kind and $\lambda_k$, $k=1,2,\ldots$ denote its positive zeros.
We choose the second and fourth zeros as the inner and outer radii of the annulus, respectively, i.e. $a = \lambda_2 \approx 11.065$ and $b = \lambda_4 \approx 17.616$.
The analytical solution \eqref{eq:annulus_analytic} at time $t=0$ plotted in Figure \ref{fig:annulus_initial}b defines the initial displacement field $u(r,\theta,0)$ .

We study the \changed{semi-discrete form} \eqref{deom3} of the annular membrane in explicit dynamics.
We employ a multivariate spline space that is an extension of the univariate spaces of the angular coordinate \changed{($\theta \in [0,2\pi]$) and radial coordinate ($r \in [a,b]$)}. 
The univariate spline spaces are free of boundary outliers \changed{due to} outlier removal boundary constraints \cite{hiemstra_outlier_2021}.
In the spline space of the angular coordinate, we also build in the periodic end-conditions.
We consider $C^{p-1}$ B-splines of different polynomial degrees $p=2$ through $5$ for multipatch discretizations with $2$ patches in the radial direction, and $4$ patches in the angular direction, and $C^0$ patch continuity (see Figure \ref{fig:annulus_initial}a).
We apply the proposed approach to enforce the $C^{p-1}$ continuity constraints \eqref{eq:continuity_constraint} at patch interfaces (see Figure \ref{fig:annulus_initial}a).
We perform uniform mesh refinement of each patch \changed{with $4,8,16$ elements, i.e. $\nele = 8,16,32$ elements in the radial coordinate, and $16,32,64$ elements in the angular coordinate.} 
We simulate up to a final time of $T = 2\pi / \lambda_2$ which is one full period of the periodic function $u(r,\theta,t)$ \eqref{eq:annulus_analytic}.
For time integration, we apply the central difference method \cite{hughes_finite_2003}.
To verify that \changed{spatial accuracy is maintained}, we choose a small, order-dependent time step of $\Delta t = \left(p/(2 \nele)\right)^p$.

\begin{figure}[h!]
    \centering
    \includegraphics[width=0.7\textwidth]{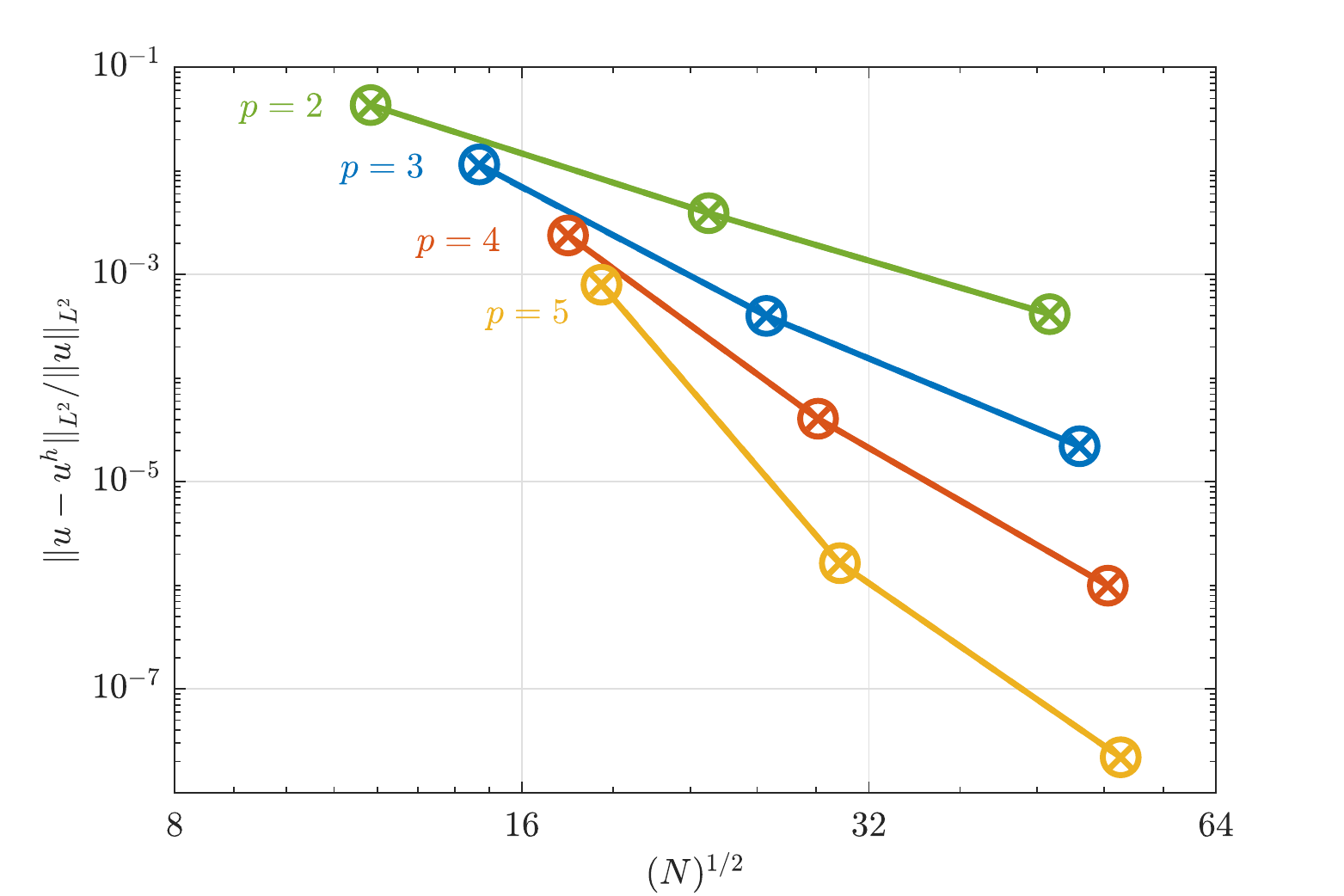} 
    \vspace{0.3cm}
    \begin{tikzpicture}
    \filldraw[black,line width=1pt, solid] (0.0,0) -- (0.2,0);
    \filldraw[black,line width=1pt] (0.3,0) [fill=none] circle (2pt);
    \filldraw[black,line width=1pt, solid] (0.4,0) -- (0.6,0);
    \filldraw[black,line width=1pt] (0.6,0) node[right]{\footnotesize standard spectrum};
    \filldraw[black,line width=1pt, dashed] (6,0) -- (6.6,0);
    \filldraw[black,line width=1pt] (6.0,0) node[right]{\footnotesize $\boldsymbol{\bigtimes}$};
    \filldraw[black,line width=1pt] (6.6,0) node[right]{\footnotesize improved spectrum};
\end{tikzpicture}
    \caption{Convergence of the relative $L^2$ error in the vertical displacement field $u$ of \textbf{the annular membrane in Figure \ref{fig:annulus_initial}}, \changed{computed with a small, order-dependent time step of $\Delta t = \left(p/(2 \nele)\right)^p$}.}
    \label{fig:convergence_annulus}
\end{figure}

Figure \ref{fig:convergence_annulus} compares the convergence behavior of the $L^2$ error in the discrete displacement field $u^h(r,\theta,t)$, when we do standard analysis (circle) and when we do analysis with the proposed approach based on perturbed eigenvalue problems (cross). 
We observe that the analysis with our perturbation approach maintains optimum spatial accuracy.

\subsection{Critical time-step size}

\begin{figure}[h!]
    \centering
    \subfloat[Membrane]{{\includegraphics[width=0.5\textwidth]{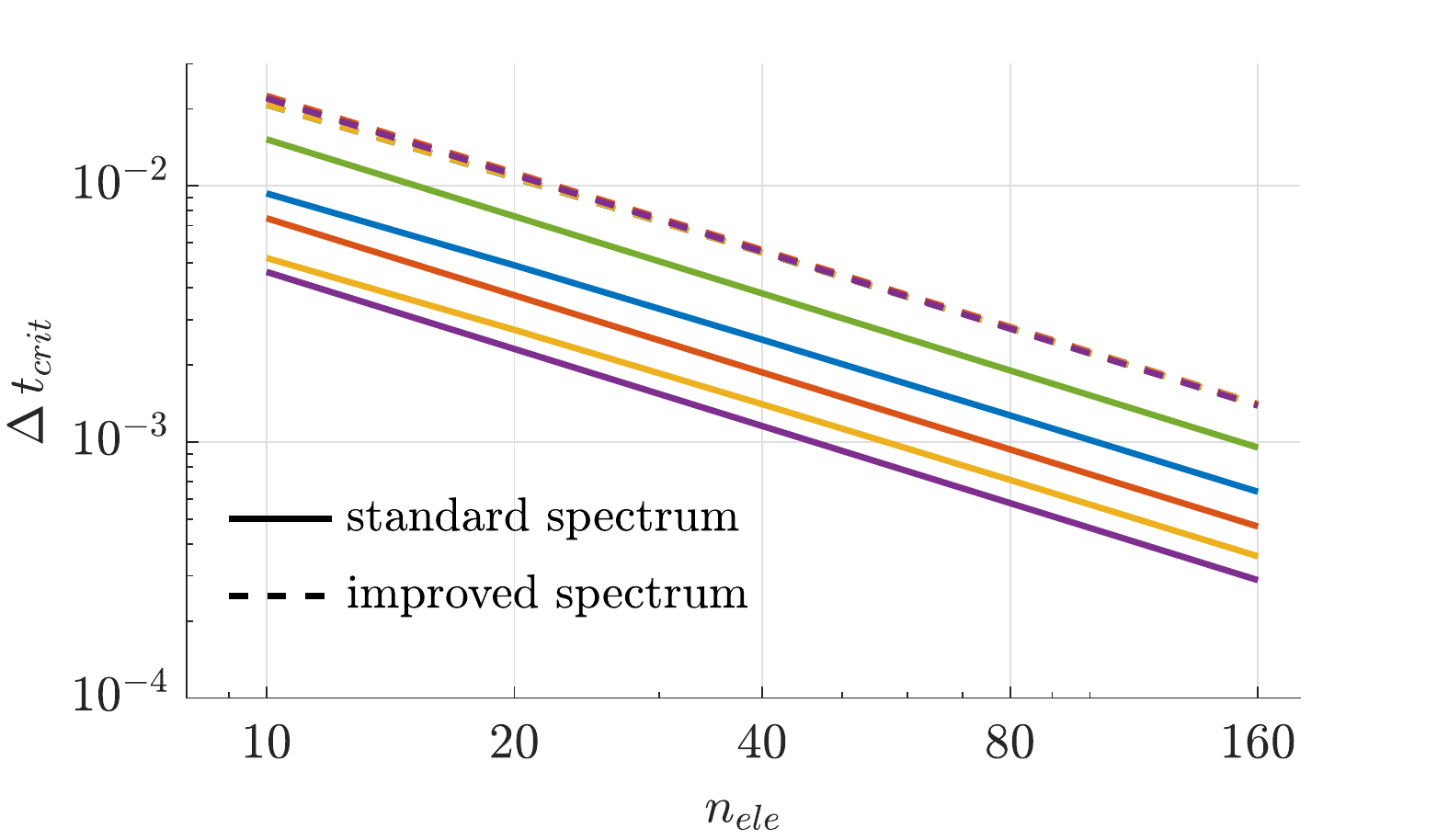} }}
    \subfloat[Plate]{{\includegraphics[width=0.5\textwidth]{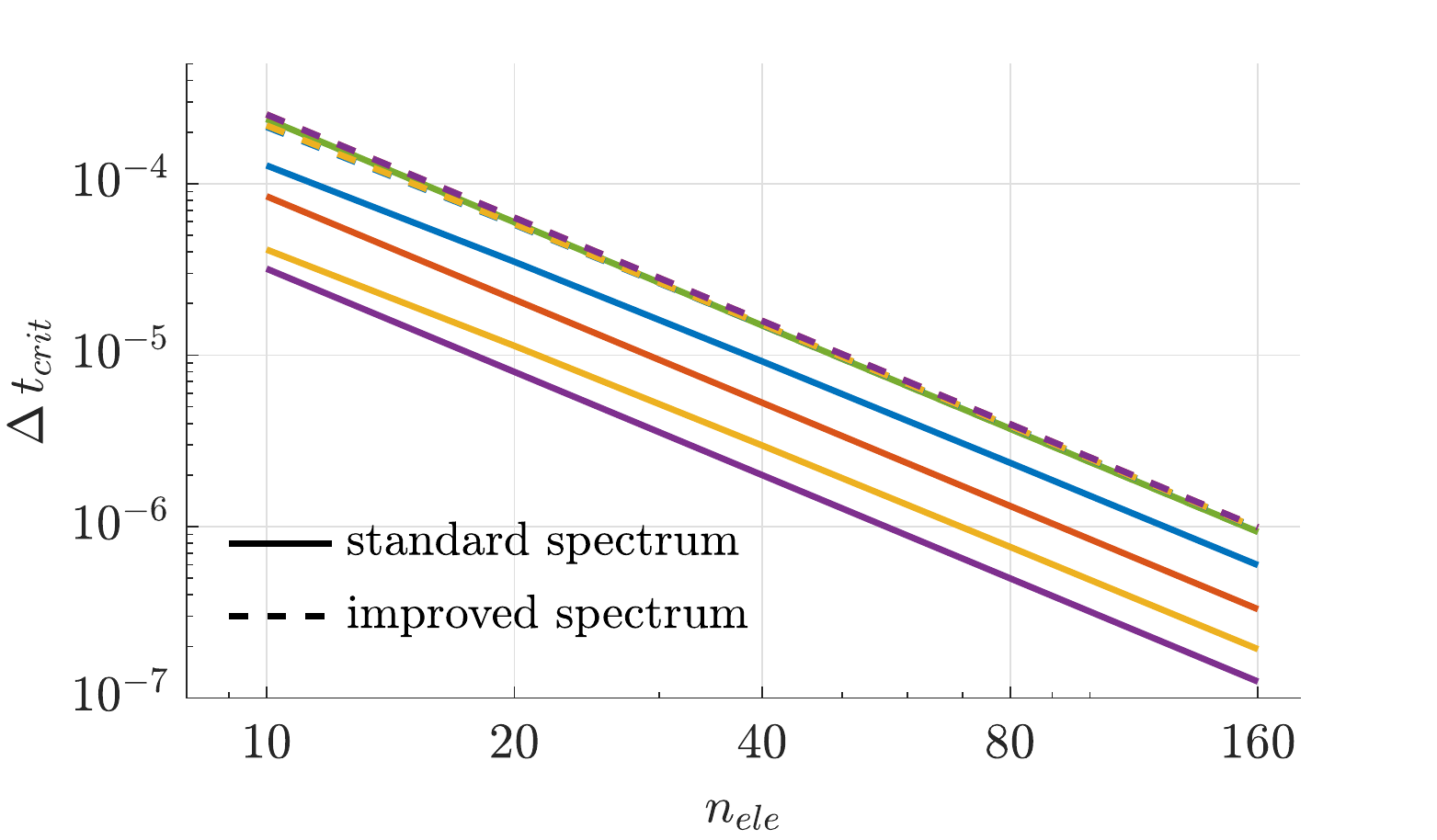} }}
    \vspace{0.3cm}
    \begin{tikzpicture}
    \filldraw[green1,line width=1pt, solid] (0,0) -- (0.6,0);
    \filldraw[green1,line width=1pt] (0.6,0) node[right]{\footnotesize $p=2$};	
    \filldraw[blue1,line width=1pt, solid] (3,0) -- (3.6,0);
    \filldraw[blue1,line width=1pt] (3.6,0) node[right]{\footnotesize $p=3$};
    \filldraw[red1,line width=1pt, solid] (6,0) -- (6.6,0);
    \filldraw[red1,line width=1pt] (6.6,0) node[right]{\footnotesize $p=4$};
    \filldraw[yellow1,line width=1pt, solid] (9,0) -- (9.6,0);
    \filldraw[yellow1,line width=1pt] (9.6,0) node[right]{\footnotesize $p=5$};
    \filldraw[purple1,line width=1pt, solid] (12,0) -- (12.6,0);
    \filldraw[purple1,line width=1pt] (12.6,0) node[right]{\footnotesize $p=6$};
\end{tikzpicture}
    \caption{Critical time step size in explicit dynamics of a \textbf{square membrane and plate} with \textbf{fixed and simply supported boundary conditions}, respectively, as a function of the mesh $N = 2\nele \times 2\nele$ using $2 \times 2$ patches, resulting from standard and improved spectrum.}
    \label{fig:time_step}
\end{figure}

\begin{figure}[h!]
    \centering
    \subfloat[Membrane]{{\includegraphics[width=0.5\textwidth]{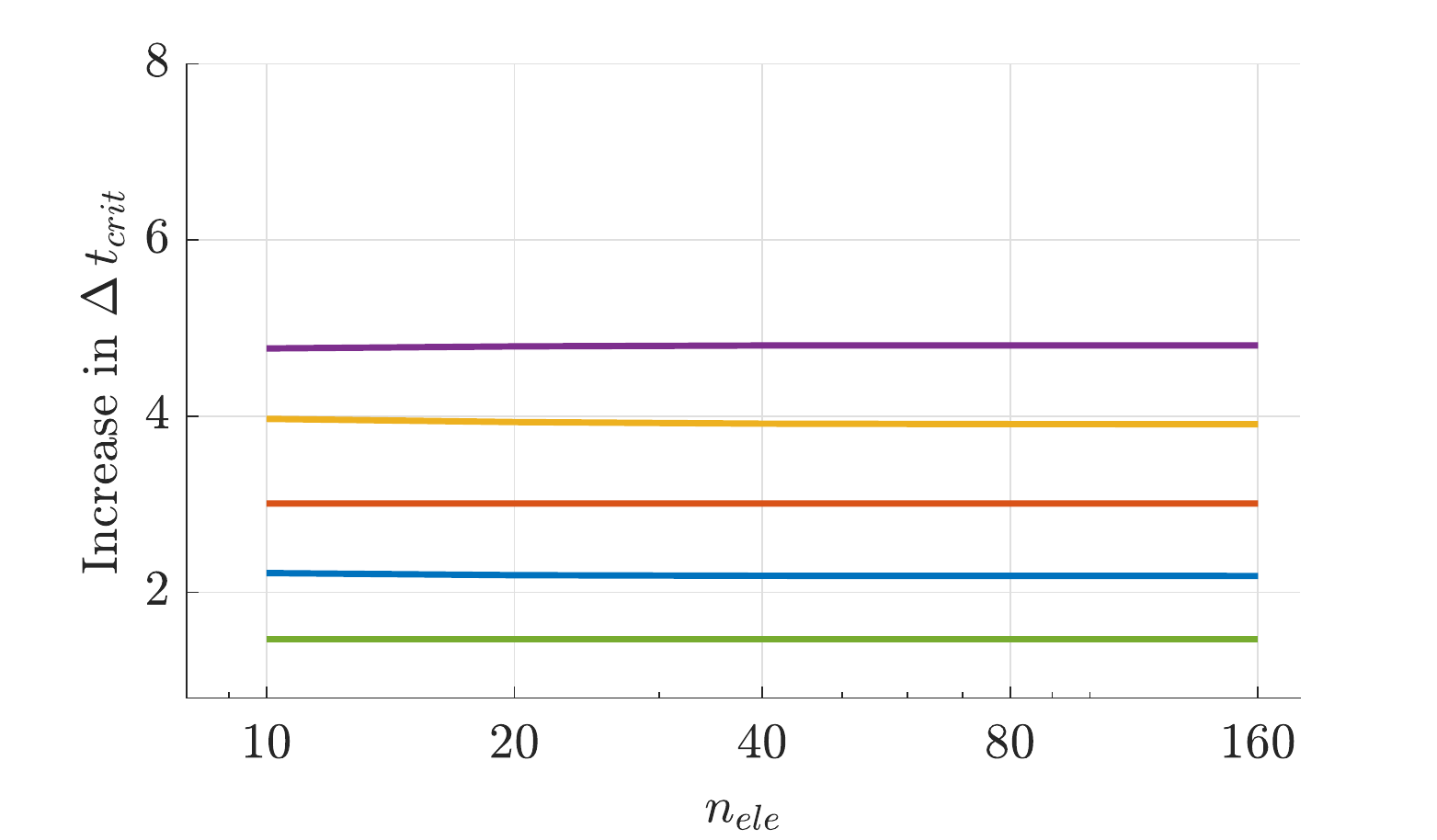} }}
    \subfloat[Plate]{{\includegraphics[width=0.5\textwidth]{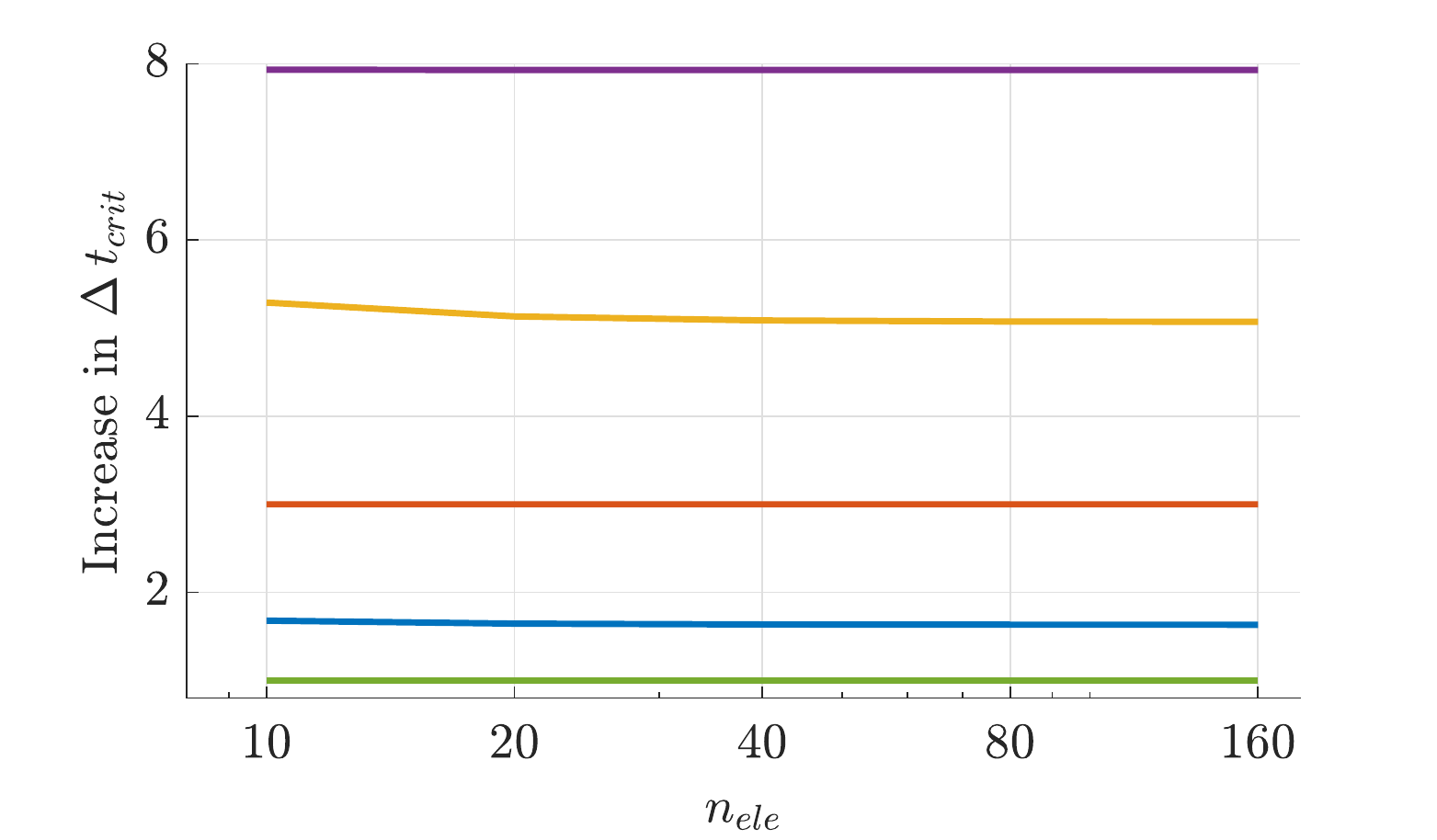} }}
    \vspace{0.3cm}
    \begin{tikzpicture}
    \filldraw[green1,line width=1pt, solid] (0,0) -- (0.6,0);
    \filldraw[green1,line width=1pt] (0.6,0) node[right]{\footnotesize $p=2$};	
    \filldraw[blue1,line width=1pt, solid] (3,0) -- (3.6,0);
    \filldraw[blue1,line width=1pt] (3.6,0) node[right]{\footnotesize $p=3$};
    \filldraw[red1,line width=1pt, solid] (6,0) -- (6.6,0);
    \filldraw[red1,line width=1pt] (6.6,0) node[right]{\footnotesize $p=4$};
    \filldraw[yellow1,line width=1pt, solid] (9,0) -- (9.6,0);
    \filldraw[yellow1,line width=1pt] (9.6,0) node[right]{\footnotesize $p=5$};
    \filldraw[purple1,line width=1pt, solid] (12,0) -- (12.6,0);
    \filldraw[purple1,line width=1pt] (12.6,0) node[right]{\footnotesize $p=6$};
\end{tikzpicture}
    \caption{Increase in critical time step size due to improved spectrum in explicit dynamics of a \textbf{square membrane and plate} with \textbf{fixed and simply supported boundary conditions}, respectively, as a function of the mesh $N = 2\nele \times 2\nele$ using $2 \times 2$ patches.}
    \label{fig:increase_time_step}
\end{figure}

For an undamped problem, the time step $\Delta t$ in the central difference method is bounded from above by the critical time step $\Delta t_{crit}$  \cite[Chapter~9,p.~492]{hughes_finite_2003}: 
\begin{align}
    \Delta t \leq \Delta t_{crit} = 2 / \omega^h_{\max} \,.
\end{align}

\begin{remark}
    The maximum eigenfrequency is obtained as part of Algorithm \ref{fig:parameter_estimation}. Hence, in our approach, we directly compute the critical time step size which is inversely proportional to the maximal frequency.
\end{remark}

Figure \ref{fig:time_step} plots the critical time-step size with respect to the number of elements per patch $\nele$, obtained with \changed{our perturbation approach  
(dashed curves)} and with standard analysis (solid curves), for the examples of a square membrane and plate with fixed and simply supported boundary conditions, respectively, as defined in the sections above. 
We consider $C^{p-1}$ B-splines of polynomial degrees $p=2$ through $6$ and multipatch discretizations of $2 \times 2$ patches with $C^0$ and $C^1$ patch continuity for the membrane and plate, respectively.
We observe that using \changed{our perturbation approach} 
allows for a significantly larger critical time step, and thus effectively reduces the associated computational cost of explicit dynamics calculations. \changed{We can also see that the our approach removes the dependency of the critical time-step size on the polynomial degree $p$, which exists for the standard analysis.} 
Figure \ref{fig:increase_time_step} shows the relative increase in the critical time-step size.

\section{Summary and conclusions}\label{sec:conclusion}

In this paper, we presented a variational approach based on perturbed eigenvalue analysis that reduces overestimated outlier frequencies due to reduced continuity at patch interfaces in isogeometric multipatch discretizations. It relies on the addition of scaled perturbation terms that weakly enforce $C^{p-1}$ continuity at patch interfaces. 
We also presented an iterative procedure to estimate effective scaling parameters for the perturbation term. It requires two input parameters (scaling factor $f>1$ and reduction factor $c \in (0,1)$) and computation of the maximum eigenfrequency and corresponding mode, which can be efficiently computed via power iteration. We demonstrated that our approach is robust with respect to the scaling factor $f > 1$, i.e.\ it reduces the outlier frequencies to approximately the same values for all $f > 1$. Furthermore, a reduction factor $c =0.9$ showed good results in all test cases.

We demonstrated numerically that the proposed approach improves spectral properties of multipatch discretizations for a variety bar, beam, membrane and 
plate. We  showed that the approach effectively addresses the outlier frequencies, while maintaining accuracy in the remainder of the spectrum and modes. We confirmed that spatial accuracy of the response was maintained in an explicit dynamics setting and showed that our approach allows for a much larger time-step size. In particular, we observed that the proposed approach removes the negative dependency of the critical time-step size on the polynomial degree $p$.

We note that our approach may be combined with the approach of \cite{deng_outlier_2021} to reduce outlier frequencies due to patch interfaces and boundaries. 
There are a number of avenues for future work. One aspect is the extension of our approach to non-uniform spline discretizations, trimmed and unfitted spline discretizations, and problems with non-smooth solution fields, where continuity constraints at patch interfaces cannot be consistently formulated. 
A second aspect is the further exploration of the case when $\alpha=0$ and $\beta<1$, which reduces the perturbation approach to the mass matrix, and to further study the resulting perturbation schemes in the context of mass lumping. One could investigate operator splitting techniques to move the added mass matrix to the right-hand side, which enables e.g. row-sum lumping of the unperturbed mass matrix. 
Another aspect is to study the performance of the proposed approach and the resulting problem conditioning in realistic scenarios with different material parameters.

\section*{Acknowledgments}
The authors gratefully acknowledge financial support from the German Research Foundation (Deutsche Forschungsgemeinschaft) through the DFG Emmy Noether Grant SCH 1249/2-1 and SCH 1249/5-1.

\bibliographystyle{elsarticle-num}
\bibliography{elsarticle-template}

\begin{thebibliography}{10}
\expandafter\ifx\csname url\endcsname\relax
  \def\url#1{\texttt{#1}}\fi
\expandafter\ifx\csname urlprefix\endcsname\relax\def\urlprefix{URL }\fi
\expandafter\ifx\csname href\endcsname\relax
  \def\href#1#2{#2} \def\path#1{#1}\fi

\bibitem{hughes_isogeometric_2005}
T.~J.~R. Hughes, J.~A. Cottrell, Y.~Bazilevs, Isogeometric analysis: {CAD},
  finite elements, {NURBS}, exact geometry and mesh refinement, Computer
  Methods in Applied Mechanics and Engineering 194~(39) (2005) 4135--4195.

\bibitem{cottrell_studies_2007}
J.~A. Cottrell, T.~J.~R. Hughes, A.~Reali, Studies of refinement and continuity
  in isogeometric structural analysis, Computer Methods in Applied Mechanics
  and Engineering 196~(41) (2007) 4160--4183.

\bibitem{hughes_finite_2014}
T.~J.~R. Hughes, J.~A. Evans, A.~Reali, Finite element and {NURBS}
  approximations of eigenvalue, boundary-value, and initial-value problems,
  Computer Methods in Applied Mechanics and Engineering 272 (2014) 290--320.

\bibitem{hughes_duality_2008}
T.~J.~R. Hughes, A.~Reali, G.~Sangalli, Duality and unified analysis of
  discrete approximations in structural dynamics and wave propagation:
  {Comparison} of p-method finite elements with k-method {NURBS}, Computer
  Methods in Applied Mechanics and Engineering 197~(49) (2008) 4104--4124.

\bibitem{puzyrev_spectral_2018}
V.~Puzyrev, Q.~Deng, V.~Calo, Spectral approximation properties of isogeometric
  analysis with variable continuity, Computer Methods in Applied Mechanics and
  Engineering 334 (2018) 22--39.

\bibitem{strang_analysis_2008}
G.~Strang, G.~Fix, An {Analysis} of the {Finite} {Element} {Method}, 2nd
  Edition, Wellesley-Cambridge Press, Wellesley, Mass, 2008.

\bibitem{hughes_finite_2003}
T.~J.~R. Hughes, The {Finite} {Element} {Method}: {Linear} {Static} and
  {Dynamic} {Finite} {Element} {Analysis}, Dover Publications, 2003.

\bibitem{cottrell_isogeometric_2006}
J.~A. Cottrell, A.~Reali, Y.~Bazilevs, T.~J.~R. Hughes, Isogeometric analysis
  of structural vibrations, Computer Methods in Applied Mechanics and
  Engineering 195~(41) (2006) 5257--5296.

\bibitem{cottrell_isogeometric_2009}
J.~A. Cottrell, T.~J.~R. Hughes, Y.~Bazilevs, Isogeometric {Analysis}: {Toward}
  {Integration} of {CAD} and {FEA}, Chichester, West Sussex, U.K. ; Hoboken,
  NJ, 2009.

\bibitem{hiemstra_outlier_2021}
R.~R. Hiemstra, T.~J.~R. Hughes, A.~Reali, D.~Schillinger, {Removal of spurious
  outlier frequencies and modes from isogeometric discretizations of second-
  and fourth-order problems in one, two, and three dimensions}, Computer
  Methods in Applied Mechanics and Engineering 387 (2021) 114115.

\bibitem{deng_outlier_2021}
Q.~Deng, V.~Calo, {A boundary penalization technique to remove outliers from
  isogeometric analysis on tensor-product meshes}, Computer Methods in Applied
  Mechanics and Engineering 383 (2021) 113907.

\bibitem{Manni2021}
C.~Manni, E.~Sande, H.~Speleers, {Application of optimal spline subspaces for
  the removal of spurious outliers in isogeometric discretizations} (2021).
\newblock \href {http://arxiv.org/abs/2106.03710} {\path{arXiv:2106.03710}}.

\bibitem{horger_penalty_2019}
T.~Horger, A.~Reali, B.~Wohlmuth, L.~Wunderlich, {A hybrid isogeometric
  approach on multi-patches with applications to Kirchhoff plates and
  eigenvalue problems}, Computer Methods in Applied Mechanics and Engineering
  348 (2019) 396--408.

\bibitem{Dittmann2019}
M.~Dittmann, S.~Schu{\ss}, B.~Wohlmuth, C.~Hesch, {Weak C n coupling for
  multipatch isogeometric analysis in solid mechanics}, International Journal
  for Numerical Methods in Engineering 118~(11) (2019) 678--699.

\bibitem{Schuss2019}
S.~Schu{\ss}, M.~Dittmann, B.~Wohlmuth, S.~Klinkel, C.~Hesch, {Multi-patch
  isogeometric analysis for Kirchhoff–Love shell elements}, Computer Methods
  in Applied Mechanics and Engineering 349 (2019) 91--116.

\bibitem{ainsworth_blended_quadrature_2010}
M.~Ainsworth, H.~A. Wajid, {Optimally Blended Spectral-Finite Element Scheme
  for Wave Propagation and NonStandard Reduced Integration}, SIAM Journal on
  Numerical Analysis 48~(1) (2010) 346--371.

\bibitem{puzyrev_quadrature_2017}
V.~Puzyrev, Q.~Deng, V.~Calo, {Dispersion-optimized quadrature rules for
  isogeometric analysis: Modified inner products, their dispersion properties,
  and optimally blended schemes}, Computer Methods in Applied Mechanics and
  Engineering 320 (2017) 421--443.

\bibitem{calo_quadrature_2019}
V.~Calo, Q.~Deng, V.~Puzyrev, {Dispersion optimized quadratures for
  isogeometric analysis}, Journal of Computational and Applied Mathematics 355
  (2019) 283--300.

\bibitem{deng_outlier2_2021}
Q.~Deng, V.~Calo, {Outlier Removal for Isogeometric Spectral Approximation with
  the Optimally-Blended Quadratures}~(017) (2021) 315--328.
\newblock \href {http://arxiv.org/abs/2102.07543} {\path{arXiv:2102.07543}}.

\bibitem{hartmann_mass_scaling_2015}
S.~Hartmann, D.~J. Benson, {Mass scaling and stable time step estimates for
  isogeometric analysis}, International Journal for Numerical Methods in
  Engineering 102~(3-4) (2015) 671--687.

\bibitem{Macek1995}
R.~W. Macek, B.~H. Aubert, {A mass penalty technique to control the critical
  time increment in explicit dynamic finite element analyses}, Earthquake
  Engineering \& Structural Dynamics 24~(10) (1995).

\bibitem{olovsson_mass_scaling_2005}
L.~Olovsson, K.~Simonsson, M.~Unosson, {Selective mass scaling for explicit
  finite element analyses}, International Journal for Numerical Methods in
  Engineering 63~(10) (2005) 1436--1445.

\bibitem{olovsson_mass_scaling_2006}
L.~Olovsson, K.~Simonsson, {Iterative solution technique in selective mass
  scaling}, Communications in Numerical Methods in Engineering 22~(1) (2006).

\bibitem{Askes2011}
H.~Askes, D.~C. Nguyen, A.~Tyas, {Increasing the critical time step:
  Micro-inertia, inertia penalties and mass scaling}, Computational Mechanics
  47~(6) (2011).

\bibitem{tkachuk_mass_scaling_2014}
A.~Tkachuk, M.~Bischoff, {Local and global strategies for optimal selective
  mass scaling}, Computational Mechanics 53~(6) (2014) 1197--1207.

\bibitem{schaeuble_mass_scaling_2017}
A.-K. Schaeuble, A.~Tkachuk, M.~Bischoff, {Variationally consistent inertia
  templates for B-spline- and NURBS-based FEM: Inertia scaling and
  customization}, Computer Methods in Applied Mechanics and Engineering 326
  (2017) 596--621.

\bibitem{cocchetti_mass_scaling_2013}
G.~Cocchetti, M.~Pagani, U.~Perego, {Selective mass scaling and critical
  time-step estimate for explicit dynamics analyses with solid-shell elements},
  Computers \& Structures 127 (2013) 39--52.

\bibitem{gonzalez_mass_tailoring_2020}
J.~A. Gonz{\'{a}}lez, K.~C. Park, {Large‐step explicit time integration via
  mass matrix tailoring}, International Journal for Numerical Methods in
  Engineering 121~(8) (2020) 1647--1664.

\bibitem{szilard_theories_2004}
R.~Szilard, Theories and {Applications} of {Plate} {Analysis}: {Classical},
  {Numerical} and {Engineering} {Methods}, Wiley, Hoboken, NJ, 2004.

\bibitem{ern_fem_2004}
A.~Ern, J.-L. Guermond, {Theory and Practice of Finite Elements},
  Springer-Verlag New York Inc., New York, USA, 2004.

\bibitem{tagliabue_error_2014}
A.~Tagliabue, L.~Ded{\`{e}}, A.~Quarteroni, {Isogeometric Analysis and error
  estimates for high order partial differential equations in fluid dynamics},
  Computers and Fluids 102 (2014) 277--303.

\bibitem{bazilevs_error_2006}
Y.~Bazilevs, L.~{Beir{\~{a}}o Da Veiga}, J.~A. Cottrell, T.~J.~R. Hughes,
  G.~Sangalli, {Isogeometric analysis: approximation, stability and error
  estimates for h-refined meshes}, Mathematical Models and Methods in Applied
  Sciences 16~(07) (2006) 1031--1090.

\bibitem{li_perturbation_2014}
R.-C. Li, {Matrix Perturbation Theory}, in: L.~Hogben (Ed.), Handbook of linear
  algebra, 1st Edition, Chapman \& Hall/CRC Press, LLC, 2007, Ch.~15.

\end{thebibliography}

\end{document}